\magnification 1200
  \input amssym
  \input miniltx
  \input pictex

  %
  \font \bbfive = bbm5
  \font \bbseven = bbm7
  \font \bbten = bbm10
  \font \eightbf = cmbx8
  \font \eighti = cmmi8 \skewchar \eighti = '177
  \font \eightit = cmti8
  \font \eightrm = cmr8
  \font \eightsl = cmsl8
  \font \eightsy = cmsy8 \skewchar \eightsy = '60
  \font \eighttt = cmtt8 \hyphenchar \eighttt = -1

  \font \sixi = cmmi6 \skewchar \sixi = '177
  \font \sixrm = cmr6
  \font \sixsy = cmsy6 \skewchar \sixsy = '60
  \font \tensc = cmcsc10

  \scriptfont \bffam = \bbseven
  \scriptscriptfont \bffam = \bbfive
  \textfont \bffam = \bbten

  \newskip \ttglue

  \def \eightpoint {\def \rm {\fam 0 \eightrm }\relax
  \textfont 0= \eightrm
  \scriptfont 0 = \sixrm \scriptscriptfont 0 = \fiverm
  \textfont 1 = \eighti
  \scriptfont 1 = \sixi \scriptscriptfont 1 = \fivei
  \textfont 2 = \eightsy
  \scriptfont 2 = \sixsy \scriptscriptfont 2 = \fivesy
  \textfont 3 = \tenex
  \scriptfont 3 = \tenex \scriptscriptfont 3 = \tenex
  \def \it {\fam \itfam \eightit }\relax
  \textfont \itfam = \eightit
  \def \sl {\fam \slfam \eightsl }\relax
  \textfont \slfam = \eightsl
  \def \bf {\fam \bffam \eightbf }\relax
  \textfont \bffam = \bbseven
  \scriptfont \bffam = \bbfive
  \scriptscriptfont \bffam = \bbfive
  \def \tt {\fam \ttfam \eighttt }\relax
  \textfont \ttfam = \eighttt
  \tt \ttglue = .5em plus.25em minus.15em
  \normalbaselineskip = 9pt
  \def \MF {{\manual opqr}\-{\manual stuq}}\relax
  \let \sc = \sixrm
  \let \big = \eightbig
  \setbox \strutbox = \hbox {\vrule height7pt depth2pt width0pt}\relax
  \normalbaselines \rm }

  \def \withfont #1#2{\font \auxfont =#1 {\auxfont #2}}

  %

  \def \ifundef #1{\expandafter \ifx \csname #1\endcsname \relax }

  \def \undefrule{\kern 2pt \vrule width 5pt height 5pt depth 0pt \kern 2pt}
  \def \possundef #1{\ifundef {#1}\undefrule {\eighttt #1}\undefrule \global \def \UndefFlag {}\else \csname #1\endcsname \fi }

  \def \TRUE {Y}
  \def \FALSE {N}

  %

  \newcount \secno \secno = 0
  \newcount \stno \stno = 0
  \newcount \eqcntr \eqcntr = 0

  \ifundef {showlabel} \global \def \showlabel {\FALSE} \fi  
  \ifundef {auxfile}   \global \def \auxfile   {\TRUE} \fi

  \def \define #1#2{\global \expandafter \edef \csname #1\endcsname {#2}}
  \def \error #1{\parindent 0pt \medskip \bf *******\hfil *******\hfil *******\hfil *******\hfil *******\hfil
    *******\hfil *******\hfil *******\break #1. \bigskip }

  \def \advseqnumbering {\global \advance \stno by 1 \global \eqcntr =0}

  \def \current {\ifnum \secno = 0 \number \stno \else \number \secno \ifnum \stno = 0 \else .\number \stno \fi \fi}

  \begingroup \catcode `\@=0 \catcode `\\=11 @global@def@textbackslash{\} @endgroup

  %
  \def \deflabel #1#2{%
    \if\TRUE\showlabel \hbox {\sixrm [[ #1 ]]} \fi
    \ifundef {#1PrimarilyDefined}%
      \define{#1}{#2}%
      \define{#1PrimarilyDefined}{#2}%
      \if\TRUE\auxfile \immediate \write 1 {\textbackslash newlabel {#1}{#2}}\fi
    \else
      \edef \old {\csname #1\endcsname}%
      \edef \new {#2}%
      \if \old \new \else \error{Duplicate definition for label ``{\tt #1}'', already defined as ``{\tt \csname #1\endcsname}''}\fi
      \fi}

  \def \label #1 {\deflabel {#1}{\current }}

  \def \equationmark #1 {\ifundef {InsideBlock}
	  \advseqnumbering
	  \eqno {(\current )}
	  \deflabel {#1}{\current }
	\else
	  \global \advance \eqcntr by 1
	  \edef \subeqmarkaux {\current .\number \eqcntr }
	  \eqno {(\subeqmarkaux )}
	  \deflabel {#1}{\subeqmarkaux }
	\fi }

  \def \split #1.#2.#3.#4;{\global \def \parone {#1}\global \def \partwo {#2}\global \def \parthree {#3}\global \def \parfour {#4}}
  \def \NA {NA}
  \def \ref #1{\split #1.NA.NA.NA;(\possundef {\parone }\ifx \partwo \NA \else .\partwo \fi )}
  \def \redundantref #1#2{\ref {#2}}

  %
  \newcount \bibno \bibno = 0

  \def \Bibitem #1 #2; #3; #4 \par{\smallbreak
    \global \advance \bibno by 1
    \item {[\possundef{#1}]} #2, {``#3''}, {#4}.\par
    \ifundef {#1PrimarilyDefined}\else
      \error{Duplicate definition for bibliography item ``{\tt #1}'', already defined in ``{\tt [\csname #1\endcsname]}''}
      \fi
	\ifundef {#1}\else
	  \edef \prevNum{\csname #1\endcsname}
	  \ifnum \bibno=\prevNum \else
		\error{Mismatch bibliography item ``{\tt #1}'', defined earlier as ``{\tt \prevNum}'' but should be ``{\tt \number\bibno}''}
		\fi
	  \fi
    \define{#1PrimarilyDefined}{#2}%
    \if\TRUE\auxfile \immediate\write 1 {\textbackslash newbib {#1}{\number\bibno}}\fi}

  \def \jrn #1, #2 (#3), #4-#5;{\sl #1, \bf #2 \rm (#3), #4--#5}
  \def \Article #1 #2; #3; #4 \par{\Bibitem #1 #2; #3; \jrn #4; \par}

  \def \references {\begingroup \bigbreak \eightpoint \centerline {\tensc References} \nobreak \medskip \frenchspacing }

  %

  \catcode `\@=11
  \def \c@itrk #1{{\bf \possundef {#1}}} 
  \def \c@ite #1{{\rm [\c@itrk{#1}]}}
  \def \sc@ite [#1]#2{{\rm [\c@itrk{#2}\hskip 0.7pt:\hskip 2pt #1]}}
  \def \du@lcite {\if \pe@k [\expandafter \sc@ite \else \expandafter \c@ite \fi }
  \def \cite {\futurelet\pe@k \du@lcite }
  \catcode `\@=12

  %
  \def \Headlines #1#2{\nopagenumbers
    \headline {\ifnum \pageno = 1 \hfil
    \else \ifodd \pageno \tensc \hfil \lcase {#1} \hfil \folio
    \else \tensc \folio \hfil \lcase {#2} \hfil
    \fi \fi }}

  \def \title #1{\medskip\centerline {\withfont {cmbx12}{\ucase{#1}}}}

  \def \Subjclass #1#2{\footnote {\null }{\eightrm #1 \eightsl Mathematics Subject Classification:  \eightrm #2.}}

  \long \def \Quote #1\endQuote {\begingroup \leftskip 35pt \rightskip 35pt
\parindent 17pt \eightpoint #1\par \endgroup }
  \long \def \Abstract #1\endAbstract {\bigskip \Quote \noindent #1\endQuote }
  
  \def \Authors #1{\bigskip \centerline {\tensc #1}}
  \def \Note #1{\footnote {}{\eightpoint #1}}
  \def \Date #1 {\Note {\it Date: #1.}}

  \def \part #1#2{\vfill \eject \null \vskip 0.3\vsize
    \withfont{cmbx10 scaled 1440}{\centerline{PART #1} \vskip 1.5cm \centerline{#2}} \vfill\eject }

  %

  \def \fix {\smallskip \noindent $\blacktriangleright $\kern 12pt}
  \def \iskip {\medskip\noindent}

  \def \ucase #1{\edef \auxvar {\uppercase {#1}}\auxvar }
  \def \lcase #1{\edef \auxvar {\lowercase {#1}}\auxvar }

  \def \section #1 \par {\global \advance \secno by 1 \stno = 0
    \goodbreak \bigbreak
    \noindent {\bf \number \secno .\enspace #1.}
    \nobreak \medskip \noindent }

  \def \state #1 #2\par {\begingroup \def \InsideBlock {} \medbreak \noindent \advseqnumbering {\bf \current .\enspace
#1.\enspace \sl #2\par }\medbreak \endgroup }

  \def \definition #1\par {\state Definition \rm #1\par }

  \long \def \Proof #1\endProof {\begingroup \def \InsideBlock {} \medbreak \noindent {\it Proof.\enspace }#1
\ifmmode \eqno \endproofmarker $$ \else \hfill $\endproofmarker $ \looseness = -1 \fi \medbreak \endgroup }

  \def \$#1{#1 $$$$ #1}
  \def \explica #1#2{\mathrel {\buildrel \hbox {\sixrm #2} \over #1}}
  \def \explain #1#2{\explica{#1}{\ref{#2}}}  
  \def \=#1{\explain {=}{#1}}

  \def \pilar #1{\vrule height #1 width 0pt}
  \def \stake #1{\vrule depth  #1 width 0pt}

  \newcount \fnctr \fnctr = 0
  \def \fn #1{\global \advance \fnctr by 1
    \edef \footnumb {$^{\number \fnctr }$}%
    \footnote {\footnumb }{\eightpoint #1\par \vskip -10pt}}

  \def \text #1{\hbox {#1}}
  \def \bool #1{[{\scriptstyle #1}]\,}
  \def \equal #1#2{\bool {#1=#2}}

  %
  
  \def \Item #1{\smallskip \item {{\rm #1}}}
  \newcount \zitemno \zitemno = 0

  \def \izitem {\global \zitemno = 0}
  \def \zitemplus {\global \advance \zitemno by 1 \relax }
  \def \rzitem {\romannumeral \zitemno }
  \def \rzitemplus {\zitemplus \rzitem }
  \def \zitem {\Item {{\rm (\rzitemplus )}}}
  \def \Zitem {\izitem \zitem }
  \def \zitemmark #1 {\deflabel {#1}{\rzitem }}

  \newcount \nitemno \nitemno = 0
  \def \initem {\nitemno = 0}
  \def \nitem {\global \advance \nitemno by 1 \Item {{\rm (\number \nitemno )}}}

  \newcount \aitemno \aitemno = -1
  \def \boxlet #1{\hbox to 6.5pt{\hfill #1\hfill }}
  \def \iaitem {\aitemno = -1}
  \def \aitemconv {\ifcase \aitemno a\or b\or c\or d\or e\or f\or g\or
h\or i\or j\or k\or l\or m\or n\or o\or p\or q\or r\or s\or t\or u\or
v\or w\or x\or y\or z\else zzz\fi }
  \def \aitem {\global \advance \aitemno by 1\Item {(\boxlet \aitemconv )}}
  \def \aitemmark #1 {\deflabel {#1}{\aitemconv }}

  \def \Case #1:{\medskip \noindent {\tensc Case #1:}}

  %
  \def \<{\left \langle\vrule width 0pt depth 0pt height 8pt }
  \def \>{\right \rangle}
  \def \({\big (}
  \def \){\big )}
  \def \ds {\displaystyle }
  \def \and {\hbox {,\quad and \quad }}

  \def \IFF {\kern 7pt\Leftrightarrow \kern 7pt}
  \def \IMPLY {\kern 7pt \Rightarrow\kern 7pt}
  \def \for #1{,\quad \forall\,#1}
  \def \endproofmarker {\square }
  \def \"#1{{\it #1}\/} \def \umlaut #1{{\accent "7F #1}}
  \def \inv {^{-1}}
  \def \*{\otimes}
  \def \caldef #1{\global \expandafter \edef \csname #1\endcsname {{\cal #1}}}
  \def \bfdef #1{\global \expandafter \edef \csname #1\endcsname {{\bf #1}}}
  \bfdef N \bfdef Z \bfdef C \bfdef R

  %

  \if\TRUE\auxfile
    \IfFileExists {\jobname.bib}{\input \jobname.bib }{\null}
    \IfFileExists {\jobname.aux}{\input \jobname.aux }{\null}
    \immediate \openout 1 \jobname.aux
    \fi

  \def\close{\if\TRUE\auxfile \closeout 1 \fi
    \ifundef {UndefFlag}\else
      \message{*** There were undefined labels ***} \iskip \tt ******* There were undefined labels *******
      \fi
    \par\vfill\supereject\end}

  %

  \def \Caixa #1{\setbox 1=\hbox {$#1$\kern 1pt}\global \edef \tamcaixa {\the \wd 1}\box 1}
  \def \caixa #1{\hbox to \tamcaixa {$#1$\hfil }}

  \def \med #1{\mathop {\textstyle #1}\limits }

  \def \medcup {\med \bigcup}

  %

  \def \src {d}	 \def\sr#1{\src(#1)}               
  \def \ran {r}	 \def\rn#1{\ran(#1)}               
  \def \vr {x}	                                   
  \def \vro {y}	                                   
  \def \ed {e}                                     
  \def \oed {f}	                                   
  \def \eproj {f}                                  
  \def \s {s}  	                                   
  \def \auto {\sigma}                                   
  \def \ts {\tilde \s }
  \def \tp {\tilde p}
  \def \tu {\tilde u}
  \def \E {{\cal E}}                               
  \def \SGE {{\cal S}_{G,E}}                       
  \def \EGE {{\cal E}}                             
  \def \SE {{\cal S}_E}                            
  \def \GpdGE {{\cal G}\tight (\SGE )}             

  \def \q {\breve }
  \def \corona {{\q G}}
  \def \lag {\ell}
  \def \O {{\cal O}}

  \def \g {{\bf g}}  
  \def \cyl #1{Z(#1)}

  \def \germ #1#2#3#4{\big [#1,#2,#3;\,#4\big ]}
  \def \trunc #1#2{#1|_{#2}}

  \def \Lin {{\cal L}}
  \def \acite [#1]{\cite [#1]{actions}}



  \def \S {{\cal S}}
  \def \E {{\cal E}}

  \def \J {{\cal J}}
  \def \tight {_{\rm tight}}


  \def \G {{\cal G}}




  \def \clos #1{\overline {#1}}



  \def \Data {G,E,\varphi}
  \def \OGE {{\cal O}_{G,E}}
  \def \OAB {{\cal O}_{A,B}}


  \def \proj {q} 
  \def \gp {v} 
  \def \rep {V} 
  \def \t {t} 
  \def \modmap {\Psi} 
  \def \repalg {\psi} 
  \def \projMod {Q} 


  \def \mathcal #1{{\cal #1}}
  \def \text #1{\hbox {#1}}
  \def \mathbb #1{{\bf #1}}
  \def \mbox #1{\hbox {#1}}
  \def \frac #1#2{{#1 \over #2}}

  \def \vspace #1{\vskip #1}
  \def \hspace #1{\hskip #1}

  \newcount \itemdpt \itemdpt = 0
  \def \uplevel {\begingroup \global \advance \itemdpt by 1 \advance \parindent by 18pt
  \ifcase \itemdpt \or \initem \or \iaitem \or \izitem \fi }
  \def \dnlevel {\par \global \advance \itemdpt by -1 \advance \parindent by -18pt \endgroup }
  \def \itm {\ifcase \itemdpt \or \nitem \or \aitem \or \zitem \fi }
  \def \notext {\vskip -15pt}

  \def \N {{\mathbb {N}}}
  \def \Ninf {{\mathbb {N}\cup\{\infty\}}}
  \def \Z {{\mathbb {Z}}}
  
  \def \OAB {{\mathcal {O}_{A,B}}}
  \def \OGE {{\mathcal {O}_{G,E}}}
  \def \Zplus {\Z ^+}
  \def \OmA {{\Omega_{A}}}
  \def \SZE {{\mathcal {S}_{\Z , E}}}
  \def \GZE {{\cal G}\tight ({\mathcal {S}_{\Z , E}})}


  \def \quoapprox {E^0{\kern -1pt/\kern -2pt}\approx }

\title{Self-Similar Graphs}
\title{A unified treatment of Katsura and}
\title{Nekrashevych C*-algebras}

  \Headlines {Inverse semigroup actions and  self-similar graphs} {R.~Exel and E.~Pardo}

  \Authors {Ruy Exel and Enrique Pardo}

  \Date {27 September 2014}

  \Subjclass {2010}{46L05, 46L55}

  \Note {\it Key words and phrases: \rm Kirchberg algebra, Katsura algebra, Nekrashevych algebra, tight representation,
inverse semigroup, groupoid, groupoid C*-algebra.}

  \Note {The first-named author was partially supported by CNPq. The second-named author was partially supported by PAI
III grants FQM-298 and P11-FQM-7156 of the Junta de Andaluc\'{\i}a and by the DGI-MICINN and European Regional
Development Fund, jointly, through Project MTM2011-28992-C02-02.}

  \Abstract Given a graph $E$, an action of a group $G$ on $E$, and a $G$-valued cocycle $\varphi$
on the  edges of $E$, we define a C*-algebra denoted $\OGE$, which is shown to be isomorphic to the tight C*-algebra
associated to a certain inverse semigroup $\SGE$ built naturally from the triple $(\Data)$.  As a tight C*-algebra, $\OGE$ is
also isomorphic to the full  C*-algebra of a naturally occurring groupoid  $\GpdGE$.  We then study the relationship
between properties of the action, of the groupoid and of the C*-algebra, with an emphasis on situations in which $\OGE$ is a
Kirchberg algebra.  Our main applications are to   Katsura algebras
  and to certain algebras constructed by Nekrashevych  from self-similar groups.  These two classes of C*-algebras are shown to be special
cases of our $\OGE$, and many of their known properties are shown to follow from our general theory.
  \endAbstract

\section Introduction

The purpose of this paper is to give a unified treatment to two classes of C*-algebras which have been studied in the
past few years from rather different points of view, namely
  Katsura algebras \cite {KatsuraOne},
  and certain algebras constructed by Nekrashevych \cite {NekraJO}, \cite {NC} from self-similar groups.

The realization that these classes are indeed closely related, as well as the fact that they could be given a unified
treatment, came to our mind as a result of our attempt to understand Katsura's algebras $\O _{A,B}$ from the point of
view of inverse semigroups. The fact, proven by Katsura in \cite {KatsuraOne}, that all Kirchberg algebras in the UCT
class may be described in terms of his $\O _{A,B}$ was, in turn, a strong motivation for that endeavor.

While studying $\O _{A,B}$, it slowly became clear to us that the two matricial parameters $A$ and $B$, present in
Katsura's construction, play very
different roles.  The reader acquainted with Katsura's work will easily recognize that the matrix $A$ is destined to be
viewed as the edge matrix of a graph, but it took us much longer to realize that $B$ should be thought of as providing
parameters for an action of the group ${\bf Z}$ on the graph given by $A$.  In trying to understand these different roles,
some interesting arithmetic popped up sparking a connection with the work done by Nekrashevych \cite {NC} on the
C*-algebra $\O _{(G,X)}$ associated to a self-similar group $(G,X)$.

While Nekrashevych's algebras contain a Cuntz algebra, Katsura's algebras contain a graph C*-algebra.  This fact alone
ought to be considered as a hint that self-similar groups lie in a much bigger class, where the group action takes place
on the path space of a graph, rather than on a rooted tree (which, incidentally, is the path space of a bouquet of
circles).

One of the first important applications of the idea of self-similarity in group theory is in constructing groups with
exotic properties \cite {Grig}, \cite {GS}.  Many of these are defined as subgroups of the group of all automorphisms of
a tree.  Having been born from automorphisms, it is natural that the theory of self-similar groups generally assumes
that the group acts \"{faithfully} on its tree (see, e.g. \cite [Definition 2.1]{NC}).

However, based on the example provided by Katsura's algebras, we decided that perhaps it is best to view the group on
its own, the action being an extra ingredient.

The main idea behind self-similar groups, namely the equation
  $$
  g(xw) = yh(w)
  \equationmark SSimilerity
  $$
  appearing in \cite [Definition 2.1]{NC}, and the subsequent notion of \"{restriction}, namely
  $$
  g|_x:= h,
  $$
  depend on faithfulness, since otherwise the group element $h$ appearing in \ref{SSimilerity} would not be unique
and therefore will not be well defined as a function of $g$ and $x$.  Working with non-faithful group actions we were
forced to postulate a functional dependence
  $$
  h = \varphi(g,x),
  $$
  and we were surprised to find that the natural properties expected of $\varphi$ are that of a group cocycle.

To be precise, the ingredients needed in our generalization of self-similar groups are: a countable discrete group $G$,
an action
  $$
  G\times E\to E
  $$
  of $G$ on a finite graph
  $E = (E^0, E^1, \ran , \src )$, and
  a one-cocycle
  $$
  \varphi:G\times E^1 \to G
  $$
  for the action of $G$ on the edges of $E$.

Starting with this data (which we assume satisfies a few other natural axioms) we construct an action of $G$ on the
space of finite paths $E^*$ which satisfies the ``self-similarity'' equation
  $$
  g(\alpha\beta) = (g\alpha)\big(\varphi(g,\alpha)\beta\big)
  \for g\in E \for \alpha,\beta\in E^*.
  $$

Adopting a philosophy similar to that embraced by Katsura and Nekrashevych, we define a C*-algebra, denoted
  $$
  \OGE ,
  $$
  in terms of generators and relations inspired by the above group action.  The study of $\OGE $ is, thus, the purpose of
this paper.

Given a self-similar group $(G,X)$, if we consider $X$ as the set of edges of a graph with a single vertex, and if we
define $\varphi(g,x) = g|_x$, then our $\OGE $ coincides with Nekrashevych's $\O _{(G,X)}$.

On the other hand, if we are given two integer $N\times N$ matrices $A$ and $B$, with $A_{i,j}\geq0$, for all $i$ and $j$,
we may form a graph $E$ with vertex set $E^0=\{1,2,\ldots,N\}$ and with $A_{i,j}$ edges from vertex $i$ to vertex $j$.

We may then use $B$ to define an action of ${\bf Z}$ on $E$, by fixing all vertices and acting on the set of edges as
follows:
  denote the set of edges in $E$ from $i$ to $j$ by
  $$
  \{e_{i,j,n}: 0\leq n<A_{i,j}\}.
  $$
  Given $m\in{\bf Z}$, and given an edge $e_{i,j,n}$, in order to define $\auto _m(e_{i,j,n})$, we first perform the
Euclidean division of $mB_{i,j}+n$ by $A_{i,j}$, say
  $$
  mB_{i,j}+n=\hat k A_{i,j} + \hat {n}
  $$
  with $0\leq\hat {n}<A_{i,j}$.  We then put
  $$
  \auto _m(e_{i,j,n}) := e_{i,j,\hat n},
  $$
  so that the group element $m$ permutes the $A_{i,j}$ edges from $i$ to $j$ in the same way that addition by
$mB_{i,j}$, modulo $A_{i,j}$, permutes the integers $\{0,1,\ldots,A_{i,j}-1\}$.

The quotient $\hat k$ in the above Euclidean division also plays an important role, namely in the definition of the
cocycle:
  $$
  \varphi(m, e_{i,j,n}):= \hat k.
  $$

In possession of the graph, the action of ${\bf Z}$, and the cocycle $\varphi$ constructed above, we apply our
construction and we find that $\OGE $ is isomorphic to Katsura's $\OAB $.

So, both Nekrashevych's and Katsura's algebras become special cases of our construction.  We therefore believe that the
project of studying such group actions on path spaces as well as the corresponding algebras is of great importance.

Taking the first few steps in this direction we have been able to describe $\OGE $ as the C*-algebra of an \'etale
groupoid $\G _{G,E} $, whose construction is remarkably similar to the groupoid associated to the relation of ``tail
equivalence with lag'' on the path space, as described by Kumjian, Pask, Raeburn and Renault in \cite {KPRR}.

The first similarity is that our groupoid $\G _{G,E} $ has the exact same unit space as the corresponding graph groupoid,
namely the infinite path space.  The second, and most surprising similarity is that $\G _{G,E} $ is also described by a
\"{lag} function, except that the values of the lag are not integer numbers, as in \cite {KPRR}, but lie in a slightly
more complicated group, namely the semi-direct product of the \"{corona} group of $G$ by the right shift automorphism
(see below for precise definitions).

We would like to stress that, like Nekrashevych's groupoid \cite [Theorem 5.1]{NC}, our groupoid $\G _{G,E} $ is constructed
as a groupoid of germs.  However, departing from Nekrashevych's techniques, we use Patterson's \cite {pat} notion of
``germs'', rather than the one employed in \cite [Section 5]{NC}.  While agreeing in many cases, such as when the action
is topologically free (see below for the precise definition), the former has a much better chance of producing Hausdorff
groupoids and, in our case, we are able to give a precise characterization of Hausdorffness in terms of a property we
call \"{pseudo freeness} (see below for the precise definition).

The techniques we use to give $\OGE $ a groupoid model bear heavily on the theory of tight representations of inverse
semigroups developed by the first named author in \cite {actions}.  In particular, from our initial data we construct an
abstract inverse semigroup $\SGE $ and show that $\OGE $ is the universal C*-algebra for tight representations of $\SGE $.

In another direction we again take inspiration from Nekrashevych \cite {NekraJO} and give a description of $\OGE $ as a
Cuntz-Pimsner algebra for a very natural correspondence $M$ over the algebra
  $$
  C(E^0) \ifundef {rtimes} \times \else \rtimes \fi G.
  $$
  As a result we are able to prove that $\OGE $ is nuclear when $G$ is amenable.

We briefly study the natural representation of the graph C*-algebra $C^*(E)$ \cite {Raeburn} into $\OGE $, which turns out to be faithful. Also, we study the natural representation of the group $G$ into $\OGE $, which turns out to be faithful when the triple $(G,E, \varphi)$ satisfies pseudo freeness, but fails in general.

Simplicity of $\OGE $ is also discussed by using our description of this algebra as a groupoid C*-algebra and employing
results from \cite {SimpleGroupoid}.  In doing so, it is crucial to determine when is $\GpdGE $ a Hausdorff, minimal essentially principal groupoid. To this end, we strongly rely on results obtained by both authors in \cite{EPFour} about characterization of minimality and essential irreducibility for the groupoid of germs of a general $\ast$-inverse semigroup. We then specialize these results to the particular context of the inverse semigroup $\SGE$. Hence, we characterize Hausdorffness of $\GpdGE $ in terms of the existence of finitely many minimal strongly fixed paths (see below for a precise definition).

Also, we characterize minimality of $\GpdGE $ in terms of weak $G$-transitivity of the graph (see Section 13 for a definition of this concept).  We then
obtain a natural generalization of the analog result obtained in \cite {ExelLaca} for Exel-Laca algebras.

We also show that being essentially principal is related to the topological freeness of the
action of $\SGE $ on the infinite path space.  In this sense, we obtain a characterization that relies on the existence of entries for any circuit of the graph, plus a formal condition which forces any element of $G$ fixing open sets to be tighly related to the existence of suitable minimal strongly fixed paths.

Moreover, we give sufficient conditions on $\GpdGE $ to guarantee its
local contractiveness  (see e.g. \cite {AdelaR} for a definition); this property turns out to be a consequence of essential principality, so that any simple algebra in the class $\OGE$ will be purely infinite simple.

With the machinery developped we are then able to give a characterization of simplicity (and so pure infinite simplicity) for $\OGE $ when $\GpdGE $
is Hausdorff.

Finally, we revisit the case of Nekrashevych and Katsura algebras, giving a picture of the properties enjoyed by these
algebras that turns out to be more general than the ones given by Nekrashevych or Katsura.

Some of the results in the present paper appeared in the preprints \cite {EP} and \cite {EPTwo}, which in turn are to be
replaced by the present work.

We would also like to mention \cite{ExelStar} and \cite{Starling}, which are strongly related to the algebras we
study here.  In \cite{ExelStar} conditions are given for $\OGE$ to be a partial crossed product and in \cite{Starling}
an interesting connection with Zappa-Sz\'ep products is made.

Part of this work was done during visits of the second named author to the Departamento de Matem\'atica da Universidade
Federal de Santa Catarina (Florian\'opolis, Brasil) and he would like to express his thanks to the host center for its
warm hospitality. Both authors thank Benjamin Steinberg for interesting discussions on topological freeness of actions,
and Hausdorffness of groupoids.

\section Groups acting on graphs

Let $E = (E^0, E^1, \ran , \src )$ be a directed graph, where $E^0$ denotes the set of \"{vertices}, $E^1$ is the set of
\"{edges}, $\ran $ is the \"{range} map, and $\src $ is the \"{source}, or \"{domain} map.

By definition, a \"{source} in $E$ is a vertex $\vr \in E^0 $, for which $\ran \inv (\vr )=\ifundef {varnothing} \emptyset
\else \varnothing \fi $.  Thus, when we say that a graph has \"{no sources}, we mean that $\ran \inv (\vr )\neq \ifundef
{varnothing} \emptyset \else \varnothing \fi $, for all $\vr \in E^0$.

  By an \"{automorphism} of $E$ we shall mean a bijective map
  $$
  \auto : E^0 \mathop {\dot \cup } E^1 \to E^0 \mathop {\dot \cup } E^1
  $$
  such that
  $\auto (E^i)\subseteq E^i$, for $i = 0,1$,
  and moreover such that
   $\ran \circ \auto = \auto \circ \ran $, and $\src \circ \auto = \auto \circ \src $, on $E^1$.
  It is evident that the collection of all automorphisms of $E$ forms a group under composition.

  By an action of a group $G$ on a graph $E$ we shall mean a group homomorphism from $G$ to the group of all
automorphisms of $E$.

If $X$ is any set, and if $\auto $ is an action of a group $G$ on $X$, we shall say that a map
  $$
  \varphi :G\times X \to G
  $$
  is a \"{one-cocycle} for $\auto $, when
  $$
  \varphi (gh, x) =
  \varphi \big(g,\auto _h(x)\big)\varphi (h,x),
  \equationmark CocycleId
  $$
  for all $g,h \in G$, and all $x \in X$.  In this case, plugging $g = h = 1$, above, we see that necessarily
  $$
  \varphi (1,x) = 1,
  \equationmark CocycleAtOne
  $$
  for every $x$.

\state {$\blacktriangleright $ Standing Hypothesis} \label StandingHyp
  \rm Throughout this work we shall let
    $G$ be a countable discrete group,
    $E$ be a finite graph with no sources,
    $\auto $ be an action of $G$ on $E$,
    and
  $$
  \varphi :G\times E^1 \to G
  $$
  be a one-cocycle for the restriction of $\auto $ to $E^1$, which moreover satisfies
  $$
  \auto _{\varphi (g,\ed )}(\vr ) = \auto _g(\vr )
  \for g \in G \for \ed \in E^1 \for \vr \in E^0.
  \equationmark ActionOfCocycleOnVertex
  $$

The assumptions that $E$ is finite and has no sources will in fact only be used in the next section and it could
probably be removed by using well known graph C*-algebra techniques.

  By a \"{path} in $E$ of \"{length} $n\geq 1$ we shall mean any finite sequence of the form
  $$
  \alpha = \alpha _1\alpha _2\ldots \alpha _n,
  $$
  where $\alpha _i \in E^1$, and $\src (\alpha _i) = \ran (\alpha _{i+1})$, for all $i$ (this is the usual convention when treating graphs
from a categorical point of view, in which functions compose from right to left).  The \"{range} of $\alpha $ is defined by
  $$
  \ran (\alpha ) = \ran (\alpha _1),
  $$
  while the \"{source} of $\alpha $ is defined by
  $$
  \src (\alpha ) = \src (\alpha _n).
  $$

  A vertex $\vr \in E^0$ will be considered a path of length zero, in which case we set $\ran (\vr ) = \src (\vr ) = \vr $.

  For every integer $n\geq 0$ we denote by $E^n$ the set of all paths in $E$ of length $n$ (this being consistent with the
already introduced notations for $E^0$ and $E^1$).  Finally, we denote by $E^*$ the sets of all finite paths, and by
$E^{\leq n}$ the set of all paths of length at most $n$, namely
  $$
  E^* = \medcup _{k\geq 0}E^k
  \and
  E^{\leq n} = \medcup _{k=0}^nE^k.
  $$

We will often employ the operation of \"{concatenation} of paths.  That is, if (and only if) $\alpha $ and $\beta $ are paths such
that $\src (\alpha )=\ran (\beta )$, we will denote by $\alpha \beta $ the path obtained by juxtaposing $\alpha $ and $\beta $.

In the special case in which $\alpha $ is a path of length zero, the concatenation $\alpha \beta $ is allowed if and only if $\alpha =\ran
(\beta )$, in which case we set $\alpha \beta =\beta $.  Similarly, when $|\beta |=0$, then $\alpha  \beta $ is defined iff $\src (\alpha )=\beta $, and then $\alpha \beta =\alpha $.

We would now like to describe a certain extension of $\auto $ and $\varphi $ to finite paths.

\state Proposition \label extendedaction
  Under the assumptions of \ref{StandingHyp}
  there exists a unique pair $(\auto ^*,\varphi ^*)$, formed by an action $\auto ^*$ of $G$ on $E^*$ (viewed simply as a set),
  and a one-cocycle $\varphi ^*$ for $\auto ^*$, such that, for every $n\geq 0$, every $g \in G$, and every $\vr \in E^0$, one has that:
  \Zitem $\auto ^*_g=\auto _g$, on $E^{\leq 1}$,
  \zitem $\varphi ^*(g,\vr ) = g$, \zitemmark CocZero
  \zitem $\varphi ^* = \varphi $, on $G\times E^1$, \zitemmark ExtPhi
  \zitem $\auto ^*_g(E^n)\subseteq E^n$, \zitemmark KeepLength
  \zitem $\ran \circ \auto ^*_g=\auto _g\circ \ran $, on $E^n$, \zitemmark MatchRange
  \zitem $\src \circ \auto ^*_g=\auto _g\circ \src $, on $E^n$, \zitemmark MatchSource
  \zitem $\auto _{\varphi ^*(g,\alpha )}(\vr )=\auto _g(\vr )$, for all $\alpha \in E^n$, \zitemmark PhiStarOnVert
  \zitem $\auto ^*_1$ is the
  identity\fn {This is evidently already included in the statement that $\auto ^*$ is an action, but we repeat it here to
aid our proof by induction.}
  on $E^n$, \zitemmark PhiOneId
  \zitem $\auto ^*_g(\alpha \beta ) = \auto ^*_g(\alpha )\ \auto ^*_{\varphi ^*(g,\alpha )}(\beta )$, provided $\alpha $ and $\beta $ are finite paths with $\alpha \beta \in E^n$,
\zitemmark MainConcat
  \zitem $\varphi ^*(g, \alpha \beta )=\varphi ^*\big(\varphi ^*(g,\alpha ),\beta \big)$, provided $\alpha $ and $\beta $ are finite paths with $\alpha \beta \in E^n$. \zitemmark LastItem

\Proof Initially notice that, once \ref{MatchRange}, \ref{MatchSource} and \ref{PhiStarOnVert} are
proved, the concatenation of the paths
  ``$\auto ^*_g(\alpha )$'' and ``$\auto ^*_{\varphi ^*(g,\alpha )}(\beta )$'',
  appearing in \ref{MainConcat}, is permitted because
  $$
  \ran \big(\auto ^*_{\varphi ^*(g,\alpha )}(\beta )\big) \={MatchRange}
  \auto _{\varphi ^*(g,\alpha )}(\ran (\beta )) \={PhiStarOnVert}
  \auto _g(\ran (\beta )) =
  \auto _g\big(\src (\alpha )\big) \={MatchSource}
  \src \big(\auto ^*_g(\alpha )\big).
  $$

For every $g$ in $G$, define $\auto ^*_g$ on $E^{\leq 1}$ to coincide with $\auto _g$.  Also, define $\varphi ^*$ on $G\times E^{\leq 1}$ by
\ref{CocZero} and \ref{ExtPhi}.  It is then clear that (i--iii) hold and it is easy to see that the remaining
properties (iv--\LastItem ) hold for all $n\leq 1$.

We shall complete the definitions of $\auto ^*$ and $\varphi ^*$ by induction, so we assume that $m\geq 1$, that
  $$
  \auto ^*_g: E^{\leq m} \to E^{\leq m}
  $$
  is defined for all $g$ in $G$, that
  $$
  \varphi ^*: G\times E^{\leq m} \to G,
  $$
  is defined,
  and that (i--\LastItem ) hold for all $n\leq m$.  We then define
  $$
  \auto ^*_g: E^{m+1} \to E^{m+1}
  $$
  for all $g$ in $G$, and
  $$
  \varphi ^*: G\times E^{m+1} \to G,
  $$
  by induction as follows.
  Given $\alpha \in E^{m+1}$, write $\alpha =\alpha ' \alpha ''$, with $\alpha '\in E^1$, and $\alpha ''\in E^m$, and put
  $$
  \auto ^*_g(\alpha ) = \auto _g(\alpha ')\auto ^*_{\varphi (g,\alpha ')}(\alpha '')
  \and
  \varphi ^*(g, \alpha )=\varphi ^*\big(\varphi (g,\alpha '),\alpha ''\big).
  \equationmark InducDef
  $$
  A quick analysis, as done in the first paragraph of this proof, shows that the concatenation of ``$\auto _g(\alpha ')$''
and ``$\auto ^*_{\varphi (g,\alpha ')}(\alpha '')\stake {7pt}$'', appearing above, is permitted.  We next verify (iv--\LastItem ),
substituting $m+1$ for $n$.

We have that the length of $\auto ^*_g(\alpha )$, as defined above, is clearly $1+m$, thus proving (iv).  With respect to
\ref{MatchRange} we have that
  $$
  \ran \big(\auto ^*_g(\alpha )\big) =\ran \big(\auto _g(\alpha ')\big) = \auto _g\big(\ran (\alpha ')\big) = \auto _g\big(\ran (\alpha )\big).
  $$
  As for \ref{MatchSource}, notice that
  $$
  \src \big(\auto ^*_g(\alpha )\big) = \src \big(\auto ^*_{\varphi (g,\alpha ')}(\alpha '')\big) = \auto _{\varphi (g,\alpha ')}\big(\src (\alpha '')\big) = \auto _g\big(\src (\alpha '')\big)
= \auto _g\big(\src (\alpha )\big).
  $$
  Given $\vr \in E^0$, we have that
  $$
  \auto _{\varphi ^*(g,\alpha )}(\vr ) =
  \auto _{\varphi ^*(\varphi (g,\alpha '),\alpha '')}(\vr ) =
  \auto _{\varphi (g,\alpha ')}(\vr ) =
  \auto _g(\vr ),
  $$
  taking care of \ref{PhiStarOnVert}.

The verification of \ref{PhiOneId} is done as follows: for $\alpha =\alpha '\alpha ''$, as in \ref{InducDef}, one has
  $$
  \auto ^*_1(\alpha ) = \auto ^*_1(\alpha '\alpha '') =
  \auto _1(\alpha ')\auto ^*_{\varphi (1,\alpha ')}(\alpha '') \={CocycleAtOne}
  \auto _1(\alpha ')\auto ^*_1(\alpha '') = \alpha '\alpha ''=\alpha .
  $$

In order to prove \ref{MainConcat}, pick paths $\alpha $ in $E^k$ and $\beta $ in $E^l$, where $k+l=m+1$, and such that $\src
(\alpha )=\ran (\beta )$.

We leave it for the reader to verify \ref{MainConcat} in the easy case in which $k=0$, that is, when $\alpha $ is a
vertex.  The case $k=1$ is also easy as it is nothing but the definition of $\auto ^*_g$ given in \ref{InducDef}.
So we may assume that $k\geq 2$.

  Writing $\alpha =\alpha '\alpha ''$, with $\alpha '\in E^1$, and $\alpha ''\in E^{k-1}$, we then have that $\alpha \beta = \alpha '\alpha ''\beta $, and hence, by
definition,
  $$
  \auto ^*_g(\alpha \beta ) = \auto _g(\alpha ')\auto ^*_{\varphi (g,\alpha ')}(\alpha ''\beta ) =
  \auto _g(\alpha ') \auto ^*_{\varphi (g,\alpha ')}(\alpha '')\ \auto ^*_{\varphi ^*(\varphi (g,\alpha '),\alpha '')}(\beta ) \$=
  \auto _g^*(\alpha '\alpha '')\ \auto ^*_{\varphi ^*(g,\alpha '\alpha '')}(\beta ).
  $$
  We remark that, in last step above, one should use the induction hypothesis in case $k\leq m$, and the definitions of
$\auto ^*$ and $\varphi ^*$, when $k=m+1$.

To verify \ref{LastItem} we again pick paths $\alpha $ in $E^k$ and $\beta $ in $E^l$, where $k+l=m+1$, and such that $\src
(\alpha )=\ran (\beta )$.  We once more leave the easy case $k=0$ to the reader and observe that the case $k=1$ follows from the
definition of $\varphi ^*$.

We may then suppose that $k\geq 2$, so we write $\alpha =\alpha '\alpha ''$, with $\alpha '\in E^1$, and $\alpha ''\in E^{k-1}$.  Then
  $$
  \varphi ^*(g,\alpha \beta ) = \varphi ^*(g,\alpha '\alpha ''\beta ) = \varphi ^*\big(\varphi (g,\alpha '),\alpha ''\beta \big) =
  \varphi ^*\Big (\varphi ^*\big(\varphi (g,\alpha '),\alpha ''\big),\beta \Big ) \$=
  \varphi ^*\Big (\varphi ^*(g,\alpha '\alpha ''\big),\beta \Big ) = \varphi ^*\Big (\varphi ^*(g,\alpha \big),\beta \Big ).
  $$

Let us now prove that $\auto ^*$ is in fact an action of $G$ on $E^n$.  We begin by proving that $\auto ^*_g\auto ^*_h =
\auto ^*_{gh}$ on $E^n$, for every $g$ and $h$ in $G$, which we do by induction on $n$.

This follows immediately from the hypothesis for $n\leq 1$, so let us assume that $n\geq 2$.  Given $\alpha \in E^n$, write $\alpha =\alpha '\alpha ''$,
with $\alpha '\in E^1$, and $\alpha ''\in E^{n-1}$. Then
  $$
  \auto ^*_g\big(\auto ^*_h(\alpha )\big) = \auto ^*_g\big(\auto ^*_h(\alpha '\alpha '')\big) =
  \auto ^*_g\big( \auto _h(\alpha ')\auto _{\varphi (h,\alpha ')}(\alpha '') \big) \$=
  \auto _g\big( \auto _h(\alpha ') \big) \auto ^*_{\varphi (g,\auto _h(\alpha '))}\big(\auto _{\varphi (h,\alpha ')}(\alpha '') \big) =
  \auto _{gh}(\alpha ') \auto ^*_{ \varphi (g,\auto _h(\alpha ')) \varphi (h,\alpha ')}(\alpha '') \$=
  \auto _{gh}(\alpha ') \auto ^*_{ \varphi (gh,\alpha ')}(\alpha '') =
  \auto ^*_{gh}(\alpha '\alpha '') = \auto ^*_{gh}(\alpha ).
  $$

  That $\alpha ^*_g$ is bijective on each $E^n$ then
  follows\fn {This is why it is useful to include \ref{PhiOneId} as a separate statement, since we may now use it
to prove bijectivity.}
  from \ref{PhiOneId}, so $\alpha ^*$ is indeed an action of $G$ on $E^n$.

  Finally, let us show that $\varphi ^*$ is a cocycle for $\auto ^*$ on $E^n$.  For this fix $g$ and $h$ in $G$ and let $\alpha
 \in E^n$.  Then, with $\alpha =\alpha '\alpha ''$, as before,
  $$
  \varphi ^*(gh, \alpha )=
  \varphi ^*(gh, \alpha '\alpha '')=
  \varphi ^*(\varphi (gh, \alpha '),\alpha '')=
  \varphi ^*\Big ( \varphi \big(g, \auto _h(\alpha ')\big) \varphi (h, \alpha '),\alpha ''\Big )\$=
  \varphi ^*\Big (\varphi \big(g, \auto _h(\alpha ')\big),\auto ^*_{\varphi (h, \alpha ')}(\alpha '')\Big ) \varphi ^*\big(\varphi (h, \alpha '),\alpha ''\big) =: (\star ).
  $$
  On the other hand, focusing on the right-hand-side of \ref{CocycleId}, notice that
  $$
  \varphi ^*(g,\auto ^*_h(\alpha ))\varphi ^*(h,\alpha ) =
  \varphi ^*\big(g,\auto ^*_h(\alpha '\alpha '')\big)\varphi ^*(h,\alpha '\alpha '') \$=
  \varphi ^*\big(g, \auto _h(\alpha ')\auto ^*_{\varphi (h,\alpha ')}(\alpha '')\big) \varphi ^*\big(\varphi (h,\alpha '),\alpha ''\big) \$=
  \varphi ^*\Big (\varphi \big(g, \auto _h(\alpha ') \big),\auto ^*_{\varphi (h,\alpha ')}(\alpha '') \Big ) \varphi ^*\big(\varphi (h,\alpha '),\alpha ''\big),
  $$
  which coincides with $(\star )$ above.  This concludes the proof.  \endProof

\bigskip

The only action of $G$ on $E^*$ to be considered in this paper is $\auto ^*$ so, from now on, we will adopt the shorthand
notation
  $$
  g\alpha = \auto ^*_g(\alpha ).
  $$
  Moreover, since $\varphi ^*$ extends $\varphi $, we will drop the star decoration and denote $\varphi ^*$ simply as $\varphi $.
  The group law, the cocycle condition, and properties
  \ref{CocZero},
  \ref{MatchRange},
  \ref{MatchSource},
  \ref{PhiStarOnVert},
  \ref{MainConcat} and
  \ref{LastItem}
  of Proposition \ref{extendedaction} may then be rewritten as follows:

\state Equations \label Equacoes
  For every $g$ and $h$ in $G$, for every $\vr \in E^0$, and for every $\alpha $ and $\beta $ in $E^*$ such that $\src (\alpha )=\ran (\beta )$,
one has that
  \iaitem
  \aitem $(gh)\alpha = g(h\alpha )$,
  \aitem $\varphi (gh, \alpha ) = \varphi \big(g,h\alpha \big)\varphi (h,\alpha ),$
  \Item {(\CocZero )} $\varphi (g,\vr ) = g$,
  \Item {(\MatchRange )} $\ran (g\alpha )=g\ran (\alpha )$,
  \Item {(\MatchSource )} $\src (g\alpha )=g\src (\alpha )$,
  \Item {(\PhiStarOnVert )} $\varphi (g,\alpha )\vr =g\vr $,
  \Item {(\MainConcat )} $g(\alpha \beta ) = (g\alpha )\ \varphi (g,\alpha )\beta $,
  \Item {(\LastItem )} $\varphi (g, \alpha \beta )=\varphi \big(\varphi (g,\alpha ),\beta \big)$.

It might be worth noticing that if $\varphi (g,\alpha )=1$, then \ref{Equacoes.\MainConcat } reads
  ``$g(\alpha \beta ) = (g\alpha )\beta $'',
  which may be viewed as an associativity property.  However associativity does not hold in general as $\varphi $ is not always
trivial, and hence parentheses must be used.

On the other hand parentheses are unnecessary in expressions of the form $\alpha g\beta $, when $\alpha ,\beta \in E^*$, and $g \in G$, since the
only possible interpretation for this expression is the concatenation of $\alpha $ with $g\beta $.

Another useful property of $\varphi $ is in order.

\state Proposition \label Inverses
  For every $g \in G$, and every $\alpha \in E^*$, one has that
  $$
  \varphi (g\inv ,\alpha ) = \varphi (g,g\inv \alpha )\inv .
  $$

\Proof We have
  $$
  1 =
  \varphi (1,\alpha ) =
  \varphi (gg\inv ,\alpha ) =
  \varphi (g,g\inv \alpha )\varphi (g\inv ,\alpha ),
  $$
  from where the conclusion follows.
  \endProof

\section The universal C*-algebra $\OGE $

As in the above section we fix a graph $E$, an action of a group $G$ on $E$, and a one-cocycle $\varphi $ satisfying \ref{StandingHyp}.

  It is our next goal to build a C*-algebra from this data but first let us recall the following notion from \cite
{Raeburn}:

\definition \label DefineCKFamily
  A \"{Cuntz-Krieger $E$-family} consists of a set
  $$
  \{p_\vr : \vr \in E^0\}
  $$
  of mutually orthogonal projections and a set
  $$
  \{\s _\ed : \ed \in E^1 \}
  $$
  of partial isometries, all lying in some C*-algebra, and satisfying
  \izitem
  \zitem $s_\ed ^* s_\ed = p_{\src (\ed )}$, for every $\ed \in E^1$, \zitemmark CkTwo
  \zitem
  $\ds
  p_\vr =\kern -7pt \sum _{\ed \in \ran \inv (\vr )}s_\ed s_\ed ^*,
  $
  for every $\vr \in E^0$ for which $\ran \inv (\vr )$ is finite and nonempty. \zitemmark CkOne

\definition \label DefineOGE
  We define $\OGE $ to be the universal unital C*-algebra generated by a set
  $$
  \{p_\vr : \vr \in E^0\}\cup \{\s _\ed : \ed \in E^1\} \cup \{u_g : g \in G\},
  $$
  subject to the following relations:
  \iaitem
  \aitem $\{p_\vr : \vr \in E^0\}\cup \{\s _\ed : \ed \in E^1\}$ is a Cuntz-Krieger $E$-family, \aitemmark CKRel
  \aitem the map $u:G\rightarrow \OGE $, defined by the rule $g\mapsto u_g$, is a unitary representation of $G$,
  \aitem $u_g\s _\ed =\s _{g\ed }u_{\varphi (g,\ed )}$, for every $g \in G$, and $\ed \in E^1$,
  \aitem $u_gp_\vr =p_{g\vr }u_g$, for every $g \in G$, and $\vr \in E^0$.

\bigskip Observe that, under our standing assumptions \ref{StandingHyp}, for every $\vr \in E^0$ we have that $\ran
\inv (\vr )$ is finite and nonempty.  So \ref{DefineCKFamily.ii} and \ref{DefineOGE.\CKRel } imply that
  $$
  u_gp_\vr u_g^* =
  \sum _{\ran (\ed )=\vr } u_g\s _\ed \s _\ed ^*u_g^* =
  \sum _{\ran (\ed )=\vr } \s _{g\ed }u_{\varphi (g,\ed )}u_{\varphi (g,\ed )}^*\s _{g\ed }^* \$=
  \sum _{\ran (\ed )=\vr } \s _{g\ed }\s _{g\ed }^* =
  \sum _{\ran (\oed )=g\vr } \s _{\oed }\s _{\oed }^* =
  p_{g\vr },
  $$
  which says that \ref{DefineOGE.d} follows from the other conditions.  We have nevertheless included it in \ref{DefineOGE} in the belief that our theory may be generalized to graphs with sources.

Our construction generalizes some well known constructions in the literature as we would now like to mention.

\state Example \rm \label Nekrashevych
  Let $(G,X)$ be a self similar group as in \cite [Definition 2.1]{NC}.
  We may then consider a graph $E$ having only one vertex and such that $E^1=X$.  If we define
  $$
  \varphi (g,\vr )=g|_\vr ,
  $$
  where, in the terminology of \cite {NC}, $g|_\vr $ is the restriction (or section) of $g$ at $\vr $,
  then the triple $(\Data )$ satisfies \ref{StandingHyp} and one may easily show that $\OGE $ is isomorphic to the
algebra ${\cal O}_{(G,X)}$ introduced by Nekrashevych in \cite {NC}.

\state Example \rm \label KatsuraExample
  As in \cite {KatsuraOne}, let us assume we are given two $N\times N$ matrices $A$ and $B$ with integer entries, and such
that $A_{i,j}\geq 0$, for all $i$ and $j$.  We may then
  consider the graph $E$ with vertex set
  $$
  E^0 = \{1,2,\ldots ,N\},
  $$
  and such that, for each pair of vertices $i,j \in E^0$, the set of edges from vertex $j$ to vertex $i$ is a set with
$A_{i,j}$ elements, say
  $$
  \{e_{i,j,n} : 0\leq n<A_{i,j}\}.
  $$

Assuming moreover that $A$ has no identically zero rows, it is easy to see that $E$ has no sources.

Define an action $\auto $ of ${\bf Z}$ on $E$, which is trivial on $E^0$, and which acts on edges as follows: given $m\in {\bf Z}$, and
  $e_{i,j,n} \in E^1$, let $(\hat k,\hat n)$ be the unique pair of integers such that
  $$
  mB_{i,j}+n=\hat k A_{i,j} + \hat {n} \and 0\leq \hat {n}<A_{i,j}.
  $$
  That is, $\hat k$ is the quotient and $\hat {n}$ is the remainder of the Euclidean division of $mB_{i,j}+n$ by
$A_{i,j}$.  We then put
  $$
  \auto _m(e_{i,j,n})=e_{i,j,\hat {n}}.
  $$
  In other words, $\auto _m$ corresponds to the addition of $mB_{i,j}$ to the variable ``$n$'' of ``$e_{i,j,n}$'', taken
modulo $A_{i,j}$.
  In turn, the one-cocycle is defined by
  $$
  \varphi (m, e_{i,j,n})= \hat k.
  $$

Observe that if $A_{i,j}=0$, then there are no edges from $j$ to $i$, so the value $B_{i,j}$ is entirely irrelevant for
the above construction.  Therefore it makes no difference to assume that
  $$
  A_{i,j}=0 \IMPLY B_{i,j}=0.
  $$

  It may then be proved without much difficulty that ${\cal O}_{\Z , E}$ is isomorphic to Katsura's \cite {KatsuraOne}
algebra $\OAB $, under an isomorphism sending each $u_m$ to the $m ^{th}$ power of the unitary
  $$
  u:=\sum _{i=1}^Nu_i
  $$
  in $\OAB $, and sending $\s _{e_{i,j,n}}$ to $\s _{i,j,n}$.

\bigskip When $N=1$, the relevant graph for Katsura's algebras is the same as the one we used above in the description
of Nekrashevych's example.  However the former is not a special case of the latter because, contrary to what is required
in \cite {NC}, the group action might not be faithful.

\state Example \rm \label CrossedProduct
  Given any finite graph $E$, and any action $\sigma $ of a group $G$ on $E$, the map $\varphi : G\times E^1 \rightarrow G$ defined by
  $$
  \varphi (g,a)=g \for g\in G \for a\in E^1
  $$
  is a one-cocycle, and the triple $(G,E, \varphi )$ satisfies \ref{StandingHyp}. By \ref{DefineOGE.c}, we have that
  $$
  u_gs_au_g^*=s_{ga},
  $$
  for any $g$ in $G$, and every $a$ in $ E^1$.
  It is therefore easy to see that
  $\OGE $ is isomorphic to the crossed product of the graph C*-algebra $C^*(E)$ \cite {Raeburn} by $G$, relative to the
natural action of $G$ on $C^*(E)$ induced by $\sigma $.
  In particular, if $\sigma $ is the trivial action, we have that $\OGE $ is the maximal tensor product of $C^*(E)$ by the full
group C*-algebra of $G$.

\state Example \rm
  Given any finite graph without sources, and any action $\sigma $ of a group $G$ on $E$ fixing the vertices, consider the map
  $\varphi : G\times E^1 \rightarrow G$ defined by
  $$
  \varphi (g,a)=1 \for g\in G \for a\in E^1.
  $$
  It is easy to see that $\varphi $
  is a one-cocycle, and that the triple $(G,E, \varphi )$ satisfies \ref{StandingHyp}. Since $E$ has no sources we have,
for any $g$ in $ G$, that
  $$
  u_g=
  \sum \limits _{x\in E^0}u_gp_x=
  \sum \limits _{x\in E^0}\sum \limits _{a\in r\inv (x)}u_gs_as_a^*\={DefineOGE.c}
  \sum \limits _{x\in E^0}\sum \limits _{a\in r\inv (x)}s_{ga}s_a^*,
  $$
  which therefore lies in the copy of $C^*(E)$ within $\OGE $.  Since the natural representation of
$C^*(E)$ in $\OGE $ is faithful by \ref{PropGraphSeaStarisSubalg}, the conclusion is that $\OGE \cong C^*(E)$.

We now return to the general case of a triple $(\Data )$ satisfying \ref{StandingHyp}.  We initially recall the
usual extension of the notation ``$\s _\ed $'' to allow for paths of arbitrary length.

\state Definition \label ExtendToPath
  Given a finite path $\alpha $ in $E^*$, we shall let $\s _\alpha $ denote the element of $\OGE $ given by:
  \izitem
  \zitem when $\alpha =\vr \in E^0$, we let $\s _\alpha = p_\vr $,
  \zitem when $\alpha \in E^1$, then $\s _\alpha $ is already defined above,
  \zitem when $\alpha \in E^n$, with $n>1$, write $\alpha = \alpha '\alpha ''$, with $\alpha ' \in E^1$, and $\alpha '' \in  E^{n-1}$, and set $\s _\alpha = \s _{\alpha '}\s
_{\alpha ''}$, by recurrence.

Commutation relation \ref{DefineOGE.c} may then be generalized to finite paths as follows:

\state Lemma \label RuleforFinitePaths
  Given $\alpha \in E^*$, and $g \in G$, one has that
  $$
  u_g\s _\alpha =\s _{g\alpha }u_{\varphi (g,\alpha )}.
  $$

\Proof
  Let $n$ be the length of $\alpha $.  When $n=0,1$, this follows from \ref{DefineOGE.d\&c}, respectively.  When $n>1$,
write $\alpha = \alpha '\alpha ''$, with $\alpha ' \in E^1$, and $\alpha '' \in E^{n-1}$.  Using induction, we then have
  $$
  u_g\s _{\alpha } =
  u_g\s _{\alpha '}\s _{\alpha ''} =
  \s _{g\alpha '}u_{\varphi (g,\alpha ')}\s _{\alpha ''} =
  \s _{g\alpha '}\s _{\varphi (g,\alpha ')\alpha ''} u_{\varphi (\varphi (g,\alpha '),\alpha '')} \$=
  \s _{(g\alpha ')\varphi (g,\alpha ')\alpha ''} u_{\varphi (g,\alpha '\alpha '')} =
  \s _{g(\alpha '\alpha '')} u_{\varphi (g,\alpha '\alpha '')} =
  \s _{g\alpha } u_{\varphi (g,\alpha )}.
  \endProof

Our next result provides a spanning set for $\OGE $.

\state Proposition \label InvSemPicture
  Let
  $$
  \S = \big \{\s _\alpha u_g\s _\beta ^* : \alpha ,\beta \in E^*,\ g \in G,\ \src (\alpha )=g \src (\beta )\big \} \cup \{0\}.
  $$
  Then $\S $ is closed under multiplication and adjoints and its closed linear span coincides with $\OGE $.

\Proof That $\S $ is closed under adjoints is clear.  With respect to closure under multiplication, let $\s _\alpha
u_g\s _\beta ^*$ and $\s _\gamma u_h\s _\delta ^*$ be elements of $\S $.

From \ref{DefineOGE.\CKRel } we know that $\s _\beta ^*\s _\gamma = 0$, unless either $\gamma =\beta \varepsilon $, or $\beta =\gamma \varepsilon $, for some $\varepsilon \in E^*$.  If
$\gamma =\beta \varepsilon $, then
  $$
  \s _\beta ^*\s _\gamma = \s _\beta ^*\s _{\beta \varepsilon } = \s _\beta ^*\s _\beta \s _\varepsilon = \s _\varepsilon ,
  $$
  and hence
  $$
  (\s _\alpha u_g\s _\beta ^*)(\s _\gamma u_h\s _\delta ^*)=
  \s _\alpha u_g\s _\varepsilon u_h\s _\delta ^*=
  \s _\alpha \s _{g \varepsilon }u_{\varphi (g,\varepsilon )}u_h\s _\delta ^* =
  \s _{\alpha g \varepsilon }u_{\varphi (g,\varepsilon )h}\s _\delta ^*.
  \equationmark ProdSpan
  $$
  Moreover, since
  $$
  \src (\alpha g \varepsilon ) =
  \src (g\varepsilon ) =
  g\src (\varepsilon ) =
  \varphi (g,\varepsilon )\src (\varepsilon ) =
  \varphi (g,\varepsilon )\src (\gamma ) =
  \varphi (g,\varepsilon )h\src (\delta ),
  $$
  we deduce that the element appearing in the right-hand-side of \ref{ProdSpan} indeed belongs to $\S $.

In the second case, namely if $\beta =\gamma \varepsilon $, then the adjoint of the term appearing in the left-hand-side of \ref{ProdSpan} is
  $$
  (\s _\delta u_{h\inv }\s _\gamma ^*)
  (\s _\beta u_{g\inv }\s _\alpha ^*),
  $$
  and the case already dealt with implies that this belongs to $\S $.  The result then follows from the fact
that $\S $ is self-adjoint.

In order to prove that $\OGE $ coincides with the closed linear span of $\S $, let $A$ denote the latter.  Given
that $\S $ is self-adjoint and closed under multiplication, we see that $A$ is a closed *-subalgebra of $\OGE $.
Since $A$ evidently contains $\s _\alpha $ for every $\alpha $ in $E^{\leq 1}$, and since it also contains $u_g$ for every $g$ in $G$,
we deduce that $A=\OGE $.  \endProof

\section The inverse semigroup $\SGE $

As before, we keep \ref{StandingHyp} in force.

In this section we will give an abstract description of the set $\S $ appearing in \ref{InvSemPicture} as
well as its multiplication and adjoint operation.  The goal is to construct an inverse semigroup from which we will
later recover $\OGE $.

\definition \label inverseSemigroup
  Over the set
  $$
  \SGE =\big \{ (\alpha ,g,\beta ) \in E^*\times G\times E^*: \src (\alpha )=g\src (\beta )\big \}\cup \{0\},
  $$
  consider a binary \"{multiplication} operation defined by
  $$
  (\alpha ,g,\beta ) (\gamma ,h,\delta ) = \left \{\matrix {
  (\alpha g \varepsilon ,\hfil \varphi (g,\varepsilon ) h,\hfil \delta ), & \text {if } \gamma =\beta \varepsilon , \cr \cr
  (\alpha ,\ g\varphi (h\inv ,\varepsilon )\inv ,\ \delta h\inv \varepsilon ), & \text {if } \beta =\gamma \varepsilon , \cr \cr
  0,\hfil \hfil \hfil & \text {otherwise,}
  }\right .
  $$
  and a unary \"{adjoint} operation defined by
  $$
  (\alpha ,g,\beta )^*:= (\beta ,g\inv , \alpha ).
  $$
  Furthermore, the subset of $\SGE $ formed by all elements $(\alpha ,g,\beta )$, with $g=1$, will be denoted by $\SE $.

It is easy to see that $\SE $ is closed under the above operations, and that it is isomorphic to the inverse semigroup
generated by the canonical partial isometries in the graph C*-algebra of $E$.

Let us begin with a simple, but useful result:

\state Lemma \label EasyMult
  Given $(\alpha ,g,\beta )$ and $(\gamma ,h,\delta )$ in $\SGE $, one has
  $$
  \beta =\gamma \IMPLY (\alpha ,g,\beta ) (\gamma ,h,\delta ) = (\alpha ,gh,\delta ).
  $$

\Proof
  Focusing on the first clause of \ref{inverseSemigroup}, write $\gamma =\beta \varepsilon $, with $\varepsilon =\src (\beta )$.  Then
  $$
  (\alpha ,g,\beta ) (\gamma ,h,\delta )=
  (\alpha g \src (\beta ),\ \varphi \big(g,\src (\beta )\big) h,\ \delta ) =
  (\alpha \src (\alpha ),\ gh,\ \delta ) = (\alpha ,gh, \delta ).
  \endProof

\state Proposition $\SGE $ is an inverse semigroup with zero.

\Proof We leave it for the reader to prove that the above operations are well defined and the multiplication is
associative.  In order to prove the statement it then suffices \cite [Theorem 1.1.3]{Lawson} to show that, for all
$y,z \in \SGE $, one has that
  \Zitem $yy^*y=y$, and
  \zitem $yy^*$ commutes with $zz^*$.

\bigskip

Given $y=(\alpha ,g,\beta ) \in \SGE $, we have by the above Lemma that
  $$
  yy^*y =
  (\alpha ,g,\beta )(\beta , g\inv , \alpha )(\alpha ,g,\beta )=
  (\alpha , 1, \alpha )(\alpha ,g,\beta ) =
  (\alpha ,g,\beta ) = y,
  $$
  proving (i).
  Notice also that
  $$
  yy^* = (\alpha , 1, \alpha )
  \equationmark IdempotForm
  $$
  is an element of the idempotent semi-lattice of $\SE $, which is a commutative set because $\SE $ is an inverse
semigroup.  Point (ii) above then follows immediately, concluding the proof.
  \endProof

As seen in \ref{IdempotForm}, the idempotent semi-lattice of $\SGE $, henceforth denoted by $\EGE $, is given by
  $$
  \EGE = \big \{(\alpha ,1,\alpha ): \alpha \in E^*\big \} \cup \{0\}.
  \equationmark ISLDef
  $$
  Evidently $\EGE $ is also the idempotent semi-lattice of $\SE $.

  For simplicity, from now on we will adopt the short-hand notation
  $$
  \eproj _\alpha = (\alpha ,1,\alpha ) \for \alpha \in E^*.
  \equationmark DefineEAlpha
  $$

The following is a standard fact in the theory of graph C*-algebras:

\state Proposition \label ProdIdemp
  If $\alpha ,\beta \in E^*$, then
  $$
  \eproj _\alpha \eproj _\beta = \left \{\matrix {
  \eproj _\alpha , & \text {if there exists $\gamma $ such that } \alpha =\beta \gamma , \cr
  \eproj _\beta , & \text {if there exists $\gamma $ such that } \alpha \gamma =\beta , \pilar {12pt} \cr
  \hfill 0, & \text {otherwise.}\hfill \pilar {12pt}
  }\right .
  $$

Recall that if $\alpha $ and $\beta $ are in $E^*$, we say that $\alpha \preceq \beta $, if $\alpha $ is a \"{prefix} of $\beta $, i.e.~if there exists
$\gamma \in E^*$, such that $\alpha \gamma = \beta $.  It therefore follows from \ref{ProdIdemp} that
  $$
  \eproj _\alpha \leq \eproj _\beta \iff \beta \preceq \alpha .
  \equationmark WrodVsIdemp
  $$

Another easy consequence of
  \redundantref {Proposition}{ProdIdemp} is that, for any two elements $e,f \in \EGE $, one has
that either $e\perp f$, or $e$ and $f$ are comparable.  It follows that
  $$
  e\Cap f \IMPLY e\leq f \text {, \ or \ } f\leq e.
  \equationmark Compara
  $$

\section Pseudo freeness and E*-unitarity

\label PseudoFreeSect Again working under \ref{StandingHyp}, suppose we are given $g$ in $G$ and a finite path $\alpha $ such that
  $$
  g\alpha = \alpha  \and \varphi (g,\alpha )=1.
  \equationmark StrFixed
  $$
  Then, given any finite path $\alpha '$ extending $\alpha $, that is a path of the form $\alpha '=\alpha \beta $, where $\beta $ is another finite path,
we have
  $$
  g\alpha ' = g(\alpha \beta ) = (g\alpha )\varphi (g,\alpha )\beta  = \alpha \beta  = \alpha ',
  $$
  and
  $$
  \varphi (g,\alpha ') = \varphi (g,\alpha \beta ) = \varphi \big(\varphi (g,\alpha ),\beta \big) = \varphi \big(1,\beta \big) = 1.
  $$
  This says that any path $\alpha '$ extending $\alpha $ also satisfies \ref{StrFixed} so, in particular, every extension of
$\alpha $ is fixed by $g$.

\definition
  If $g\in G$ and $\alpha $ is a finite path satisfying \ref{StrFixed},
  we will say that $\alpha $ is \"{strongly fixed} by $g$.  In addition, if no proper prefix of $\alpha $ is strongly fixed by $g$,
we will say that $\alpha $ is a \"{minimal} strongly fixed path for $g$.

The following result is an easy consequence of the discussion above:

\state Proposition \label StrFixElts
  Given $g$ in $G$, let $M_g$ be the set of all minimal strongly fixed paths for $g$.  Then the set of all strongly
fixed paths for $g$ is given by
  $$
  \textstyle \bigsqcup \limits _{\mu \in M_g}\{\mu \gamma : \gamma \in E^*,\ \src (\mu ) = \ran (\gamma )\},
  $$
  where the square cup stands for disjoint union.

Let us now introduce terminology to describe situations in which nontrivial strongly fixed paths do not exist.

\definition \label EssFree
  We will say that $(\Data )$ is \"{pseudo free}\fn
    {In a preprint version of this work we have used the term \"{residually free} to refer to the concept presently
    being defined, but this apparently conflicts with a well established notion in group theory.}
  if, whenever $(g,\ed ) \in G\times E^1$, is such that $g\ed = \ed $, and $\varphi (g,\ed )=1$, then $g=1$.

Notice that pseudo freeness is equivalent to the fact that an edge is never a strongly fixed path for a nontrivial group
element.  In fact we may boost this up to finite paths as follows:

\state Proposition \label EssFreePath
  Suppose that $(\Data )$ is pseudo free and that a finite path $\alpha $ of nonzero length is strongly fixed for some $g$ in
$G$.  Then $g=1$.

\Proof Arguing by contradiction, assume that there is a counter-example to the statement, meaning that there is a
strongly fixed path $\alpha $ for a nontrivial group element $g$.  Then, as already mentioned, $\alpha $ has a minimal strongly
fixed prefix, so we may assume without loss of generality that $\alpha $ itself is minimal.

By \ref{Equacoes.\CocZero }, $\alpha $ can't be a vertex, and neither can it be an edge, by hypothesis.  So $|\alpha |\geq 2$, and
we may then write $\alpha =\beta \gamma $, with $\beta ,\gamma \in E^*$, and $|\beta |,|\gamma |<|\alpha |$.  Then
  $$
  \beta \gamma = \alpha = g\alpha = g(\beta \gamma ) = (g\beta )\varphi (g,\beta )\gamma ,
  $$
  whence $\beta =g\beta $, and $\gamma = \varphi (g,\beta )\gamma $, by length considerations.  Should $\varphi (g,\beta )=1$, the pair $(g,\beta )$ would be a smaller
counter-example to the statement, violating the minimality of $\alpha $.  So we have that $\varphi (g,\beta )\neq 1$.  In addition,
  $$
  \varphi \big(\varphi (g,\beta ),\gamma \big) = \varphi (g,\beta \gamma ) = \varphi (g,\alpha ) = 1.
  $$
  It follows that $\big(\varphi (g,\beta ),\gamma \big)$ is a counter-example to the statement, again violating the minimality of $\alpha $.  This is a
contradiction and hence no counter-example exists whatsoever, concluding the proof.
  \endProof

An apparently stronger version of pseudo freeness is in order.

\state Proposition \label VarStar
  Suppose that $(\Data )$ is pseudo free.  Then, for all $g_1,g_2 \in G$, and $\alpha \in E^*$, one has that
  $$
  g_1\alpha =g_2\alpha \text { \ and \ } \varphi (g_1,\alpha )=\varphi (g_2,\alpha ) \IMPLY g_1=g_2.
  $$

\Proof
  Defining $g=g_2\inv g_1$, observe that
  $
  g\alpha =\alpha ,
  $
  and we claim that $\varphi (g,\alpha ) = 1$.
  In fact,
  $$
  \varphi (g,\alpha ) =
  \varphi \big(g_2\inv g_1,\alpha \big) =
  \varphi \big(g_2\inv ,g_1\alpha \big) \varphi \big(g_1,\alpha \big) \={Inverses} $$$$ =
  \varphi \big(g_2,g_2\inv g_1\alpha \big)\inv \varphi (g_1,\alpha ) =
  \varphi (g_2,\alpha )\inv \varphi (g_1,\alpha ) =1,
  $$
  so it follows that $g=1$, which is to say that $g_1=g_2$.
  \endProof

We will now determine conditions under which $\SGE $ is E*-unitary.  In order to do so we first need to understand when
does an element $s$ of $\SGE $ dominate a nonzero idempotent $e$, which in turn must necessarily have the form
$e=(\gamma ,1,\gamma )$, as seen in \ref{IdempotForm}.  If $s$ indeed dominates a nonzero idempotent, it is clear that $s$ is
itself nonzero, so $s$ must have the form $(\alpha ,g,\beta )$.

\state Proposition \label LemDominateIdempotent
  Let $\alpha $, $\beta $ and $\gamma $ be finite paths in $E$, and let $g\in G$ be such that $\src (\alpha )=g\src (\beta )$, so that $s:=(\alpha ,g,\beta )$ is
a general nonzero element of $\SGE $ and $e:=(\gamma ,1,\gamma )$ is a general nonzero idempotent element of $\SGE $.  Then $e\leq s$,
  if and only
  \izitem
  \zitem $\alpha =\beta $,
  \zitem $\gamma =\alpha \tau $, for some finite path $\tau $,
  \zitem $\tau $ is strongly fixed by $g$.

\Proof In order to prove the ``if'' part, we have
  $$
  se =
  (\alpha ,g,\beta )(\gamma ,1,\gamma ) =
  (\alpha ,g,\alpha )(\alpha \tau ,1,\gamma ) =
  (\alpha g\tau ,\varphi (g,\tau ),\gamma ) \$=
  (\alpha \tau ,1,\gamma ) =
  (\gamma ,1,\gamma ) =
  e,
  $$
  proving that $e\leq s$.  Conversely, assuming that $e\leq s$, we have $se=e$, so in particular $se\neq 0$, and hence by the
definition of the multiplication on $\SGE $, either $\gamma $ is a prefix of $\beta $ or vice versa.

In case $\beta $ is a prefix of $\gamma $, we may write $\gamma =\beta \tau $, for some finite path $\tau $, and then
  $$
  (\gamma ,1 , \gamma ) = e = se =
  (\alpha , g,\beta )(\beta \tau  ,1 , \gamma ) =
  (\alpha g\tau  ,\varphi (g,\tau ) , \gamma ),
  $$
  so we conclude that
  $$
  \alpha g\tau =\gamma = \beta \tau \and \varphi (g,\tau )=1.
  $$
  So $\alpha =\beta $, $g\tau =\tau $ and the statement is proved.

On the other hand, if $\gamma $ is a prefix of $\beta $, we may write $\beta =\gamma \varepsilon $ and, again according to the definition of the
multiplication on $\SGE $,
  the third coordinate of the product $se$ will be $\gamma \varepsilon $, from where we conclude that $\gamma =\gamma \varepsilon $.  So $|\varepsilon |=0$ and then $\gamma =\beta $,
which in particular means that $\beta $ is a prefix of $\gamma $, and the proof follows as above.  \endProof

\state Proposition \label EssFreeEUnitary
  $\SGE $ is an E*-unitary inverse semigroup if and only if $(\Data )$ is pseudo free.

\Proof Let $s$ be an element of $\SGE $ which dominates a nonzero idempotent element $e$.  As discussed above, we
necessarily have
  $$
  s=(\alpha ,g,\beta ) \and e = (\gamma ,1,\gamma ),
  $$
  where $\alpha $, $\beta $ and $\gamma $ are finite paths in $E$, and $\src (\alpha )=g\src (\beta )$.  Then, by \ref{LemDominateIdempotent}
we conclude that $g\tau =\tau $, and $\varphi (g, \tau )=1$ so, assuming that $(\Data )$ is pseudo free, we have $g=1$.  Moreover by \ref{LemDominateIdempotent.i} we see that $\alpha =\beta $, so
  $$
  s = (\alpha ,g,\beta ) = (\alpha ,1,\alpha ),
  $$
  which is idempotent as desired.
  In order to prove the converse, let $(g,\ed )\in G\times E^1$, be such that $g\ed =\ed $, and $\varphi (g,\ed )=1$.  Then the element
  $$
  s:=\big(\src (\ed ),g,\src (\ed )\big)
  $$
  lies in $\SGE $ because
  $$
  g\src (\ed ) = \src (g\ed ) = \src (\ed ).
  $$
  Moreover
  observe that $s$ dominates the nonzero idempotent element $(\ed ,1,\ed )$, since
  $$
  s (\ed ,1,\ed ) =
  \big(\src (\ed ),g,\src (\ed )\big) \big(\src (\ed )\ed ,1,\ed \big) =
  (\src (\ed )g\ed ,\varphi (g,\ed ),\ed ) = (\ed ,1,\ed ).
  $$
  So, under the hypothesis that $\SGE $ is E*-unitary, we conclude that $s$ is idempotent, which is to say that $g=1$.
This proves that $(\Data )$ is pseudo free.  \endProof

\section Tight representations of $\SGE $

As before, we keep \ref{StandingHyp} in force.

It is the main goal of this section to show that $\OGE $ is the universal C*-algebra for tight representations of $\SGE
$.

Recall from
  \redundantref {Proposition}{ProdIdemp} that $\eproj _\alpha \leq \eproj _{\src (\alpha )}$, for every $\alpha \in E^*$, so we see that the set
  $$
  \{\eproj _\vr : \vr \in E^0\}
  \equationmark FiniteCover
  $$
  is a \"{cover} \cite [Definition 11.5]{actions} for $\EGE $.

\state Proposition \label TightRep
  The map
  $$
  \pi : \SGE \to \OGE ,
  $$
  defined by $\pi (0)=0$, and
  $$
  \pi (\alpha ,g,\beta ) = \s_\alpha u_g\s_\beta ^*,
  $$
  is a tight \cite [Definition 13.1]{actions} representation.

\Proof We leave it for the reader to show that $\pi $ is in fact multiplicative and that it preserves adjoints.

In order to prove that $\pi $ is tight, we shall use the characterization given in \cite [Proposition 11.8]{actions},
observing that $\pi $ satisfies condition (i) of \cite [Proposition 11.7]{actions} because, with respect to the cover
\ref{FiniteCover}, we have that
  $$
  \bigvee_{\vr \in E^0} \pi (\eproj_\vr ) =
  \bigvee_{\vr \in E^0} \pi (\vr ,1,\vr ) =
  \bigvee_{\vr \in E^0} p_\vr =
  \sum _{\vr \in E^0} p_\vr =
  1,
  $$
  by \ref{DefineOGE.\CKRel }.
  So we assume that $\{\eproj_{\alpha ^1},\ldots ,\eproj_{\alpha ^n}\}$ is a cover for a given $\eproj_\beta $, where $\alpha ^1,\ldots \alpha ^n,\beta \in E^*$, and we
need to show that
  $$
  \bigvee_{i=1}^n\pi (\eproj_{\alpha ^i}) \geq \pi (\eproj_\beta ).
  \equationmark TightGoal
  $$

In particular, for each $i$, we have that $\eproj_{\alpha ^i}\leq \eproj_\beta $, so by \ref{WrodVsIdemp} there exists $\gamma ^i \in E^*$
such that $\alpha ^i=\beta \gamma ^i$.

We shall prove \ref{TightGoal} by induction on the variable
  $$
  L = \min_{1\leq i\leq n}|\gamma ^i|.
  $$

  If $L=0$, we may pick $i$ such that $|\gamma ^i|=0$, and then necessarily $\gamma ^i= \src (\beta )$, in which case $\alpha ^i=\beta $, and \ref{TightGoal} is trivially true.

Assuming that $L\geq 1$,
  let $\vr := \src (\beta )$.  Observe that $x$ is not a source either because this is part of our standing assumptions
\ref{StandingHyp}, or simply because $x$ is the range of every $\gamma ^i$.  In any case let us write
  $$
  \ran \inv (\vr ) = \{\ed_1,\ldots ,\ed_k\},
  $$
  and observe that
  $$
  \pi (\eproj_\beta ) = \s_\beta \s_\beta ^* = \s_\beta p_\vr \s_\beta ^* \={DefineOGE.\CKRel }
  \sum _{j=1}^k\s_\beta \s_{\ed_j}\s_{\ed_j}^*\s_\beta ^* =
  \sum _{j=1}^k \pi (\eproj_{\beta \ed_j}).
  $$
  In order to prove \ref{TightGoal} it is therefore enough to show that
  $$
  \bigvee_{i=1}^n\pi (\eproj_{\alpha ^i}) \geq \pi (\eproj_{\beta \ed_j}),
  \equationmark SmallGoal
  $$
  for all $j=1,\ldots ,k$.  Fixing $j$ we claim that $\eproj_{\beta \ed_j}$ is covered by the set
  $$
  Z = \big \{\eproj_{\alpha ^i}: 1\leq i\leq n,\ \eproj_{\alpha ^i}\leq \eproj_{\beta \ed_j}\big \}.
  $$

  In order to see this let $y$ be a nonzero element in $\EGE $ such that $y\leq \eproj_{\beta \ed_j}$.  Then $y\leq \eproj_{\beta }$,
and so $y\Cap \eproj_{\alpha ^i}$ for some $i$.  Thus, to prove the claim it is enough to check that $\eproj_{\alpha ^i}$ lies in
$Z$.
  Observe that
  $$
  y \eproj_{\beta \ed_j} \eproj_{\alpha ^i} =
  y \eproj_{\alpha ^i} \neq 0,
  $$
  which implies that $\eproj_{\beta \ed_j} \Cap \eproj_{\alpha ^i}$.

By \ref{Compara} we have that
  $\eproj_{\beta \ed_j}$ and $\eproj_{\alpha ^i}$ are comparable, so either $\beta \ed_j\preceq \alpha ^i$ or $\alpha ^i\preceq \beta \ed_j$, by
\ref{WrodVsIdemp}.
  Since we are under the hypothesis that $L\geq 1$, and hence that
  $$
  |\alpha ^i| =
  |\beta ^i| + |\gamma ^i| \geq
  |\beta |+1 =
  |\beta \ed_j|,
  $$
  we must have that $\beta \ed_j\preceq \alpha ^i$, from where we deduce that $\eproj_{\alpha ^i}\leq \eproj_{\beta  \ed_j}$, proving our
claim.

Employing the induction hypothesis we then deduce that
  $$
  \bigvee_{z \in Z}\pi (z) \geq \pi (\eproj_{\beta \ed_j}),
  $$
  verifying \ref{SmallGoal}, and thus concluding the proof.
  \endProof

We would now like to prove that the representation $\pi $ above is in fact the \"{universal} tight representation of $\SGE
$.

  \def \ts {\tilde \s }
  \def \tp {\tilde p}
  \def \tu {\tilde u}

\state Theorem Let $A$ be a unital C*-algebra and let $\rho :\SGE \to A$ be a tight representation.  Then there exists a
unique unital *-homomorphism $\psi :\OGE \to A$, such that the diagram
  \hfill \break
  \vbox {
  \beginpicture
  \setcoordinatesystem units <0.0040truecm, -0.0040truecm> point at 3000 0
  \setplotarea x from -1500 to 1000, y from 500 to 1000
  \put {
    $\matrix {
      \SGE & \buildrel \ds \pi \over \longrightarrow & \OGE \cr \cr
      && \ \Big \downarrow \psi \cr \cr
      && A \
      }$
    } at 500 500
  \arrow <0.11cm> [0.3,1.2] from 300 400 to 650 650 \put {$\rho $} at 400 600
  \endpicture
  }
  \hfill \break commutes.

\Proof We will initially prove that the elements
  $$
  \matrix {
  \tp _\vr := \rho (\vr ,1,\vr ), \hfill & \forall \vr \in E^0, \hfill \cr \cr
  \ts _\ed := \rho \big(\ed , 1, \src (\ed )\big), \hfill & \forall \ed \in E^1, \hfill \cr \cr
  \tu _g := \ds \sum _{\vr \in E^0}\rho (\vr , g, g\inv \vr ), \hfill & \forall g \in G, \hfill }
  $$
  satisfy relations \ref{DefineOGE.a--d}.
  Since the $\eproj _\vr $ defined in \ref{DefineEAlpha} are mutually orthogonal idempotents in $\SGE $, it is
clear that the $\tp _\vr $ are mutually orthogonal projections.  Evidently the $\ts _\ed $ are partial isometries so, in
order to check \ref{DefineOGE.a}, we must only verify \ref{DefineCKFamily.\CkTwo } and \ref{DefineCKFamily.\CkOne }.  With respect to the former, let $\ed \in E^1$.  Then
  $$
  \ts _\ed ^* \ts _\ed =
  \rho \big(( \src (\ed ), 1, \ed )(\ed , 1, \src (\ed )\big) =
  \rho \big(\src (\ed ), 1, \src (\ed )\big) =
  \tp _{\src (\ed )},
  $$
  proving \ref{DefineCKFamily.\CkTwo }.  In order to prove \ref{DefineCKFamily.\CkOne }, let $\vr $ be a
vertex such that $\ran \inv (\vr )$ is nonempty and write
  $$
  \ran \inv (\vr ) = \big \{\ed _1,\ldots ,\ed _n\big \}.
  $$

  Putting $q_i = (\ed _i,1,\ed _i)$,
  we then claim that the set
  $$
  \big \{q_1,\ldots ,q_n\big \}
  $$
  is a cover for $q:= (\vr ,1,\vr )$.  In order to prove this we must show that, if the nonzero idempotent $f$ is
dominated by $q$, then $f\Cap q_i$ for some $i$.

Let $f=(\alpha ,1,\alpha )$ by \ref{ISLDef} and notice that
  $$
  0 \neq f = fq = (\alpha ,1,\alpha )(\vr ,1,\vr ).
  $$
  So $\alpha $ and $\vr $ are comparable, and this can only happen when $\vr =\ran (\alpha )$.  If $|\alpha | = 0$ then necessarily $\alpha =\vr
$, so $f=q$, and it is clear that $f\Cap q_i$ for all $i$.  On the other hand, if $|\alpha |\geq 1$, we write
  $$
  \alpha =\alpha '\alpha '',
  $$
  with $\alpha ' \in E^1$, so that $\ran (\alpha ') = \ran (\alpha )=\vr $, and hence $\alpha '=\ed _i$, for some $i$.
  Therefore
  $$
  fq_i =
  (\alpha ,1,\alpha )(\ed _i,1,\ed _i) =
  (\alpha ,1,\alpha )(\alpha ',1,\alpha ') =
  (\alpha ,1,\alpha )\neq 0,
  $$
  so $f\Cap q_i$, proving the claim.  Since $\rho $ is a tight representation, we deduce that
  $$
  \rho (q) = \bigvee _{i=1}^n\rho (q_i),
  $$
  but since the $q_i$ are easily seen to be pairwise orthogonal, their supremum coincides with their sum, whence
  $$
  \tp _\vr = \rho (q) = \sum _{i=1}^n\rho (q_i) = \sum _{i=1}^n\rho (\ed _i,1,\ed _i) \$=
  \sum _{i=1}^n \rho \big((\ed _i,1,\src (\ed _i))\ (\src (\ed _i),1,\ed _i) \big) =
  \sum _{i=1}^n \ts _{\ed _i}\ts _{\ed _i}^*,
  $$
  thus verifying \ref{DefineCKFamily.\CkOne }, and hence proving \ref{DefineOGE.a}.

With respect to \ref{DefineOGE.b}, let us first prove that $\tu _1 = 1$.  Considering the subsets of $\EGE $ given
by
  $$
  X = \ifundef {varnothing} \emptyset \else \varnothing \fi , \quad Y = \ifundef {varnothing} \emptyset \else \varnothing \fi \and Z = \big \{(\vr ,1,\vr ): \vr \in E^0\big \},
  $$
  notice that, according to \cite [Definition 11.4]{actions}, one has that
  $$
  \EGE ^{X,Y}=\EGE ,
  $$
  and that $Z$ is a cover for $\EGE ^{X,Y}$, as seen in \ref{FiniteCover}. By the tightness condition \cite
[Definition 11.6]{actions} we have
  $$
  \bigvee _{z \in Z}\rho (z) \geq \bigwedge _{\vr \in X} \rho (\vr ) \wedge \bigwedge _{y \in Y} \neg {\rho (y)}.
  $$
  As explained in the discussion following \cite [Definition 11.6]{actions}, the right-hand-side above must be
interpreted as 1 because $X$ and $Y$ are empty.  On the other hand, since the $\rho (z)$ are pairwise orthogonal, the
supremum in the left-hand-side above becomes a sum, so
  $$
  1 = \sum _{z \in Z}\rho (z) = \sum _{\vr \in E^0}\rho (\vr , 1, \vr ) = \tu _1.
  $$

In order to prove that $\tu $ is multiplicative, let $g$ and $h$ be in $G$.  Then
  $$
  \tu _g \tu _h =
  \sum _{\vr ,\vro \in E^0}\rho \big((\vr , g, g\inv \vr ) (\vro , h, h\inv \vro )\big) \$=
  \sum _{\vr \in E^0}\rho \big((\vr , g, g\inv \vr ) (g\inv \vr , h, h\inv g\inv \vr )\big) =
  \sum _{\vr \in E^0}\rho \big(\vr , gh, (gh)\inv \vr \big)=
  \tu _{gh}.
  $$

We next claim that $\tu _g^* =\tu _{g\inv }$, for all $g$ in $G$.  To prove it we compute
  $$
  \tu _g^* =
  \sum _{\vr \in E^0}\rho (\vr , g, g\inv \vr )^* =
  \sum _{\vr \in E^0}\rho (g\inv \vr , g\inv , \vr ) = \cdots
  $$
  which, upon the change of variables $\vro =g\inv \vr $, becomes
  $$
  \cdots = \sum _{\vro \in E^0}\rho (\vro , g\inv , gy) = \tu _{g\inv }.
  $$

This shows that $\tu $ is a unitary representation, verifying \ref{DefineOGE.b}.  Turning now our attention to
\ref{DefineOGE.c}, let $g \in G$ and $\ed \in E^1$.  Then
  $$
  \tu _g\ts _\ed =
  \sum _{\vr \in E^0}\rho (\vr , g, g\inv \vr )\, \rho \big(\ed , 1, \src (\ed )\big) =
  \rho \big(g\ran (\ed ), g, \ran (\ed )\big)\,\rho \big(\ed , 1, \src (\ed )\big) \$=
  \rho \big(\ran (g\ed )g\ed , \varphi (g,\ed ), \src (\ed )\big) =
  \rho \big(g\ed , \varphi (g,\ed ), \src (\ed )\big) = (\star ).
  $$
  On the other hand
  $$
  \ts _{g\ed }\tu _{\varphi (g,\ed )} =
  \rho \big(g\ed ,1,\src (g\ed )\big) \sum _{\vr \in E^0} \rho \big(\vr , \varphi (g,\ed ), \varphi (g,\ed )\inv \vr \big) \$=
  \rho \big(g\ed ,1,\src (g\ed )\big) \,\rho \big(\src (g\ed ), \varphi (g,\ed ), g\inv \src (g\ed )\big) \$=
  \rho \big(g\ed ,1,\src (g\ed )\big)\,\rho \big(\src (g\ed ), \varphi (g,\ed ), \src (\ed )\big) =
  \rho \big(g\ed , \varphi (g,\ed ), \src (\ed )\big),
  $$
  which coincides with $(\star )$ and hence proves \ref{DefineOGE.c}.  We leave the proof of \ref{DefineOGE.d} to
the reader after which the universal property of $\OGE $ intervenes to provide us with a *-homomorphism
  $$
  \psi :\OGE \to A
  $$
  sending
  $$
  p_\vr \mapsto \tp _\vr ,\quad \s _\ed \mapsto \ts _\ed \and u_g \mapsto \tu _g.
  $$

Now we must show that
  $$
  \psi \big(\pi (\gamma )\big)=\rho (\gamma ) \for \gamma \in \SGE .
  \equationmark PsiPiRho
  $$
  We will first do so for the following special cases:
  \Zitem $\gamma = (\vr ,1,\vr )$, for $\vr \in E^0$,
  \zitem $\gamma = \big(\ed , 1, \src (\ed )\big)$, for $\ed \in E^1$,
  \zitem $\gamma = (\vr ,g,g\inv \vr )$, for $\vr \in E^0$, and $g \in G$.

\bigskip \noindent In case (i) we have
  $$
  \psi (\pi (\gamma )) =
  \psi (\pi (\vr ,1,\vr )) =
  \psi (p_\vr ) =
  \tp _\vr =
  \rho (\vr ,1,\vr ) =
  \rho (\gamma ).
  $$
  As for (ii),
  $$
  \psi (\pi (\gamma )) =
  \psi \big(\pi \big(\ed , 1, \src (\ed )\big)\big) =
  \psi (\s _\ed ) =
  \ts _\ed =
  \rho \big(\ed , 1, \src (\ed )\big) =
  \rho (\gamma ).
  $$
  Under (iii),
  $$
  \psi (\pi (\gamma )) =
  \psi \big(\pi (\vr ,g,g\inv \vr ) \big) =
  \psi (p_\vr u_g p_{g\inv \vr }) =
  \psi (p_\vr u_g) =
  \tp _\vr \tu _g \$=
  \rho (\vr ,1,\vr ) \sum _{\vro \in E^0}\rho (\vro , g, g\inv \vro ) =
  \sum _{\vro \in E^0}\rho \big((\vr ,1,\vr ) (\vro , g, g\inv \vro )\big) =
  \rho (\vr , g, g\inv \vr ) =
  \rho (\gamma ).
  $$

In order to prove \ref{PsiPiRho}, it is now clearly enough to check that the *-sub-semigroup of $\SGE $ generated
by the elements mentioned in (i--iii), above, coincides with $\SGE $.

Denoting this *-sub-semigroup by ${\cal T}$, we will first show that $\big(\alpha ,1,\src (\alpha )\big)$ is in ${\cal T}$, for every $\alpha \in E^*$.  This is
evident for $|\alpha |\leq 1$, so we suppose that $\alpha =\alpha '\alpha ''$, with $\alpha ' \in E^1$, and $\ran (\alpha '')=\src (\alpha ')$.  We then have by induction
that
  $$
  {\cal T}\ni
  \big(\alpha ',1,\src (\alpha ')\big) \big(\alpha '',1,\src (\alpha '')\big) =
  \big(\alpha' \alpha '',1,\src (\alpha '')\big) =
  \big(\alpha ,1,\src (\alpha )\big).
  $$

Considering a general element $(\alpha ,g,\beta ) \in \SGE $,
  let $\vr =\src (\alpha )$, so that $g\inv \vr =\src (\beta )$, and notice that
  $$
  {\cal T}\ni
  \big(\alpha ,1,\src (\alpha )\big) (\vr , g, g\inv \vr ) \big(\beta ,1,\src (\beta )\big)^* \$=
  \big(\alpha ,1,\src (\alpha )\big) \big(\src (\alpha ), g, \src (\beta )\big) \big(\src (\beta ),1,\beta \big) =
  (\alpha ,g,\beta ),
  $$
  which proves that ${\cal T}=\SGE $, and hence that \ref{PsiPiRho} holds.

To conclude we observe that the uniqueness of $\psi $ follows from the fact that $\OGE $ is generated by the $p_\vr $, the
$\s _\ed $, and the $u_g$.
  \endProof

Given an inverse semigroup $\S $ with zero, recall from \cite [Theorem 13.3]{actions} that $\G \tight (\S )$
(denoted simply as $\G \tight $ in \cite {actions}) is the groupoid of germs for the natural action of $\S $ on the
space of tight filters over its idempotent semi-lattice.  Moreover the C*-algebra of $\G \tight (\S )$ is universal
for tight representations of $\S $.

\state Corollary \label UniversalTightAlgebra
Under the assumptions of \ref{StandingHyp} one has that $\OGE $ is isomorphic to the C*-algebra of
the groupoid $\GpdGE $.

\Proof Follows from \cite [Theorem 13.3]{actions} and the uniqueness of universal C*-algebras. \endProof

We should notice that our requirement that $G$ be countable in \ref{StandingHyp} is only used in the above proof,
since the application of \cite [Theorem 13.3]{actions} depends on the countability of $\SGE $.

\section The Lag Group

It is our next goal to give a concrete description of $\GpdGE $, similar to the description given of the groupoid
associated to a row-finite graph in \cite [Definition 2.3]{KPRR}.  The crucial ingredient there is the notion of
\"{tail equivalence with lag}.  In this section we will construct a group where our generalized \"{lag} function will
take values.

Let $G$ be a group.  Within the infinite cartesian product\fn {For the purpose of this cartesian product we adopt the
convention that ${\bf N}=\{1,2,3,\ldots \}$.}
  $$
  G^\infty = \prod _{n \in {\bf N}} G
  $$
  consider the infinite direct sum
  $$
  G^{(\infty )} = \bigoplus _{n \in {\bf N}} G
  $$
  formed by the elements $g=(g_n)_{n \in {\bf N}} \in G^\infty $ which are eventually trivial, that is, for which there exists $n_0$ such
that $g_n=1$, for all $n\geq n_0$.
  It is clear that $G^{(\infty )}$ is a normal subgroup of $G^\infty $.

\definition Given a group $G$, the \"{corona} of $G$ is the quotient group
  $$
  \corona = G^\infty /G^{(\infty )}.
  $$

  Consider the \"{left} and \"{right shift} endomorphisms of $G^\infty $
  $$
  \lambda ,\rho : G^\infty \to G^\infty
  $$
  given for every $\g =(\g _n)_{n \in {\bf N}} \in G^\infty $, by
  $$
  \lambda (\g )_n = \g _{n+1} \for n \in {\bf N},
  $$
  and
  $$
  \rho (\g )_n = \left \{ \matrix {1, & \hbox { if } n=0, \cr \g _{n-1}, & \hbox { if } n\geq 1.} \right .
  $$
  It is readily seen that $G^{(\infty )}$ is invariant under both $\lambda $ and $\rho $, so these pass to the quotient providing
endomorphisms
  $$
  \q \lambda , \q \rho : \corona \to \corona .
  \equationmark DefineQuotientAuto
  $$

For every $\g =(\g _n)_{n \in {\bf N}} \in G^\infty $, we have that
  $$
  \lambda (\rho (\g ))= \g
  \and
  \rho (\lambda (\g )) = (1,\g _2,\g _3,\ldots ) \equiv \g ,
  \equationmark ShiftIsAuto
  $$
  where we use ``$\equiv $'' to refer to the equivalence relation determined by the normal subgroup $G^{(\infty )}$.
  Therefore both $\q \lambda \q \rho $ and $\q \rho \q \lambda $ coincide with the identity, and hence $\q \lambda $ and $\q \rho $ are each other's
inverse.  In particular, they are both automorphisms of $\corona $.

Iterating $\q \rho $ therefore gives an action of ${\bf Z}$ on $\corona $.

\definition Given any countable discrete group $G$, the \"{lag group} associated to $G$ is the semi-direct product group
  $$
  \corona \ifundef {rtimes} \times \else \rtimes \fi _{\q \rho } {\bf Z}.
  $$

The reason we call this the ``lag group'' is that it will play a very important role in the next section, as the co-domain
for our \"{lag} function.

\section The tight groupoid of $\SGE $

\label TightGpdSectn We would now like to give a detailed description of the groupoid $\GpdGE $.  As already mentioned
this is the groupoid of germs for the natural action of $\SGE $ on the space of tight filters over the idempotent
semi-lattice $\EGE $ of $\SGE $.  See \cite [Section 4]{actions} for more details.

  By an \"{infinite path} in $E$ we shall mean any infinite sequence of the form
  $$
  \xi =\xi _1\xi _2\ldots ,
  $$
  where $\xi _i \in E^1$, and $\src (\xi _i) = \ran (\xi _{i+1})$, for all $i$.  The set of all infinite paths will be denoted by
$E^\infty $.  Given an infinite path
  $$
  \xi =\xi _1\xi _2\ldots \in E^\infty ,
  $$
  and an integer $n\geq 0$, denote by $\trunc \xi n$ the finite path of length $n$ given by
  $$
  \trunc \xi n = \left \{\matrix {\xi _1\xi _2\ldots \xi _n, & \hbox { if } n\geq 1, \cr \cr
  \ran (\xi _1), & \hbox { if } n=0.}\right .
  $$

\state Proposition \label ActionOnInfWords
  There is a unique action
  $$
  (g,\xi ) \in G\times E^\infty \mapsto g\xi \in E^\infty
  $$
  of $G$ on $E^\infty $ such that,
  $$
  \trunc {(g\xi )}n = g (\trunc \xi n),
  $$
  for every $g \in G$, $\xi \in E^\infty $, and $n \in {\bf N}$.

\Proof Left to the reader. \endProof

Recall from \ref{DefineEAlpha} that, for any finite path $\alpha \in E^*$, we denote by $\eproj _\alpha $ the idempotent element
$(\alpha ,1,\alpha )$ in $\EGE $.  Thus, given an infinite path $\xi \in E^\infty $, we may look at the subset
  $$
  {\cal F}_\xi = \{\eproj _{\xi |_n}: n \in {\bf N}\} \subseteq \EGE ,
  $$
  which turns out to be an ultra-filter \cite [Definition]{actions} over $\EGE $.  Denoting the set of all
ultra-filters over $\EGE $ by $\widehat \EGE _\infty $, as in \cite [Definition 12.8]{actions}, one may also show
  \cite [Proposition 19.11]{actions}
  that the correspondence
  $$
  \xi \in E^\infty \mapsto {\cal F}_\xi \in \widehat \EGE _\infty
  $$
  is bijective, and we will use it to identify
  $
  E^\infty
  $
  and
  $
  \widehat \EGE _\infty .
  $
  Furthermore, this correspondence may be proven to be a homeomorphism if $E^\infty $ is equipped with the product topology.

Since $E$ is finite, $E^\infty $ is compact by Tychonov's Theorem, and consequently so is $\widehat \EGE _\infty $.  Being the
closure of $\widehat \EGE _\infty $ within $\pilar {11pt}\widehat \EGE $ \cite [Theorem 12.9]{actions}, the space $\widehat
\EGE \tight $ formed by the tight filters therefore necessarily coincides with $\widehat \EGE _\infty $.

Identifying $\widehat \EGE \tight $ with $E^\infty $, as above, we may transfer the canonical action of $\SGE $ from the
former to the latter resulting in the following: to each element
  $(\alpha ,g,\beta ) \in \SGE $,
  we associate the partial homeomorphism of $E^\infty $ whose domain is the \"{cylinder}
  $$
  \cyl \beta := \{\eta \in E^\infty : \eta =\beta \xi , \hbox { for some } \xi \in E^\infty \},
  \equationmark DefineCylinder
  $$
  and which sends each $\eta =\beta \xi \in \cyl \beta $ to $\alpha g\xi $, where the meaning of ``$g\xi $'' is as in
  \redundantref {Proposition}{ActionOnInfWords}.

  As before we will not use any special symbol to indicate this action, using module notation instead:
  $$
  (\alpha ,g,\beta )\eta = \alpha g\xi
  \for (\alpha ,g,\beta )\in \SGE \for \eta =\beta \xi \in \cyl \beta .
  \equationmark ActionOfSGE
  $$

Before we proceed let us at least check that $\alpha g\xi $ is in fact an element of $E^\infty $, which is to say that
  $
  \src (\alpha ) = \ran (g\xi ).
  $
  Firstly, for every element $(\alpha ,g,\beta ) \in \SGE $, we have that $\src (\alpha )=g\src (\beta )$.  Secondly, if $\eta =\beta \xi \in E^\infty $, then $\src
(\beta ) = \ran (\xi )$.  Therefore
  $$
  \ran (g\xi ) = g\ran (\xi ) = g\src (\beta ) = \src (\alpha ).
  $$

  This leads to a first, more or less concrete description of $\GpdGE $.

\state Proposition \label FirstDescrGpd
  Under \ref{StandingHyp}, one has that $\GpdGE $ is isomorphic to the groupoid of germs for the above action of
$\SGE $ on $E^\infty $.

Our aim is nevertheless a much more precise description of this groupoid.  Recall from \cite [Definition 4.6]{actions}
that the germ of an element $s \in \SGE $ at a point $\xi $ in the domain of $s$ is denoted by
  $[s,\xi ]$.  If $s=(\alpha ,g,\beta )$, this would lead to the somewhat awkward notation
  $[(\alpha ,g,\beta ),\xi ]$, which from now on will instead be written as
  $$
  \germ \alpha g\beta \xi .
  $$

Thus the groupoid $\GpdGE $, consisting of all germs for the action of $\SGE $ on $E^\infty $, is given by
  $$
  \GpdGE =
  \Big \{\germ \alpha g\beta \xi : (\alpha ,g,\beta ) \in \SGE ,\ \xi \in \cyl \beta \Big \}.
  \equationmark FirstModel
  $$

Let us now prove a useful criterion for equality of germs.

\state Proposition \label EqualGerms
  Suppose that $(\Data )$ is pseudo free and let us be given elements $(\alpha _1,g_1,\beta _1)$ and $(\alpha _2,g_2,\beta _2)$ in $\SGE $, with
$|\beta _1|\leq |\beta _2|$, as well as infinite paths $\eta _1$ in $\cyl {\beta _1}$, and $\eta _2$ in $\cyl {\beta _2}$.  Then
  $$
  \germ {\alpha _1}{g_1}{\beta _1}{\eta _1} = \germ {\alpha _2}{g_2}{\beta _2}{\eta _2}
  $$
  if and only if there is a finite path $\gamma $ and an infinite path $\xi $, such that
  \Zitem $\alpha _2 = \alpha _1 g_1\gamma ,$
  \zitem $g_2 = \varphi (g_1,\gamma ),$
  \zitem $\beta _2 = \beta _1\gamma ,$
  \zitem $\eta _1=\eta _2=\beta _1\gamma \xi .$

\Proof Assuming that the germs are equal, we have by \cite [Definition 4.6]{actions} that
  $$
  \eta _1=\eta _2 =: \eta ,
  $$
  and there is an idempotent $(\delta ,1,\delta ) \in \EGE $, such that $\eta \in \cyl \delta $, and
  $$
  (\alpha _1,g_1,\beta _1) (\delta ,1,\delta ) = (\alpha _2,g_2,\beta _2) (\delta ,1,\delta ).
  \equationmark SameGerm
  $$

  It follows that $\eta =\delta \zeta $, for some $\zeta \in E^\infty $.  Upon replacing $\delta $ by a longer prefix of $\eta $, we may assume that $|\delta |$ is
as large as we want.  Furthermore the element of $\SGE $ represented by the two sides of
  \ref{SameGerm}
  is evidently nonzero because the partial homeomorphism associated to it under our action has $\eta $ in its domain.  So,
focusing on \ref{inverseSemigroup}, we see that $\beta _1$ and $\delta $ are comparable, and so are $\beta _2$ and $\delta $.

Assuming that $|\delta |$ exceeds both $|\beta _1|$ and $|\beta _2|$, we may then write
  $\delta = \beta _1\varepsilon _1 = \beta _2\varepsilon _2$,
  for suitable $\varepsilon _1$ and $\varepsilon _2$ in $E^*$.  But since $|\beta _1|\leq |\beta _2|$, this in turn implies that $\beta _2=\beta _1\gamma $, for some $\gamma \in E^*$, hence
proving (iii).
  Therefore
  $
  \delta =\beta _1\gamma \varepsilon _2,
  $
  so
  $$
  \eta = \delta \zeta = \beta _1\gamma \varepsilon _2\zeta ,
  $$
  and (iv) follows once we choose $\xi = \varepsilon _2\zeta $.  Moreover, equation \ref{SameGerm} reads
  $$
  (\alpha _1,g_1,\beta _1) (\beta _1\gamma \varepsilon _2,1,\beta _1\gamma \varepsilon _2) = (\alpha _2,g_2,\beta _1\gamma ) (\beta _1\gamma \varepsilon _2,1,\beta _1\gamma \varepsilon _2).
  $$

  Computing the products according to \ref{inverseSemigroup}, we get
  $$
  \big(\alpha _1g_1(\gamma \varepsilon _2),\ \varphi (g_1,\gamma \varepsilon _2),\ \beta _1\gamma \varepsilon _2\big) = \big(\alpha _2g_2\varepsilon _2,\ \varphi (g_2,\varepsilon _2),\ \beta _1\gamma \varepsilon _2\big),
  $$
  from where we obtain
  $$
  \alpha _2g_2\varepsilon _2 = \alpha _1g_1(\gamma \varepsilon _2) = \alpha _1(g_1\gamma )\varphi (g_1,\gamma )\varepsilon _2,
  \equationmark EqOne
  $$
  and
  $$
  \varphi (g_2,\varepsilon _2) = \varphi (g_1,\gamma \varepsilon _2) = \varphi \big(\varphi (g_1,\gamma ),\varepsilon _2\big).
  \equationmark EqTwo
  $$

  Since $|g_2\varepsilon _2| = |\varepsilon _2| = |\varphi (g_1,\gamma )\varepsilon _2|$, we deduce from \ref{EqOne} that
  $$
  g_2\varepsilon _2 = \varphi (g_1,\gamma )\varepsilon _2,
  \equationmark DeduceOne
  $$
  and hence also that
  $$
  \alpha _2 = \alpha _1g_1\gamma ,
  $$
  proving (i).  In view of \ref{EqTwo} and \ref{DeduceOne}, point (ii) follows from
  \redundantref {Proposition}{VarStar}.

Conversely, assume (i--iv) and let us prove equality of the above germs.  Setting $\delta = \beta _1\gamma $, we have by (iv) that
  $$
  \eta := \eta _1 = \eta _2 \in \cyl \delta ,
  $$
  so it suffices to verify \ref{SameGerm}, which the reader could do without any difficulty.  \endProof

Proposition \ref{EqualGerms} then says that the typical situation in which an equality of germs takes place is
  $$
  \germ {\alpha }{g}{\beta }{\beta \gamma \xi } = \germ {\alpha g\gamma }{\varphi (g,\gamma )}{\beta \gamma }{\beta \gamma \xi }.
  $$

Our next two results are designed to offer convenient representatives of germs.

\state Proposition \label AnyLength
  Given any germ $u$, there exists an integer $n_0$, such that for every $n\geq n_0$,
  \Zitem there is a representation of $u$ of the form $u=\germ {\alpha _1} {g_1} {\beta _1} {\beta _1\xi _1}$, with $|\alpha _1|=n$.
  \zitem there is a representation of $u$ of the form $u=\germ {\alpha _2} {g_2} {\beta _2} {\beta _2\xi _2}$, with $|\beta _2|=n$.

\Proof
  Write $u=\germ \alpha g \beta \eta $, and choose any $n_0\geq \max \{|\alpha |,|\beta |\}$.  Then, for every $n\geq n_0$ we may write $\eta =\beta \gamma \xi $, with $\gamma \in E^*$,
$\xi \in E^\infty $, and such that $|\gamma |=n-|\alpha |$ (resp.~$|\gamma |=n-|\beta |$).  Therefore
  $$
  u = \germ \alpha g \beta {\beta \gamma \xi } = \germ {\alpha g\gamma } {\varphi (g,\gamma )} {\beta \gamma } {\beta \gamma \xi },
  $$
  and we have
  $
  |\alpha g\gamma | = |\alpha | + |g\gamma | = |\alpha | + |\gamma | = n
  $
  (resp.~
  $
  |\beta \gamma | = |\beta | + |\gamma | = n
  $).
  \endProof

\state Corollary \label MultPairs
  Given $u_1$ and $u_2$ in $\GpdGE $, such that $(u_1,u_2) \in \GpdGE ^{(2)}$, that is, such that the multiplication $u_1u_2$ is allowed
or, equivalently, such that $\src (u_1)=\ran (u_2)$, then there are representations of $u_1$ and $u_2$ of the form
  $$
  u_1 = \germ {\alpha _1} {g_1} {\alpha _2} {\alpha _2g_2\xi } \and
  u_2 = \germ {\alpha _2} {g_2} {\beta } {\beta \xi },
  $$
  and in this case
  $$
  u_1u_2 = \germ {\alpha _1} {g_1g_2} {\beta } {\beta \xi }.
  $$

\Proof
  Using
  \redundantref {Proposition}{AnyLength}, write
  $$
  u_i =
  \germ {\alpha _i} {g_i} {\beta _i} {\beta _i\xi _i},
  $$
  with $|\beta _1|=|\alpha _2|$.
  By virtue of $(u_1,u_2)$ lying in $\GpdGE ^{(2)}$, we have that
  $$
  \beta _1\xi _1 = (\alpha _2,g_2,\beta _2)(\beta _2\xi _2) = \alpha _2g_2\xi _2,
  $$
  so in fact $\beta _1 = \alpha _2$, and $\xi _1 = g_2\xi _2$.  Then
  $$
  u_1 =
  \germ {\alpha _1} {g_1} {\beta _1} {\beta _1\xi _1} =
  \germ {\alpha _1} {g_1} {\alpha _2} {\alpha _2g_2\xi _2},
  $$
  and it suffices to put $\xi =\xi _2$, and $\beta = \beta _2$.

With respect to the last assertion we have that $u_1u_2 = [s;\beta \xi ]$, where $s$ is the element of $\SGE $ given by
  $$
  s =
  (\alpha _1, g_1, \alpha _2) (\alpha _2, g_2,\beta ) \={EasyMult}
  (\alpha _1, g_1g_2, \beta ),
  $$
  concluding the proof.
  \endProof

Having extended the action of $G$ to the set of infinite paths in Proposition \ref{ActionOnInfWords}, one may ask
whether it is possible to do the same for the cocycle $\varphi $.  The following is an attempt at this which however produces a
map taking values in the infinite product $G^\infty $, rather than in $G$.

  \definition We will denote by $\Phi $, the map
  $$
  \Phi : G\times E^\infty \to G^\infty
  $$
  defined by the rule
  $$
  \Phi (g,\xi )_n=\varphi (g,\trunc \xi {n-1}),
  $$
  for $g\in G$, $\xi \in E^\infty $, and $n\geq 1$.

Recall that we are indexing the elements of $G^\infty $ on the set $\{1,2,3,\ldots \}$, so the first coordinate of $\Phi (g,\xi )$ is
  $$
  \Phi (g,\xi )_1 =
  \varphi (g,\trunc \xi 0) =
  \varphi \big(g,\ran (\xi )\big) \={Equacoes.ii} g.
  $$

We wish to view $\Phi $ as some sort of cocycle but, unfortunately, property \ref{Equacoes.\LastItem } does not hold
quite as stated.  On the fortunate side, a suitable modification of this relation, involving the left shift endomorphism
$\lambda $ of $G^\infty $, is satisfied:

\state Proposition \label LastWithShift
  Let $\alpha $ be a finite path and let $\xi $ be an infinite path such that $\src (\alpha )=\ran (\xi )$.  Then, for every $g$ in $G$,
one has that
  $$
  \Phi \big(\varphi (g,\alpha ),\xi \big) =
  \lambda ^{|\alpha |}\big(\Phi \big(g,\alpha \xi )\big).
  $$

\Proof For all $n\geq 1$, we have
  $$
  \Phi \big(\varphi (g,\alpha ),\xi \big)_n =
  \varphi \big(\varphi (g,\alpha ),\trunc \xi {n-1}\big) =
  \varphi \big(g,\alpha (\trunc \xi {n-1})\big) \$=
  \varphi \big(g,\trunc {(\alpha \xi )}{n-1+|\alpha |}\big) =
  \lambda ^{|\alpha |}\big(\Phi (g,\alpha \xi )\big)_n.
  \endProof

Another reason to think of $\Phi $ as a cocycle is as follows:

\state Proposition \label CocycleIdForCapPhi
  For every $\xi \in E^\infty $, and every $g,h \in G$, we have that
  $$
  \Phi (gh,\xi ) = \Phi \big(g,h\xi \big)\Phi (h,\xi ).
  $$

\Proof We have for all $n \in {\bf N}$, that
  $$
  \Phi (gh,\xi )_n =
  \varphi (gh,\xi |_{n-1}) \={Equacoes.b}
  \varphi \big(g,h(\xi |_{n-1})\big)\varphi (h,\xi |_{n-1}) \={ActionOnInfWords} $$$$=
  \varphi \big(g,(h\xi )|_{n-1}\big)\varphi (h,\xi |_{n-1}) =
  \Phi \big(g,h\xi \big)_n\Phi (h,\xi )_n.
  \endProof

The following elementary fact might perhaps justify the choice of ``$n-1$'' in the definition of $\Phi $.

\state Proposition \label ActionCocycle
  Given $g \in G$, and $\xi \in E^\infty $, one has that
  $$
  (g\xi )_n = \Phi (g,\xi )_n\,\xi _n.
  $$

  \Proof
  By
  \redundantref {Proposition}{ActionOnInfWords} we have that $(g\xi )|_n = g(\xi |_n)$, so the $n^{th}$ coordinate of $g\xi $ is
also the $n^{th}$ coordinate of $g(\xi |_n)$.  In addition we have that
  $$
  g(\xi |_n) =
  g(\xi |_{n-1}\xi _n) \={Equacoes.\MainConcat }
  g(\xi |_{n-1})\varphi (g,\xi |_{n-1})\xi _n,
  $$
  so
  $$
  (g\xi )_n = \varphi (g,\xi |_{n-1})\xi _n = \Phi (g,\xi )_n\,\xi _n.
  \endProof

We now wish to define a homomorphism (sometimes also called a one-cocycle) from $\GpdGE $ to the lag group $\corona
\ifundef {rtimes} \times \else \rtimes \fi _\rho {\bf Z}$, by means of the rule
  $$
  \germ \alpha g\beta {\beta \xi } \mapsto \Big (\rho ^{|\alpha |}\big(\Phi (g,\xi )\big),|\alpha |-|\beta |\Big ).
  $$

  As it is often the case for maps defined on groupoid of germs, the above tentative definition uses a representative of
the germ, so some work is necessary to prove that the definition does not depend on the choice of representative.  The
technical part of this task is the content of our next result.

\state Lemma \label LagWellDefined
  Suppose that $(\Data )$ is pseudo free.  For each $i=1,2$, let us be given $(\alpha _i,g_i,\beta _i)$ in $\SGE $, as well as
$\eta _i=\beta _i\xi _i \in \cyl {\beta _i}$.  If
  $$
  \germ {\alpha _1}{g_1}{\beta _1}{\eta _1} = \germ {\alpha _2}{g_2}{\beta _2}{\eta _2},
  $$
  then
  $$
  \rho ^{|\alpha _1|}\big(\Phi (g_1,\xi _1)\big) \equiv \rho ^{|\alpha _2|}\big(\Phi (g_2,\xi _2)\big)
  $$
  modulo $G^{(\infty )}$.

\Proof
  Assuming without loss of generality that $|\beta _1|\leq |\beta _2|$, we may use
  \redundantref {Proposition}{EqualGerms} to write
  $$
  \alpha _2 = \alpha _1 g_1\gamma , \quad g_2 = \varphi (g_1,\gamma ), \quad \beta _2 = \beta _1\gamma \and \eta  _1=\eta _2=\beta _1\gamma \xi ,
  $$
  for suitable $\gamma \in E^*$ and $\xi \in E^\infty $.  Then necessarily $\xi _1 = \gamma \xi $, and $\xi _2=\xi $, and
  $$
  \rho ^{|\alpha _2|}\big(\Phi (g_2,\xi _2)\big) =
  \rho ^{|\alpha _1|+|\gamma |}\Big ( \Phi \big(\varphi (g_1,\gamma ),\xi \big)\Big ) \={LastWithShift} $$ $$=
  \rho ^{|\alpha _1|}\rho ^{|\gamma |}\lambda ^{|\gamma |}\big(\Phi \big(g_1,\gamma \xi )\big) \explain \equiv {ShiftIsAuto}
  \rho ^{|\alpha _1|}\big(\Phi \big(g_1,\xi _1)\big).
  \endProof

\fix Due to our reliance on Proposition \ref{EqualGerms} and Lemma \ref{LagWellDefined}, from now on and until
the end of this section we will assume, in addition to \ref{StandingHyp}, that $(\Data )$ is pseudo free.

\bigskip If $\g $ is in $G^\infty $, we will denote by $\q \g $ its class in the quotient group $\corona $.  Likewise we will
denote by $\q \Phi $ the composition of $\Phi $ with the quotient map from $G^\infty $ to $\corona $.

  \beginpicture
  \setcoordinatesystem units <0.0040truecm, -0.0040truecm> point at 3000 0
  \setplotarea x from -1500 to 1000, y from -200 to 500
  \put {$G\times E^\infty \longrightarrow \ G^\infty \longrightarrow \corona $} at 510 000
  \put {$\Phi $} at 440 -90
  \setquadratic
  \plot 100 120 525 250 950 120 /
  \arrow <0.11cm> [0.3,1.2] from 950 120 to 960 113
  \put {$\q \Phi $} at 525 350
  \endpicture

  \bigskip
  It then follows from
  \redundantref {Lemma}{LagWellDefined} that the correspondence
  $$
  \germ \alpha g\beta {\beta \xi } \in \GpdGE \mapsto \q \rho ^{|\alpha |}\big(\q \Phi (g,\xi )\big) \in \corona
  $$
  is a well defined map.  This is an important part of the one-cocycle we are about to introduce.

\state Proposition \label DefineLag
  The correspondence
  $$
  \lag : \germ \alpha g\beta {\beta \xi } \mapsto \Big (\q \rho ^{|\alpha |}\big(\q \Phi (g,\xi )\big),|\alpha |-|\beta |\Big )
  $$
  gives a well defined map
  $$
  \lag : \GpdGE \to \corona \ifundef {rtimes} \times \else \rtimes \fi _\rho {\bf Z},
  $$
  which is moreover a one-cocycle.  From now on $\lag $ will be called the \"{lag function}.

  \Proof
  By the discussion above we have that the first coordinate of the above pair is well defined.  On the other hand, in
the context of Proposition \ref{EqualGerms} one easily sees that $|\alpha _1|-|\beta _1| = |\alpha _2|-|\beta _2|$, so the second coordinate
is also well defined.

In order to show that $\lag $ is multiplicative, pick $(u_1,u_2) \in \GpdGE ^{(2)}$.  We may then use
  \redundantref {Corollary}{MultPairs} to write
  $$
  u_1 = \germ {\alpha _1} {g_1} {\alpha _2} {\alpha _2g_2\xi }
  \and
  u_2 = \germ {\alpha _2} {g_2} {\beta } {\beta \xi }.
  $$
  So
  $$
  \lag (u_1)\lag (u_2) =
  \Big (\rho ^{|\alpha _1|}\big(\Phi (g_1,g_2\xi )\big),\ |\alpha _1|-|\alpha _2|\Big )
  \Big (\rho ^{|\alpha _2|}\big( \Phi (g_2,\xi )\big),\ |\alpha _2|-|\beta |\Big )
  \$=
  \Big (\rho ^{|\alpha _1|}\big(\Phi (g_1,g_2\xi )\big) \ \rho ^{|\alpha _1|}\big( \Phi (g_2,\xi )\big),\ \ |\alpha _1|-|\alpha _2|+|\alpha _2|-|\beta |\Big ) \$=
  \Big (\rho ^{|\alpha _1|}\Big (\Phi (g_1,g_2\xi ) \Phi (g_2,\xi )\Big ),\ |\alpha _1|-|\beta |\Big ) \={CocycleIdForCapPhi}
  \Big (\rho ^{|\alpha _1|}\big(\Phi (g_1g_2,\xi )\big),\ |\alpha _1|-|\beta |\Big ) \$=
  \lag \big(\germ {\alpha _1} {g_1g_2} {\beta } {\beta \xi }\big) \={MultPairs}
  \lag (u_1u_2).
  \endProof

The main relevance of this one-cocycle is that, together with the domain and range maps, it uniquely describes the
elements of $\GpdGE $, as we will now show.

\state Proposition \label FInjective
  Given
  $
  u_1,u_2 \in \GpdGE ,
  $
  one has that
  $$
  \left .\matrix {
    \src (u_1)=\src (u_2) \cr
    \ran (u_1)=\ran (u_2) \pilar {13pt} \cr
    \lag (u_1)=\lag (u_2) \pilar {14pt} }
    \right \} \IMPLY u_1=u_2.
  $$

\Proof Using
  \redundantref {Proposition}{AnyLength}, write $u_i = \germ {\alpha _i}{g_i}{\beta _i}{\beta _i\xi _i}$, for $i=1,2$, with $|\beta _1|=|\beta _2|$.
Since
  $$
  \beta _1\xi _1 = \src (u_1) = \src (u_2) = \beta _2\xi _2,
  $$
  we conclude that
  $
  \beta _1=\beta _2,
  $
  and
  $$
  \xi _1=\xi _2=:\xi .
  $$

By focusing on the second coordinate of $\lag (u_i)$, we see that $|\alpha _1|-|\beta _1|=|\alpha _2|-|\beta _2|$, and hence $|\alpha _1|=|\alpha _2|$.
Moreover, since
  $$
  \alpha _1g_1\xi = \alpha _1g_1\xi _1 = \ran (u_1) = \ran (u_2) = \alpha _2g_2\xi _2 = \alpha _2g_2\xi ,
  $$
  we see that $\alpha _1=\alpha _2$, and
  $$
  g_1\xi = g_2\xi .
  \equationmark EqualGXi
  $$

The fact that $\lag (u_1) = \lag (u_2)$ also implies that
  $$
  \q \rho ^{|\alpha _1|}\big(\q \Phi (g_1,\xi )\big) =
  \q \rho ^{|\alpha _2|}\big(\q \Phi (g_2,\xi )\big),
  $$
  and since $\alpha _1=\alpha _2$, we conclude that
  $
  \q \Phi (g_1,\xi ) = \q \Phi (g_2,\xi ),
  $
  and hence that there exists an integer $n_0$ such that
  $$
  \varphi (g_1,\xi |_n) = \varphi (g_2,\xi |_n)
  \for n\geq n_0.
  $$
  By \ref{EqualGXi} we also have that
  $
  g_1(\xi |_n) = g_2(\xi |_n),
  $
  so \ref{VarStar} gives $g_1=g_2$, whence $u_1=u_2$.
  \endProof

As a consequence of the above result we see that the map
  $$
  F:\GpdGE \to E^\infty \times (\corona \ifundef {rtimes} \times \else \rtimes \fi _{\q \rho } {\bf Z}) \times E^\infty
  $$
  defined by the rule
  $$
  F(u) = \big(\ran (u),\lag (u),\src (u)\big),
  \equationmark DefineF
  $$
  is one-to-one.

Observe that the co-domain of $F$ has a natural groupoid structure, being the cartesian product of the lag group
$\corona \ifundef {rtimes} \times \else \rtimes \fi _{\q \rho } {\bf Z}$ by the graph of the transitive equivalence relation on $E^\infty $.

Putting together
  \redundantref {Proposition}{DefineLag} and
  \redundantref {Proposition}{FInjective} we may now easily prove:

\state Corollary $F$ is a groupoid homomorphism (functor), hence establishing an isomorphism from $\GpdGE $ to the range
of $F$.

The range of $F$ is then the concrete model of $\GpdGE $ we are after.  But, before giving a detailed description of it,
let us make a remark concerning notation: since the co-domain of $F$ is a mixture of cartesian and semi-direct products,
the standard notation for its elements would be something like $\big(\eta ,(u,p),\zeta \big)$, for $\eta ,\zeta \in  E^\infty $, $u \in \corona $, and $p \in {\bf Z}$.
As part of our effort to avoid heavy notation we will instead denote such an element by
  $$
  \big(\eta ;u,p;\zeta \big).
  $$

\state Proposition
  The range of $F$ is precisely the subset of $E^\infty \times (\corona \ifundef {rtimes} \times \else \rtimes \fi _{\q \rho } {\bf Z}) \times E^\infty $, formed by the elements
  $
  (\eta ;\q \g ,p-q;\zeta ),
  $
  where $\eta ,\zeta \in E^\infty $, $\g \in G^\infty $, and $p,q \in {\bf N}$, are such that, for all $n\geq 1$,
  \Zitem $\g _{n+p+ 1} = \varphi (\g _{n+p},\zeta _{n+q})$,
  \zitem $\eta _{n+p} = \g _{n+p}\zeta _{n+q}$.

\Proof
  Pick a general element $\germ \alpha g\beta {\beta \xi } \in \GpdGE $ and, recalling that
  $$
  F(\germ \alpha g\beta {\beta \xi }) = \Big (\alpha g\xi ;\ \q \rho ^{|\alpha |}\big(\q \Phi (g,\xi )\big),\ |\alpha |-|\beta |;\ \beta \xi \Big ),
  \equationmark GeneralRangeF
  $$
  let
  $
  \eta =\alpha g\xi , \ \g = \rho ^{|\alpha |}\big(\Phi (g,\xi )\big), \ p=|\alpha |,\ q=|\beta |, \hbox { and } \zeta =\beta \xi ,
  $
  so that the element depicted in \ref{GeneralRangeF} becomes $(\eta ;\q \g ,p-q;\zeta )$, and we must now verify (i) and
(ii).
   For all $n\geq 1$, one has that
  $$
  \g _{n+|\alpha |} = \Phi (g,\xi )_n = \varphi (g,\xi |_{n-1}),
  $$
  so
  $$
  \eta _{n+p} =
  (\alpha g\xi )_{n+|\alpha |} =
  (g\xi )_n \={ActionCocycle}
  \varphi (g,\xi |_{n-1})\xi _n =
  \g _{n+|\alpha |}(\beta \xi )_{n+|\beta |} =
  \g _{n+p}\zeta _{n+q},
  $$
  proving (ii).  Also,
  $$
  \g _{n+p+ 1} =
  \g _{n+|\alpha |+ 1} =
  \varphi \big(g,\xi |_{n}\big) =
  \varphi (g,\xi |_{n-1}\xi _{n}) =
  \varphi \big(\varphi (g,\xi |_{n-1}),\xi _{n}\big) \$=
  \varphi \big(\g _{n+|\alpha |},(\beta \xi )_{n+|\beta |}\big) =
  \varphi \big(\g _{n+p},\zeta _{n+q}\big),
  $$
  proving (i) and hence showing that the range of $F$ is a subset of the set described in the statement.

Conversely, pick $\eta ,\zeta \in E^\infty $, $\g \in G^\infty $, and $p,q \in {\bf N}$ satisfying (i) and (ii), and let us show that the element
  $(\eta ;\q \g ,p-q;\zeta )$ lies in the range of $F$.  Let
  $$
  g = \g _{p+1},\quad
  \alpha = \eta |_p \and \beta =\zeta |_q,
  $$
  so $\zeta =\beta \xi $ for a unique $\xi \in E^\infty $.  We then claim that $\germ \alpha g\beta {\beta \xi }$ lies in $\GpdGE $.  In order to see this notice
that
  $$
  g\src (\beta ) =
  g\src (\zeta _q) =
  g\ran (\zeta _{q+1}) =
  \ran (\g _{p+1}\zeta _{q+1}) \explica{=}{(ii)}
  \ran (\eta _{p+1}) =
  \src (\eta _p) =
  \src (\alpha ),
  $$
  so $(\alpha ,g,\beta ) \in \SGE $, and therefore $\germ \alpha g\beta {\beta \xi }$ is indeed a member of $\GpdGE $.
  The proof will then be concluded once we show that
  $$
  F(\germ \alpha g\beta {\beta \xi }) = (\eta ;\q \g ,p-q;\zeta ),
  $$
  which in turn is equivalent to showing that
  \iaitem
  \aitem $\alpha g\xi = \eta ,$
  \aitem $\q \rho ^{|\alpha |}\big(\q \Phi (g,\xi )\big) = \q \g ,$
  \aitem $|\alpha |-|\beta | = p-q,$
  \aitem $\beta \xi = \zeta $.

\medskip Before proving these points we will show that
  $$
  \varphi (\g _{p+1},\xi |_n) = \g _{n+p+1} \for n\geq 0.
  \eqno {(\dagger )}
  $$
  This is obvious for $n=0$.  Assuming that $n\geq 1$ and using induction, we have
  $$
  \varphi (\g _{p+1},\xi |_n) =
  \varphi (\g _{p+1},\xi |_{n-1}\xi _n) =
  \varphi \big(\varphi (\g _{p+1},\xi |_{n-1}),\xi _n\big) \$=
  \varphi (\g _{n+p},\zeta _{n+q}) \explica={(i)} \g _{n+p+1},
  $$
  verifying $(\dagger )$.

Addressing (a) we have to prove that $(\alpha g\xi )_k = \eta _k$, for all $k\geq 1$, but given that $\alpha $ is defined to be $\eta |_p$, this is
trivially true for $k\leq p$.  On the other hand, for $k=n+p$, with $n\geq 1$, we have
  $$
  (\alpha g\xi )_k =
  (\alpha g\xi )_{n+p} =
  (g\xi )_n \={ActionCocycle}
  \varphi (g,\xi |_{n-1})\xi _n \$=
  \varphi (\g _{p+1},\xi |_{n-1})\xi _n \explica={($\scriptstyle \dagger $)}
  \g _{n+p}\zeta _{n+q} \explica={(ii)}
  \eta _{n+p} =
  \eta _k,
  $$
  proving (a).  Focusing on (b) we have for all $n\geq 1$ that
  $$
  \rho ^{|\alpha |}\big(\Phi (g,\xi )\big)_{p+n} =
  \Phi (g,\xi )_n =
  \varphi (\g _{p+1},\xi |_{n-1}) \explica={($\scriptstyle \dagger $)}
  \g _{n+p},
  $$
  proving that $\rho ^{|\alpha |}\big(\Phi (g,\xi )\big) \equiv \g $, modulo $G^{(\infty )}$, hence taking care of (b).  The last two points, namely (c) and
(d) are trivial and so the proof is concluded.
  \endProof

As an immediate consequence we get a very precise description of the algebraic structure of $\GpdGE $:

\state Theorem \label ConcreteModel
  Suppose that $(\Data )$ satisfies the conditions of \ref{StandingHyp} and is moreover pseudo free.  Then $\GpdGE
$ is isomorphic to the sub-groupoid of \/
  $E^\infty \times (\corona \ifundef {rtimes} \times \else \rtimes \fi _{\q \rho } {\bf Z}) \times E^\infty $ given by
  $$
  \G _{G,E} = \left \{
  \matrix {
  (\eta ;\q \g ,p-q;\zeta ) \ \in & E^\infty \times (\corona \ifundef {rtimes} \times \else \rtimes \fi _{\q \rho } {\bf Z}) \times E^\infty : \hfill \cr
  & \g \in G^\infty ,\ p,q \in {\bf N}, \hfill \pilar {12pt} \cr
  & \g _{n+p+1} = \varphi (\g _{n+p},\zeta _{n+q}), \hfill \pilar {12pt} \cr
  & \eta _{n+p} = \g _{n+p}\zeta _{n+q}, \hbox { for all }n\geq 1 \pilar {15pt}
  }
  \right \}.
  $$

\medskip Recall from \cite {KPRR} that the C*-algebra of every graph is a groupoid C*-algebra for a certain groupoid
constructed from the graph, and informally called the groupoid for the \"{tail equivalence with lag}.

Viewed through the above perspective, our groupoid may also deserve such a denomination, except that the lag is not just
an integer as in \cite {KPRR}, but an element of the lag group $\corona \ifundef {rtimes} \times \else \rtimes \fi _\rho  {\bf Z}$ precisely described by the lag function
$\lag $ introduced in Proposition \ref{DefineLag}.

\section The topology of $\GpdGE $

It is now time we look at the topological aspects of $\GpdGE $.  In fact what we will do is simply transfer the topology
of $\GpdGE $ over to $\G _{G,E} $ via $F$.  Not surprisingly $F$ will turn out to be an isomorphism of topological groupoids.

Recall from \cite [Proposition 4.14]{actions} that, if $\S $ is an inverse semigroup acting on a locally compact
Hausdorff topological space $X$, then the corresponding groupoid of germs, say $\G $, is topologized by means of the
basis consisting of sets of the form
  $$
  \Theta (s,U),
  $$
  where $s \in \S $, and $U$ is an open subset of $X$, contained in the domain of the partial homeomorphism attached to
$s$ by the given action.  Each $\Theta (s,U)$ is in turn defined by
  $$
  \Theta (s,U) = \Big \{[s,x]\in \G : x\in U\Big \}.
  \equationmark BasisOpenSets
  $$
  See \cite [4.12]{actions} for more details.

If we restrict the choice of the $U$'s above to a predefined basis of open sets of $X$, e.g.~the collection of all
cylinders in $E^\infty $ in the present case, we evidently get the same topology on the groupoid of germs.  Therefore,
referring to the model of $\GpdGE $ presented in \ref{FirstModel}, we see that a basis for its topology consists of
the sets of the form
  $$
  \Theta (\alpha ,g,\beta ;\gamma ) := \Big \{\germ \alpha g\beta \xi \in \GpdGE : \xi \in \cyl \gamma \Big \},
  \equationmark Slice
  $$
  where $(\alpha ,g,\beta ) \in \SGE $, and $\gamma \in E^*$.
  We may clearly suppose that $|\gamma | \geq |\beta |$ and, since $\Theta (\alpha ,g,\beta ;\gamma ) = \ifundef {varnothing} \emptyset \else \varnothing \fi
$, unless $\beta $ is a prefix of $\gamma $, we may also assume that $\gamma =\beta \varepsilon $, for some $\varepsilon \in E^*$.

In this case, given any $\germ \alpha g\beta \xi \in \Theta (\alpha ,g,\beta ;\gamma )$, notice that $\xi \in \cyl \gamma $, and
  $$
  (\alpha ,g,\beta ) (\gamma ,1,\gamma ) = \big(\alpha g\varepsilon ,\varphi (g,\varepsilon ),\gamma \big),
  $$
  from where one concludes that
  $$
  \germ \alpha g\beta \xi = \germ {\alpha g\varepsilon }{\varphi (g,\varepsilon )}\gamma \xi ,
  $$
  for all $\xi \in \cyl \gamma $, and hence also that
  $$
  \Theta (\alpha ,g,\beta ;\gamma ) = \Theta \big(\alpha g\varepsilon ,\varphi (g,\varepsilon ),\gamma ; \gamma \big).
  $$

  This shows that any set of the form \ref{Slice} coincides with another such set for which $\beta =\gamma $.  We may
therefore do away with this repetition and redefine
  $$
  \Theta (\alpha ,g,\beta ) := \Big \{\germ \alpha g\beta \xi \in \GpdGE : \xi \in \cyl \beta \Big \}.
  \equationmark ReSlice
  $$

We have therefore shown:

\state Proposition
  The collection of all sets of the form $\Theta (\alpha ,g,\beta )$, where $(\alpha ,g,\beta )$ range in $\SGE $, is a basis for the topology of
$\GpdGE $.

We may now give a precise description of the topology of $\GpdGE $, once it is viewed from the alternative point of view
of Theorem \ref{ConcreteModel}:

\state Proposition For each $(\alpha ,g,\beta )$ in $\SGE $, the image of\/ $\Theta (\alpha ,g,\beta )$ under $F$ coincides with the set
  $$
  \Omega (\alpha ,g,\beta ) :=
  \left \{
  \matrix {
  (\eta ;\q \g ,k;\zeta ) \in \G _{G,E} :
  &
    \eta \in \cyl \alpha ,\
    \g \in G^\infty ,\
    k = |\alpha |-|\beta |,\
    \zeta \in \cyl \beta , \hfill \cr
  & \g _{1+|\alpha |}=g, \hfill \pilar {12pt} \cr
  & \g _{n+|\alpha |+1} = \varphi (\g _{n+|\alpha |},\zeta _{n+|\beta |}), \hfill \pilar {11pt} \cr
  & \eta _{n+|\alpha |} = \g _{n+|\alpha |}\zeta _{n+|\beta |}, \hbox { for all } n\geq 1 \hfill \pilar {13pt}
  }
  \right \},
  $$
  and hence the collection of all such sets
  form the basis for a topology on $\G _{G,E} $, with respect to which the latter is isomorphic to $\GpdGE $ as topological
groupoids.

  \Proof Left for the reader. \endProof

We may now summarize the main results obtained so far:

\state Theorem \label GroupoidModel
  Suppose that $(\Data )$ satisfies the conditions of \ref{StandingHyp} and is moreover pseudo free.  Then $\OGE $
is *-isomorphic to the C*-algebra of the groupoid $\G _{G,E} $ described in \ref{ConcreteModel}, once the latter is
equipped with the topology generated by the basis of open sets $\Omega (\alpha ,g,\beta )$ described in \ref{ReSlice}, for all
$(\alpha ,g,\beta )$ in $\SGE $.

\section $\OGE $ as a Cuntz-Pimsner algebra

Inspired by Nekrashevych's paper \cite {NekraJO}, we will now give a description of $\OGE $ as a Cuntz-Pimsner algebra
\cite {Pimsner}.  With this we will also be able to prove that $\OGE $ is nuclear and that $\GpdGE $ is amenable when
$G$ is an amenable group.  As before, we will work under the conditions of \ref{StandingHyp}.

We begin by introducing the algebra of coefficients over which the relevant Hilbert bimodule, also known as a
correspondence, will later be constructed.

Since the action of $G$ on $E$ preserves length by \ref{extendedaction.\KeepLength }, we see that the set of
vertices of $E$ is $G$-invariant, so we get an action of $G$ on $E^0$ by restriction.  By dualization $G$ acts on the
algebra $C(E^0)$ of complex valued
  functions\fn {Notice that, since $E^0$ is a finite set, $C(E^0)$ is nothing but ${\bf C}^{|E^0|}$.}
  on $E^0$.
  We may therefore form the crossed-product C*-algebra
  $$
  A = C(E^0) \ifundef {rtimes} \times \else \rtimes \fi G.
  $$

  Since $C(E^0)$ is a unital algebra, there is a canonical unitary representation of $G$ in the crossed product, which we
will denote by $\{\gp _g\}_{g\in G}$.

  On the other hand, $C(E^0)$ is also canonically isomorphic to a subalgebra of $A$ and we will therefore identify these
two algebras without further warnings.

  For each $\vr $ in $E^0$, we will denote the characteristic function of the singleton $\{\vr \}$ by $\proj _\vr $, so
that $\{\proj _\vr : \vr \in E^0\}$ is the canonical basis of $C(E^0)$, and thus $A$ coincides with the closed linear span of
the set
  $$
  \big \{\proj _\vr \gp _g: \vr \in E^0,\ g\in G\big \}.
  \equationmark LinSpanForCP
  $$

For later reference, notice that the covariance condition in the crossed product reads
  $$
  \gp _g \proj _\vr = \proj _{g\vr } \gp _g
  \for \vr \in E^0 \for g\in G.
  \equationmark CovarCondCrossProdinho
  $$

Our next step is to construct a correspondence over $A$.  In preparation for this we denote by $A^\ed $ the right ideal
of $A$ generated by $\proj _{\src (\ed )} $, for each $\ed \in E^1$.  In technical terms
  $$
  A^\ed = \proj _{\src (\ed )} A.
  $$

  With the obvious right $A$-module structure, and the inner product defined by
  $$
  \langle y,z\rangle = y^*z \for y,z\in A^\ed ,
  $$
  one has that $A^\ed $ is a right Hilbert $A$-module.  Notice that this is not necessarily a full Hilbert module since
$\langle A^\ed ,A^\ed \rangle $ is the two-sided ideal of $A$ generated by $\proj _{\src (\ed )} $, which might be a proper ideal in
some cases.

As already seen in \ref{LinSpanForCP}, $A$ is spanned by the elements of the form $\proj _\vr \gp _g$.  Therefore
$A^\ed $ is spanned by the elements of the form
  $
  \proj _{\src (\ed )} \proj _\vr \gp _g,
  $
  but, since the $\proj $'s are mutually orthogonal, this is either zero or equal to $\proj _{\src (\ed )} \gp _g$.
Therefore we see that
  $$
  A^\ed = \overline {\hbox {span}} \{\proj _{\src (\ed )} \gp _g : g\in G\}.
  $$

Introducing the right Hilbert $A$-module which will later be given the structure of a correspondence over $A$, we define
  $$
  M = \bigoplus _{\ed \in E^1} A^e.
  $$

Observe that if $\vr $ is a vertex which is the source of many edges, say
  $$
  \src \inv (\vr ) = \{\ed _1,\ed _2,\ldots ,\ed _n\},
  $$
  then
  $$
  A^{\ed _i} = \proj _{\src (\ed _i)} A = \proj _\vr A,
  $$
  for all $i$, so that $\proj _\vr A$ appears many times as a direct summand of $M$.  However these copies of $\proj
_\vr A$ should be suitably distinguished, according to which edge $\ed _i$ is being considered.

On the other hand, notice that if $\src \inv (\vr )= \ifundef {varnothing} \emptyset \else \varnothing \fi $, then $\proj _\vr A$ does not appear among the summands of $M$,
at all.

Addressing the fullness of $M$, observe that
  $$
  \langle M,M\rangle = \sum _{{\vr \in E^0 \atop \src \inv (\vr )\neq \ifundef {varnothing} \emptyset \else \varnothing \fi }}A \proj _\vr A,
  $$
  so, when $E$ has no \"{sinks}, that is, when $\src \inv (\vr )$ is nonempty for every $\vr $, one has that $M$ is
full.

Given $\ed \in E^1$, the element $\proj _{\src (\ed )}$, when viewed as an element of $A^\ed \subseteq M$, will play a very special
role in what follows, so we will give it a special notation, namely
  $$
  \t _\ed := \proj _{\src (\ed )}.
  \equationmark DefineTe
  $$

There is a small risk of confusion here in the sense that, if $\ed _1 ,\ed _2 \in E^1$ are such that
  $$
  \vr := \src (\ed _1)=\src (\ed _2),
  $$
  then \ref{DefineTe} assigns $\proj _\vr $ to both $\t _{\ed _1}$ and $\t _{\ed _2}$.  However the coordinate in
which $\proj _\vr $ appears in $\t _{\ed _i}$ is determined by the corresponding $\ed _i$, so if $\ed _1\neq \ed _2$, then $\t
_{\ed _1}\neq \t _{\ed _2}$.

In order to completely dispel any confusion, here is the technical definition:
  $$
  \t _\ed = (m_\oed )_{\oed \in E^1},
  $$
  where
  $$
  m_\oed = \left \{\matrix { \proj _{\src (\ed )}, &\hbox {if } \oed =\ed ,\hfill \cr \pilar {12pt} 0, &\hbox
{otherwise.}}\right .
  $$

  We should notice that
  $$
  \t _\ed \proj _{\src (\ed )} = \t _\ed ,
  \equationmark TedProj
  $$
  and that any element $y\in M$ may be written uniquely as
  $$
  y=\sum _{\ed \in E^1} \t _\ed y_\ed ,
  \equationmark UniqueDecInM
  $$
  where each $y_\ed \in A^\ed $.

As the next step in constructing a correspondence over $A$,
  we would now like to define a certain *-homomorphism from $A$ to the algebra $\Lin (M)$ of adjointable linear
operators on $M$.
  Since $A$ is a crossed product algebra, this will be accomplished once we produce a covariant representation
  $(\repalg ,\rep )$ of the C*-dynamical system $\big(C(E^0),G\big)$ on $M$.  We begin with the group representation $\rep $.

\definition For each $g\in G$, let $\rep _g$ be the linear operator on $M$ given by
  $$
  \rep _g\Big (\sum _{\ed \in E^1} \t _\ed y_\ed \Big ) = \sum _{\ed \in E^1} \t _{g\ed } \gp _{\varphi (g,\ed )}y_\ed ,
  $$
  whenever $y_\ed \in A^\ed $, for each $\ed $ in $E^1$.

By the uniqueness in \ref{UniqueDecInM}, it is clear that $\rep _g$ is well defined.

\state Proposition Each $\rep _g$ is a unitary operator on $M$.  Moreover, the correspondence $g\mapsto \rep _g$ is a unitary
representation of $G$.

\Proof Pick $g$ in $G$.  We begin by claiming that the two sides in the expression defining $\rep _g$, above, coincide
whenever the $y_\ed $ are in $A$, even if $y_\ed $ does not belong to $A^\ed $.  Since $\rep _g$ is clearly additive, we
only need to check that
  $$
  \rep _g(\t _\ed y) = \t _{g\ed } \gp _{\varphi (g,\ed )}y
  \for y\in A.
  $$
  Observing that $\t _\ed = \t _\ed \proj _{\src (\ed )}$, we have
  $$
  \rep _g(\t _\ed y) =
  \rep _g(\t _\ed \proj _{\src (\ed )}y) =
  \t _{g\ed } \gp _{\varphi (g,\ed )}\proj _{\src (\ed )}y \$=
  \t _{g\ed } \proj _{\varphi (g,\ed )\src (\ed )}\gp _{\varphi (g,\ed )}y \={Equacoes.\PhiStarOnVert }
  \t _{g\ed } \proj _{\src (g\ed )}\gp _{\varphi (g,\ed )}y =
  \t _{g\ed } \gp _{\varphi (g,\ed )}y,
  $$
  proving the claim.  One therefore concludes that $\rep _g$ is right-$A$-linear.

We next claim that, for all $\ed ,\oed \in E^1$, one has
  $$
  \langle \rep _g(\t _\ed ),\t _\oed \rangle =
  \langle \t _\ed ,\rep _{g\inv }(\t _\oed )\rangle .
  \equationmark ClaimForAdjoint
  $$
  We have
  $$
  \langle \rep _g(\t _\ed ),\t _{\oed }\rangle =
  \langle \t _{g\ed } \gp _{\varphi (g,\ed )},\t _{\oed }\rangle =
   \gp _{\varphi (g,\ed )}^* \langle \t _{g\ed },\t _{\oed }\rangle =
  \equal {g\ed }{\oed } \gp _{\varphi (g,\ed )}\inv \proj _{\src (g\ed )} \$=
  \equal {g\ed }{\oed } \proj _{\varphi (g,\ed )\inv \src (g\ed )} \gp _{\varphi (g,\ed )}\inv \={Equacoes.\PhiStarOnVert }
  \equal {g\ed }{\oed } \proj _{\src (\ed )} \gp _{\varphi (g,\ed )}\inv = (\star ).
  $$
  Starting from the right-hand-side of \ref{ClaimForAdjoint}, we have
  $$
  \langle \t _\ed ,\rep _{g\inv }(\t _{\oed })\rangle =
  \langle \t _\ed ,\t _{g\inv \oed } \gp _{\varphi (g\inv ,\oed )}\rangle =
  \equal {e}{g\inv \oed } \proj _{\src (\ed )} \gp _{\varphi (g\inv ,\oed )} \$=
  \equal {ge}{\oed } \proj _{\src (\ed )} \gp _{\varphi (g,g\inv \oed )\inv } =
  \equal {ge}{\oed } \proj _{\src (\ed )} \gp _{\varphi (g,\ed )}\inv ,
  $$
  which agrees with $(\star )$ above, and hence proves claim \ref{ClaimForAdjoint}. If $y,z\in A$, we then have that
  $$
  \langle \rep _g(\t _\ed y ),\t _\oed z\rangle =
  y^*\langle \rep _g(\t _\ed ),\t _\oed \rangle z =
  y^*\langle \t _\ed ,\rep _{g\inv }(\t _\oed )\rangle z =
  \langle \t _\ed y,\rep _{g\inv }(\t _\oed z)\rangle ,
  $$
  from where one sees that $\langle \rep _g(\xi ),\eta \rangle = \langle \xi ,\rep _{g\inv }(\eta )\rangle $, for all $\xi ,\eta \in M$, hence proving that $\rep _g$ is an
adjointable operator with
  $
  \rep _g^* = \rep _{g\inv }.
  $

Let us next prove that
  $$
  \rep _g\rep _h = \rep _{gh}
  \for g,h\in G.
  $$
  By $A$-linearity it is enough to prove that these operators coincide on the set formed by the $\t _\ed $'s, which is a
generating set for $M$.  We thus compute
  $$
  \rep _g\big(\rep _h(\t _\ed )\big) =
  \rep _g\big(\t _{h\ed } \gp _{\varphi (h,\ed )}\big) =
  \t _{gh\ed } \gp _{\varphi (g,h\ed )} \gp _{\varphi (h,\ed )}=
  \t _{gh\ed } \gp _{\varphi (gh,\ed )} =
  \rep _{gh}(\t _\ed ).
  $$

Since it is evident that $\rep _1$ is the identity operator on $M$ we obtain, as a consequence, that
  $
  \rep _g\inv =
  \rep _{g\inv } =
  \rep _g^*,
  $
  so each $\rep _g$ is unitary and the proof is concluded.  \endProof

In order to complete our covariant pair we must now construct a *-homomorphism from $C(E^0)$ to $\Lin (M)$.  With this in
mind we give the following:

\definition For every $\vr $ in $E^0$, let
  $$
  M_\vr = \bigoplus _{\ed \in \ran \inv (\vr )} A^\ed ,
  $$
  which we view as a complemented sub-module of $M$.  In addition, we let $\projMod _\vr $ be the orthogonal projection
from $M$ to $M_\vr $, so that
  $$
  \projMod _\vr (\t _\ed y) = \equal {\ran (\ed )}{\vr }\t _\ed y
  \for \ed \in E^1 \for y\in A.
  \equationmark FormulaForQ
  $$

Observe that the $\projMod _\vr $ are pairwise orthogonal projections and that
  $
  \sum _{\vr \in E^0} \projMod _\vr = 1.
  $

\definition
  Let $\repalg :C(E^0) \to \Lin (M)$ be the unique unital *-homomorphism such that
  $$
  \repalg (\proj _\vr ) = \projMod _\vr
  \for \vr \in E^0.
  $$

From our working hypothesis that $E$ has no sources, we see that for every $\vr $ in $E^0$, there is some $\ed \in E^1$ such
that $\ran (\ed )=\vr $.  So
  $$
  \projMod _\vr (\t _\ed )=\t _\ed ,
  $$
  whence $\projMod _\vr \neq 0$.  Consequently $\repalg $ is injective.

\state Proposition The pair $(\repalg ,\rep )$ is a covariant representation of the C*-dynamical system $\big(C(E^0),G\big)$ in
$\Lin (M)$.

\Proof All we must do is check the covariance condition
  $$
  \rep _g \repalg (y) = \repalg \big(\auto _g(y)\big) \rep _g
  \for g\in G \for y\in C(E^0),
  $$
  where $\auto $ is the name we temporarily give to the action of $G$ on $C(E^0)$.  Since $C(E^0)$ is spanned by the
$\proj _\vr $, it suffices to consider $y=\proj _\vr $, in which case the above identity becomes
  $$
  \rep _g \projMod _\vr =
  \projMod _{g\vr } \rep _g.
  \equationmark ShortCovar
  $$

  Furthermore $M$ is generated, as an $A$-module, by the $\t _\ed $, for $\ed \in E^1$, so we only need to verify this on
the $\t _\ed $.  We have
  $$
  \rep _g \big(\projMod _\vr (\t _\ed )\big) =
  \equal {\ran (\ed )}{\vr } \rep _g (\t _\ed ) =
  \equal {\ran (\ed )}{\vr } \t _{g\ed } \gp _{\varphi (g,e)},
  $$
  while
  $$
  \projMod _{g\vr } \big(\rep _g (\t _\ed )\big) =
  \projMod _{g\vr } \big(\t _{g\ed } \gp _{\varphi (g,e)}\big) =
  \equal {\ran (g\ed )}{g\vr } \t _{g\ed } \gp _{\varphi (g,e)},
  $$
  verifying \ref{ShortCovar} and concluding the proof.
  \endProof

It follows from \cite [Proposition 7.6.4 and Theorem 7.6.6]{Pedersen} that there exists a *-ho\-mo\-mor\-phism
  $$
  \modmap : C(E^0) \ifundef {rtimes} \times \else \rtimes \fi G \to \Lin (M),
  $$
  such that
  $$
  \modmap (\proj _\vr )=\projMod _\vr \for \vr \in E^0,
  $$ and $$
  \modmap (\gp _g)=\rep _g \for g\in G.
  $$

Equipped with the left-$A$-module structure provided by $\modmap $, we then have that $M$ is a correspondence over $A$.

For later reference we record here a few useful calculations involving the left-module structure of $M$.

\state Proposition \label LeftModCalc
  Let $g\in G$, $\ed \in E^1$, and $\vr \in E^0$.  Then
  \Zitem $\gp _g \t _{\ed } = \t _{g\ed } \gp _{\varphi (g,\ed )}$,
  \zitem $\proj _\vr \gp _g \t _{\ed } = \equal {\ran (g\ed )}{\vr } \t _{g\ed } \gp _{\varphi (g,\ed )}.$

\Proof We have
  $$
  \gp _g \t _{\ed } =
  \modmap (\gp _g) \t _{\ed } =
  \rep _g (\t _{\ed } ) =
  \t _{g\ed } \gp _{\varphi (g,\ed )},
  $$
  proving (a).  Also
  $$
  \proj _\vr \gp _g \t _{\ed } =
  \modmap (\proj _\vr ) (\gp _g \t _{\ed }) =
  \projMod _\vr (\t _{g\ed } \gp _{\varphi (g,\ed )})=
  \equal {\ran (g\ed )}{\vr } \t _{g\ed } \gp _{\varphi (g,\ed )}.
  \endProof

It is our next goal to prove that $\OGE $ is naturally isomorphic to the Cuntz-Pimsner algebra associated to the
correspondence $M$, which we denote by $\O _M$.  As a first step, we identify a certain Cuntz-Krieger $E$-family.

\state Proposition \label RelInOM
  The following relations hold within $\O _M$.
  \iaitem
  \aitem For every $\vr \in E^0$, one has that $\sum _{\ed \in \ran \inv (\vr )}\t _\ed \t _\ed ^* = \proj _\vr $.
  \aitem $\sum _{\ed \in E^1}\t _\ed \t _\ed ^* = 1$.
  \aitem The set
  $
  \{\proj _\vr : \vr \in E^0\}\cup \{\t _\ed : \ed \in E^1\}
  $
  is a Cuntz-Krieger $E$-family.

\Proof We first claim that, for every $\vr \in E^0$, and every $m\in M$, one has that
  $$
  \sum _{\ed \in \ran \inv (\vr )}\t _\ed \t _\ed ^* m = \proj _\vr m.
  $$
  To prove it, it is enough to consider the case in which $m = \t _\oed $, for $\oed \in E^1$, since these generate $M$.  In
this case we have
  $$
  \sum _{\ed \in \ran \inv (\vr )}\t _\ed \t _\ed ^* \t _\oed =
  \equal {\ran (\oed )}{\vr } \t _\oed \t _\oed ^* \t _\oed =
  \equal {\ran (\oed )}{\vr } \t _\oed \={FormulaForQ}
  \projMod _\vr (\t _\oed ) =
  \proj _\vr \t _\oed ,
  $$
  proving the claim.  This says that the pair
  $
  \big(\proj _\vr , \sum _{\ed \in \ran \inv (\vr )}\t _\ed \t _\ed ^*\big)
  $
  is a redundancy or, adopting the terminology of \cite {Pimsner}, that the generalized compact operator
  $$
  \sum _{\ed \in \ran \inv (\vr )}\Omega _{\t _\ed , \t _\ed }
  $$
  is mapped to $\modmap (\proj _\vr )$ via $\Psi ^{(1)}$.
  Therefore
  $$
  \proj _\vr = \sum _{\ed \in \ran \inv (\vr )}\t _\ed \t _\ed ^*,
  $$
  in $\O _M$, proving (a).  Point (b) then follows from the fact that
  $
  \sum _{\vr \in E^0} \proj _\vr =1.
  $

Focusing now on (c), it is evident that $\{\proj _\vr : \vr \in E^0\}$ is a family of mutually orthogonal projections.
Moreover, for each $e\in E^1$, we have
  $$
  \t _\ed ^* \t _\ed =
  \langle \t _\ed , \t _\ed \rangle =
  \proj _{\src (\ed )},
  $$
  proving \ref{DefineCKFamily.\CkTwo } and also that $\t _\ed $ is a partial isometry.  Property \ref{DefineCKFamily.\CkOne } also holds in view of (a), so the proof is concluded.  \endProof

\state Proposition \label OneMapForCuntzPimsner
  There exists a unique surjective *-homomorphism
  $$
  \Lambda :\OGE \to \O _M
  $$
  such that
  $\Lambda (p_\vr ) = \proj _\vr $,
  $\Lambda (\s _\ed ) = \t _\ed $, and
  $\Lambda (u_g) = \gp _g$.

\Proof By the universal property of $\OGE $, in order to prove the existence of $\Lambda $ it is enough to check that the
  $\proj _\vr $,
  $\t _\ed $, and
  $\gp _g$
  satisfy the conditions of Definition \ref{DefineOGE}.

  Condition \ref{DefineOGE.a} has already been proved above while
  \ref{DefineOGE.b} is evidently true since $\gp $ is a representation of $G$ in $C(E^0) \ifundef {rtimes} \times  \else
\rtimes \fi G \subseteq \O _M$.
  Condition
  \ref{DefineOGE.c} is precisely \ref{LeftModCalc.i}, while
  \ref{DefineOGE.d} was taken care of in \ref{CovarCondCrossProdinho}.

Since $A$ is spanned by the $\proj _\vr $ and the $\gp _g$ by \ref{LinSpanForCP}, and since $M$ is generated over
$A$ by the $\t _\ed $, we see that $\O _M$ is spanned by the set
  $$
  \{\proj _\vr , \t _\ed , \gp _g : \vr \in E^0, \ \ed \in E^1, \ g\in G\},
  $$
  so $\Lambda $ is surjective.
  \endProof

Let us now prove that $\Lambda $ is invertible by providing an inverse to it.  Since $A$ is the crossed product C*-algebra
$C(E^0) \ifundef {rtimes} \times \else \rtimes \fi G$, one sees that \ref{DefineOGE.a\&d} guarantees the existence of a *-homomorphism
  $$
  \theta _A: A \to \OGE ,
  $$
  sending the $\proj _\vr $ to the $p_\vr $, and the $\gp _g$ to the $u_g$.  For each $\ed $ in $E^1$, consider the
linear mapping
  $$
  \theta _M : M \to \OGE ,
  $$
  given, for every $m = (m_\ed )_{\ed \in E^1} \in M$, by
  $$
  \theta _M(m)= \sum _{e\in E^1}\s _\ed \theta _A(m_\ed ) \in \OGE .
  $$
  Notice that $\theta _M (\t _\ed )= \s _\ed $, for all $\ed \in E^1$, because
  $$
  \theta _M (\t _\ed )=
  \s _\ed \theta _A(\proj _{\src (\ed )}) =
  \s _\ed p_{\src (\ed )} =
  \s _\ed .
  $$

\state Lemma The pair $(\theta _A,\theta _M)$ is a representation of the correspondence $M$ in the sense of \cite [Theorem
3.4]{Pimsner}, meaning that for all $y\in A$ and all $\xi ,\xi '\in M$,
  \izitem
  \zitem $\theta _M(\xi )\theta _A(y) = \theta _M(\xi y),$
  \zitem $\theta _A(y)\theta _M(\xi ) = \theta _M(y\xi ),$
  \zitem $\theta _M(\xi )^*\theta _M(\xi ') = \theta _A(\langle \xi ,\xi '\rangle ).$

\Proof
  Considering the various spanning sets at our disposal, we may assume that
  $y =
  \proj _\vr \gp _g$,
  that $\xi = \t _{\ed } z$, and $\xi ' = \t _{\ed '} z'$,
  with $\vr \in E^0$, $g\in G$, $\ed ,\ed '\in E^1$, $z\in \proj _{\src (\ed )} A$, and $z'\in \proj _{\src (\ed ')} A$.  We then have
  $$
  \theta _M(\xi )\theta _A(y) =
  \theta _M(\t _{\ed } z)\theta _A(y) =
  \s _{\ed } \theta _A(z)\theta _A(y) =
  \s _{\ed } \theta _A(zy) =
  \theta _M(\t _{\ed } zy) =
  \theta _M(\xi y),
  $$
  proving (i).  As for (ii), we have
  $$
  \theta _A(y)\theta _M(\xi ) =
  \theta _A(\proj _\vr \gp _g)\theta _M(\t _{\ed } z) =
  p_\vr u_g \s _{\ed } \theta _A(z) =
  p_\vr \s _{g\ed } u_{\varphi (g,\ed )} \theta _A(z) \$=
  \equal {\ran (g\ed )}{\vr } \s _{g\ed } \theta _A(\gp _{\varphi (g,\ed )} z) =
  \equal {\ran (g\ed )}{\vr } \theta _M(\t _{g\ed } \gp _{\varphi (g,\ed )}z\big) \={LeftModCalc.ii}
  \theta _M(\proj _\vr \gp _g \t _{\ed }z) =
  \theta _M(y\xi ),
  $$
  proving (ii).  Focusing now on (iii), we have
  $$
  \theta _M(\xi )^*\theta _M(\xi ') = (\s _{\ed } \theta _A(z)\big)^*\s _{\ed '} \theta _A(z') =
  \equal \ed {\ed '} \theta _A(z)^*p_{\src (e)} \theta _A(z') \$=
  \equal \ed {\ed '} \theta _A(z^*\proj _{\src (e)} z') =
  \theta _A(\langle \xi ,\xi '\rangle ).
  \endProof

It is well known \cite [Theorem 3.4]{Pimsner} that the Toeplitz algebra for the correspondence $M$, usually denoted
${\cal T}_M$, is universal for representations of $M$, so there exists a *-ho\-mo\-mor\-phism
  $$
  \Theta _0:{\cal T}_M\to \OGE ,
  $$
  coinciding with $\theta _A$ on $A$ and with $\theta _M$ on $M$.

\state Theorem \label CPPicture
  The map $\Theta _0$, defined above, factors through $\O _M$, providing a *-isomor\-phism
  $$
  \Theta :\O _M\to \OGE ,
  $$
  such that $\Theta (\proj _\vr ) = p_\vr $, $\Theta (\t _\ed )=\s _\ed $, and $\Theta (\gp _g) = u_g$, for all $\vr \in E^0$, $\ed \in E^1$, and
$g\in G$.

\Proof
  The factorization property
  follows immediately from \ref{RelInOM.b} and an easy modification of \cite [Proposition 7.1]{ExelVesshik} to
Cuntz-Pimsner algebras.

In order to prove that $\Theta $ is an isomorphism, observe that $\Theta \circ \Lambda $ coincides with the identity map on the generators of
$\OGE $, by \ref{OneMapForCuntzPimsner}, and hence $\Theta \circ \Lambda =id$.  The result then follows from the fact that $\Lambda $ is
surjective.
  \endProof

\state Corollary \label Amenabuilidade
  If $G$ is amenable then $\OGE $ is nuclear.

\Proof
  The amenability of $G$ ensures that $C(E^0) \ifundef {rtimes} \times \else \rtimes \fi G$ is nuclear.  The result then follows from Theorem \ref{CPPicture},
the fact that Toeplitz-Pimsner algebras over nuclear coefficient algebras is nuclear \cite [Theorem 4.6.25]{BO}, and
so are quotients of nuclear algebras \cite [Theorem 9.4.4]{BO}.  \endProof

\state Remark \rm Since $E^0$ is finite, the nuclearity of $C(E^0) \ifundef {rtimes} \times \else \rtimes \fi G$ is equivalent to the amenability of $G$.  However,
if the present construction is generalized for infinite graphs, one could produce examples of non amenable groups acting
amenably on $E^0$, in which case $C(E^0) \ifundef {rtimes} \times  \else \rtimes \fi G$ would be nuclear.  The proof of
Corollary \ref{Amenabuilidade} could then be adapted to prove that $\OGE $ is nuclear.

\state Corollary \label Newamenabuilidade
If $G$ is amenable, then $\GpdGE $ is an amenable groupoid. If moreover $(\Data )$ is pseudo free, then its sibling $\G _{G,E} $ is an amenable
groupoid.

\Proof For $\GpdGE $, it follows from \ref{Amenabuilidade}, \ref{UniversalTightAlgebra} and \cite [Theorem 5.6.18]{BO}. For $\G _{G,E} $, it follows from \ref{Amenabuilidade}, \ref{GroupoidModel} and \cite [Theorem 5.6.18]{BO}.
 \endProof

Nekrashevych has proven in \cite [Theorem 5.6]{NC}, that a certain groupoid of germs, denoted ${\cal D}_G$,
constructed in the context of self-similar groups, is amenable under the hypothesis that the group is \"{contracting}
and \"{self-replicating}.  Even though there are numerous differences between ${\cal D}_G$ and $\G _{G,E} $, including a
different notion of \"{germs} and Nekrashevych's requirement that group actions be \"{faithful}, we believe it should be
interesting to try to generalize Nekrashevych's result to our context.

\section Representing $C^*(E)$ and $G$ into $\OGE $

In this section we will study natural representations of the graph C*-algebra $C^*(E)$ and of the group $G$ in $\OGE $.
As before, we keep \ref{StandingHyp} in force.

Given that
  $$
  \{p_\vr : \vr \in E^0\}\cup \{\s _\ed : \ed \in E^1\}
  $$
  is a Cuntz-Krieger $E$-family, the universal property of the graph C*-algebra $C^*(E)$ \cite {Raeburn} provides for
the existence of a *-homomorphism
  $$
  \iota : C^*(E)\to \OGE ,
  $$
  sending the canonical Cuntz-Krieger $E$-family of $C^*(E)$ to the corresponding one within $\OGE $.

\state Proposition \label PropGraphSeaStarisSubalg
  The *-homomorphism $\iota $ above is injective.

  \Proof
  Using the universal property of $\OGE $, it is easy to see that, for each complex number $z$, with $|z|=1$, there is a
*-homomorphism
  $$
  \gamma _z:\OGE \to \OGE ,
  $$
  satisfying
  $$
  \gamma _z(p_x)=p_x,\quad
  \gamma _z(s_\ed )=zs_\ed \and
  \gamma _z(u_g)=u_g,
  $$
  for all $x\in E^0$, $\ed \in E^1$ and $g\in G$.  It is also easy to see that the correspondence $z \to  \gamma _z$ defines an action of
the circle group on $\OGE $, and moreover that $\iota $ is covariant relative to this action on $\OGE $, on the one hand, and
the standard gauge action on $C^*(E)$, on the other.  In order to prove the injectivity of $\iota $ we may then apply the
gauge invariant uniqueness Theorem \cite [Theorem 2.2]{Raeburn}, which requires, in addition, that we verify that the
$p_x$ are nonzero.

To prove this we observe that, in the groupoid model of $\OGE $ given by \ref{GroupoidModel}, for each $x$ in $E^0$,
the element $p_x$ is the characteristic function of the cylinder $\cyl {x}$, seen as a subset of $E^\infty $, which in turn is
the unit space of the groupoid $\G _{G,E} $.  Since $E$ has no sources, we have that $\cyl {x}$ is nonempty, whence $p_x$ is
nonzero, as required.  This concludes the proof.
  \endProof

  In particular,
  \redundantref {Proposition}{PropGraphSeaStarisSubalg} implies that $\SE $ is $\ast $-isomorphic to the
inverse semigroup of $\OGE $ generated by $\{s_a: a\in E^1\}$.

With respect to the injectivity of the representation of $G$ into $\OGE $, we have to work a bit more to obtain a result
in the line of
  \redundantref {Proposition}{PropGraphSeaStarisSubalg}.

\state Lemma \label LempiInjImplyUInjUinj
  Let $\pi : \SGE \rightarrow \OGE $ and $u:G\rightarrow \OGE $ be the natural maps. If $\pi $ is injective, the so is $u$.

  \Proof
  Let $g\in G$ such that $u_g=u_1$. For any $x\in E^0$ we have $\pi (x,g,g\inv x)=p_xu_g$ and $\pi  (x,1,x)=p_xu_1=p_x$. Since
$u_g=u_1$, we get $(x,g,g\inv x)=(x,1,x)\in \SGE $, whence $g=1$.
  \endProof

We need to recall some extra definitions. Let $\mathcal {G}$ be an \'etale groupoid, i.e. a topological groupoid whose
unit space $\mathcal {G}^{(0)}$ is locally compact and Hausdorff in the relative topology, and such that the range map
$r:\mathcal {G}\rightarrow \mathcal {G}^{(0)}$ is a local homeomorphism (and then so is the source map $d:\mathcal
{G}\rightarrow \mathcal {G}^{(0)}$). An open set $U\subset \mathcal {G}$ is a slice if the restrictions of $r$ and $d$ to $U$
are injective (see e.g. \cite {pat}). In particular, $\mathcal {G}^{(0)}$ is a slice \cite [Proposition 3.4]{actions}, and the collection of all slices forms a basis for the topology of $\mathcal {G}$ \cite [Proposition 3.5]{actions}.

\definition
  We denote by ${\mathcal {S}\ell }(G,E)$ the set of all compact slices. It is well known (see e.g. \cite [Proposition
2.2.4]{pat}) that ${\mathcal {S}\ell }(G,E)$ forms a $\ast $-inverse semigroup with the operations
  $$UV=\{uv: u\in U, v\in V, (u,v)\in \mathcal {G}^{(2)}\}, \text { and } U^*=\{u\inv : u\in U\}.$$ Moreover, if $\mathcal {U}_{G,E}=\{
1_U: U \text { is a compact slice}\}\subseteq C^*(\GpdGE )$ is the semigroup formed by their characteristic functions, then
$$\mathcal {U}_{G,E}\cong {\mathcal {S}\ell }(G,E). \equationmark IsoInvSemgrp $$

Fix the canonical action $\theta $ of $\SGE $ on $E^{\infty }$. Given any $(\alpha , g, \beta )\in \SGE $, notice that the domain $\text
{Dom}(\theta _{(\alpha , g, \beta )})$ of the partial homeomorphism of $E^{\infty }$ given by the action of $(\alpha , g, \beta )$ is $\cyl {\beta }$. Now,
given $(\alpha , g, \beta )\in \SGE $ and any open set $U\subseteq \cyl {\beta }$, set (see Section 9) $$\Theta ((\alpha , g, \beta ), U)=\{[\alpha , g, \beta ; \eta ]: \eta \in U\}.$$
According to \cite [Proposition 4.18]{actions}, for every $(\alpha , g, \beta )\in \SGE $ and every open set $U\subseteq \cyl {\beta }$, $\Theta ((\alpha , g,
\beta ), U)$ is a slice (in fact, they form a basis for the topology of $\GpdGE $). Then, we have the following result

\state Lemma \label LemVariousIsos
  ${\mathcal {S}\ell }(G,E)=\langle \Theta ((\alpha , g, \beta ), \cyl {\beta }): (\alpha , g, \beta )\in \SGE \rangle $.

    \Proof
  By \cite [Proposition 5.13(7)]{Steinberg}, $C^*(\GpdGE )$ is generated by $$\{1_{\Theta ((\alpha , g, \beta ), \cyl {\beta })}: (\alpha , g,
\beta )\in \SGE \}.$$ Thus, the result holds by \ref{IsoInvSemgrp}.
  \endProof

The next result is the key point for proving the injectivity of $u:G\rightarrow \OGE $.

\state Lemma \label LemFirstIsoRF
  If $(G,E, \varphi )$ is pseudo free, then the map
  $$
  \matrix { \rho :& \SGE & \rightarrow & {\mathcal {S}\ell }(G,E) \cr
 & (\alpha , g, \beta ) & \mapsto & \Theta {((\alpha , g, \beta ), \cyl {\beta })} }
   $$
   is a $\ast $-semigroup isomorphism.

  \Proof
  The surjectivity of $\rho $ derives from
  \redundantref {Lemma}{LemVariousIsos}.

Now, let $(\alpha , g, \beta ),(\gamma , h, \eta )\in \SGE $ such that $\Theta ((\alpha , g, \beta ), \cyl {\beta })=\Theta  ((\gamma , h, \eta ),\cyl {\eta })$. Then, for any $\omega \in \cyl
{\beta }$ we have $[\alpha , g, \beta ; \omega ]=[\gamma , h, \eta ; \omega ]$. Since $(G,E, \varphi )$ is pseudo free, by
  \redundantref {Proposition}{EqualGerms}
  there exists $\tau \in E^*$ such that $\gamma =\alpha \cdot g\tau $, $\eta =\beta \tau $ and $h=\varphi (g, \tau  )$. If $(\alpha , g, \beta )\ne (\gamma , h, \eta )$, then we can pick $\delta \ne
\tau $ in $E^*$ and $\widehat {\omega  }=\beta \delta \widetilde {\omega }\in \cyl {\beta }$. Thus, $\widehat {\omega }\in \cyl {\eta }$ but $[(\gamma , h, \eta ), \widehat {\omega }]$
is not defined, contradicting the hypothesis. Hence, $(\alpha , g, \beta )=(\gamma , h, \eta )$, whence $\rho $ is injective, as desired.
\endProof

\state Proposition \label PropKeyIso
  There exists a $\ast $-isomorphism $\phi : \OGE \rightarrow C^*(\GpdGE )$ such that $\phi (p_x)=1_{\Theta ((x,1,x), \cyl {x})}$ for
every $x\in E^0$, $\phi (s_a)=1_{\Theta ((a, 1, d(a)), \cyl {a})}$ for every $a\in E^1$, and $\phi (u_g)= \sum  \limits _{x\in E^0}1_{\Theta ((x, g, g\inv
x), Z(g\inv x))}$ for every $g\in G$.

  \Proof
  Notice that $u_g=\sum \limits _{x\in E^0}u_gp_x$. Then it is direct but tedious to check that $$\left \{1_{\Theta  ((x,1,x), \cyl
{x})}: x\in E^0\right \}\cup \left \{1_{\Theta ((a,1,d(a)), \cyl {a})}: a\in E^1 \right \}\cup $$ $$\cup \left \{\sum \limits _{x\in E^0}1_{\Theta ((x, g, g\inv x),
Z(g\inv x))} : g\in G\right \}$$ satisfy the defining relations for $\OGE $. Thus, by the Universal Property of $\OGE $, the
map $\varphi $ is an $\ast $-homomorphism. Notice that $\varphi $ is the homomorphism given in \cite [Theorem 13.3]{actions}, and so
is injective. Surjectivity is due to \cite [Proposition 5.13(7)]{Steinberg}.
  \endProof

\state Corollary \label CorpiIsInj
  If $(G,E, \varphi )$ is pseudo free, then $\pi : \SGE \rightarrow \OGE $ is injective.

  \Proof
  The composition map $$ \matrix {
 \SGE & \rightarrow & \OGE & \rightarrow & C^*(\GpdGE )\cr
 (\alpha , g , \beta ) & \mapsto & s_{\alpha }u_gs_{\beta }^* & \mapsto & 1_{\Theta {((\alpha , g, \beta ), \cyl {\beta })}} } $$ is injective by
  \redundantref {Lemma}{LemVariousIsos} and
  \redundantref {Lemma}{LemFirstIsoRF}. By
  \redundantref {Proposition}{PropKeyIso},
  $$
  \OGE \rightarrow C^*(\GpdGE )
  $$
  is injective. Thus, $\pi $ is injective. \endProof

Hence, we conclude

\state Proposition \label PropGisSubgroup
  If $(G,E, \varphi )$ is pseudo free, then $u:G\rightarrow \OGE $ is injective.

  \Proof
  By Corollary \ref{CorpiIsInj}, $\pi $ is injective. Then, so is $u$ by
  \redundantref {Lemma}{LempiInjImplyUInjUinj}.
  \endProof

\state Remark
  In particular,
  \redundantref {Proposition}{PropGisSubgroup} implies that, if $(G,E, \varphi )$ is pseudo free, then $\SGE $ is $\ast
$-isomorphic to the inverse semigroup of $\OGE $ generated by $\{s_a: a\in E^1\}\cup \{p_xu_g: x\in E^0, g\in  G\}$.

Proposition \ref{PropGisSubgroup} provides the best situation possible, as the next example shows:

\state Example \rm
  Let $E$ be the graph with only one vertex and one edge, and let $G$ be any noncommutative group. Fix the trivial
action of $G$ on $E$, and let $\varphi $ be the one-cocycle of $G$ defined by $\varphi (g, a)=1$ for every $g\in G, a\in  E^1$. Then, it is
easy to see that $\OGE \cong C^*(E)\cong C(\mathbb {T})$, which is a commutative C*-algebra, so that it cannot contain
any faithful copy of $G$.

\section The Hausdorff property for $\GpdGE $

  Again considering a triple $(\Data )$ satisfying \ref{StandingHyp},
we will now give a characterization of the Hausdorff property for the tight groupoid of $\SGE$.  The first result we may present  in this direction  is:

\state Proposition
  If $(\Data )$ is pseudo free, then $\GpdGE $ is a Hausdorff groupoid.

\Proof
  If $(\Data )$ is pseudo free, then $\SGE $ is E*-unitary by \ref{EssFreeEUnitary}, so $\GpdGE $ is Hausdorff by   \cite[Corollary 3.17]{EPFour}.
This could also be obtained from \cite [Propositions 6.4 and 6.2]{actions}.
  \endProof

The converse of the above result is not true: as we will see in Example \ref {KatsuNoPSbutHaus}, there are examples in
which $(\Data )$ fails to be pseudo free but still $\GpdGE $ is Hausdorff.

  This may be interpreted as saying that the above assumption that $(\Data)$ is pseudo free is a much too strong
hypothesis which one would therefore like to relax.

  On the other hand, recall from \ref{EssFreePath} that the failure of pseudo freeness for $(\Data)$ is equivalent to
the existence of strongly fixed paths for nontrivial group elements.  The result below consists in allowing a limited
amount of minimal strongly fixed paths, and hence a limited number of counter-examples for pseudo freeness, without
harming Hausdorffnes.

\state Theorem \label MainHausdorff
  Assuming that $(\Data )$ satisfies \ref{StandingHyp}, the following are equivalent:
  \iaitem
  \aitem for every $g$ in $G$, there are at most finitely many minimal strongly fixed paths for $g$,
  \aitem $\GpdGE $ is Hausdorff.

\Proof
  We will of course use \cite[Theorem 3.16]{EPFour}.  So, given $s$ in $\SGE $, we must provide a finite cover for $\J_s$.
Since such a cover exists by trivial reasons when $\J _s$ is empty, let us assume that $s$ dominates at least one
nonzero idempotent element.
  By \ref{LemDominateIdempotent} we then have that $s$ necessarily has the form
  $$
  s=(\alpha ,g,\alpha ),
  $$
  and the set of nonzero idempotent elements dominated by $s$ is given by
  $$
  \J _s=\big \{(\alpha \tau ,1,\alpha \tau ): \tau \in E^*,\ \src (\alpha ) = \ran (\tau ),\ \tau \text { is strongly fixed by } g\big \}.
  $$

  Using \ref{StrFixElts} we may further describe $\J _s$ as
  $$
  \J _s =
  \big \{(\alpha \mu \gamma ,1,\alpha \mu \gamma ): \mu \in M_g,\ \gamma \in E^*,\ \src (\alpha ) = \ran (\mu ),\ \src (\mu ) = \ran (\gamma )\big \},
  \equationmark DecomposeJs
  $$
  where $M_g$ is the set of all minimal strongly fixed paths for $g$.

Assuming (a), we have that $M_g$ is finite and then it is clear that
  $$
  \big \{(\alpha \mu ,1,\alpha \mu ): \mu \in M_g,\ \src (\alpha ) = \ran (\mu ) \big \}
  $$
  is a finite cover for $\J _s$, whence $\GpdGE $ is Hausdorff by \cite[Theorem 3.16]{EPFour}.

Assuming (b), let $g\in G$, and for each vertex $\vr $ in $E^0$, denote by $M^\vr _g$ the set of all minimal strongly fixed
paths for $g$ whose range coincides with $\vr $.  Since $E$ is finite, in order to prove that $M_g$ is finite, it is
enough to check that each $M_g^\vr $ is finite.

If $M_g^\vr $ is empty, there is nothing to do, so let us assume the contrary.  Given any $\mu $ in $M_g^\vr $, we then
have that
  $$
  \vr = \ran (\mu ) = \ran (g\mu ) = g\ran (\mu ) = g\vr ,
  $$
  so $\vr $ is fixed by $g$.  Consequently $s:=(\vr ,g,\vr )$ lies in $\SGE $ and, assuming (b), we have by \cite[Theorem 3.16]{EPFour} that $\J _s$ admits a finite cover which, in view of \ref{DecomposeJs}, must necessarily be of
the form
  $$
  \big \{(\mu _i\gamma _i,1,\mu _i\gamma _i)\}_{i=1}^n,
  $$
  where the $\mu _i\in M_g^\vr $, and the $\gamma _i$ are paths with $\src (\mu _i)=\ran (\gamma _i)$.  We then claim that the $\mu _i$ apearing
above exhaust $M_g^\vr $, meaning that
  $$
  M_g^\vr = \{\mu _i: i=1,\ldots ,n\}.
  \equationmark MgEqualsCover
  $$

  To see this, let $\mu \in M_g^\vr $, so that $(\mu ,1,\mu )\leq s$, by \ref{LemDominateIdempotent}, and hence $(\mu ,1,\mu )\in \J _s$.
For some $i$, one would then have that
  $$
  (\mu _i\gamma _i,1,\mu _i\gamma _i)(\mu ,1,\mu )\neq 0,
  $$
  in which case either $\mu $ is a prefix of $\mu _i\gamma _i$, or vice versa.  This implies that either $\mu $ is a prefix of $\mu _i$,
or vice versa, but since both $\mu $ and $\mu _i$ are minimal, we must have $\mu =\mu _i$, proving \ref{MgEqualsCover}, and
hence that $M_g^\vr $ is finite.  Consequently $M_g$, which decomposes as the disjoint union of the $M_g^\vr $, is also
finite.  This verifies (a) and hence concludes the proof.
  \endProof

\section Minimality for $\GpdGE $

In this section we will study conditions under which $\GpdGE $ is minimal. Some of the results we obtain here are analog to
those proved in \cite {ExelLaca} for the case of partial actions of groups.

Given a triple $(\Data )$ satisfying \ref{StandingHyp},
there are two relations among vertices in $E^0$ which are relevant for the question at hand.  First of all let us say
that
  $$
  x \rightharpoonup y
  $$
  provided there exists a path $\alpha $ in $E^*$ such that $\src (\alpha )=x$ and $\ran (\alpha )=y$.  Notice that this relation is
reflexive (take $\alpha $ to be $x$) and transitive (take the concatenation of the relevant paths).  However this is neither
symmetric nor antisymmetric, hence it is not an equivalence relation nor an order relation.

 The other relation we have in mind is simply the orbit relation, defined by
  $$
  x\sim y
  $$
  when there exists $g$ in $G$ such that $gx=y$.  Unlike ``$\rightharpoonup $'', it is well known that ``$\sim $'' is an
equivalence relation.

We may then consider the smallest transitive relation extending both ``$\rightharpoonup $'' and ``$\sim $'', by saying that
vertices $x$ and $y$ are related when one may find a sequence of vertices $x_0,x_1,\ldots ,x_{2n}$ such that
  $$
  x =
  x_0 \rightharpoonup x_1 \sim x_2 \rightharpoonup x_3 \sim \ldots  \sim x_{2n-2} \rightharpoonup x_{2n-1} \sim x_{2n} = y.
  \equationmark MultiRelation
  $$

  The situation is in fact not so complicated due to the following:

\state Proposition \label TwoRelations
  Let $x$ and $y$ be vertices in $E^0$.  Then the following are equivalent;
  \izitem
  \zitem there exists a vertex $u$ such that $x \rightharpoonup u \sim y$,
  \zitem there exists a vertex $v$ such that $x \sim v \rightharpoonup y$.

\Proof
  The fact that $x \rightharpoonup u \sim y$ means that there exists a path $\alpha $ in $E^*$ such that $\src (\alpha ) =x$, and
$\ran (\alpha )=u$, and there exists some $g$ in $G$ such that $gu=y$.  Considering the path $\beta =g\alpha $, and the vertex $v=gx$,
notice that
  $$
  \src (\beta ) = \src (g\alpha ) = g\src (\alpha ) = gx = v,
  $$
  while
  $$
  \ran (\beta ) = \ran (g\alpha ) = g\ran (\alpha ) = gu = y,
  $$
  so $x\sim v \rightharpoonup y$.  Conversely, assuming (ii) we have that $gx=v=\src (\beta )$ and $\ran (\beta )=y$, for suitable
$g$ in $G$ and $\beta $ in $E^*$.  Defining $u=g\inv y$, and $\alpha =g\inv \beta $, we then have that
  $$
  \src (\alpha ) = g\inv \src (\beta ) = x,
  $$
  and
  $$
  \ran (\alpha ) = g\inv \ran (\beta ) = g\inv y = u,
  $$
  so $x \rightharpoonup u \sim y$.
  \endProof

\definition Given $x$ and $y$ in $E^0$, we will say that
  $$
  x\gg y
  $$
  if the equivalent conditions of \ref{TwoRelations} are satisfied.

Observe that ``$\gg $'' coincides with the relation defined in \ref{MultiRelation}, thanks to \ref{TwoRelations}, and hence it is clearly transitive.  It is also evident that ``$\gg $'' is reflexive but, again, it is neither
symmetric nor antisymmetric.  Nevertheless we will view it as a \"{defective order relation}, in the sense that it
satisfies all of the postulates of a (partial) order relation but for antisymmetry.

Anytime we have such a defective order relation, it is possible to turn it into a bona fide partial order by
identifying elements whenever antisymmetry fails.  By this we mean that two vertices $x$ and $y$ in $E^0$ will be called
equivalent, in symbols
  $$
  x\approx y
  $$
  whenever $x\gg y$ and $y\gg x$.  Writing $[x]$ for the equivalence class of each $x$ in $E^0$, the set of all
equivalent classes, namely
  $$
  {E^0\over \approx } = \big \{[x]: x\in E^0\big \}
  $$
  becomes a partially ordered set via the well defined order relation
  $$
  [x] \geq  [y] \iff x \gg y.
  $$

\definition
  Under the assumptions of \ref{StandingHyp}, we will say that:
  \izitem
  \zitem $E$ is \"{$G$-transitive} if, for any two vertices $x$ and $y$ in $E^0$, one has that $x\gg y$,
  \zitem $E$ is \"{weakly $G$-transitive} if, given any infinite path $\xi $, and any vertex $x$ in $E^0$, there is some
vertex $v$ along $\xi $ such that $v\gg x$.

The notion of $G$-transitivity generalizes the well known notion of transitivity in graphs.  When it holds, $E^0$ has a
single equivalence class.

On the other hand, weak $G$-transitivity is inspired by the notion of cofinality introduced in \cite [Section 3]{KPRR}, (see also \cite [Definition 37.16]{book}).  The reader is however warned that the notions of weak $G$-transitivity
and cofinality may only be reconciled upon a reversal of the direction of the edges in $E^1$, following the new trend in
graph algebras started by Katsura (see the penultimate paragraph of the introduction in \cite {KatsuraFundRes}).

It is evident that every $G$-transitive graph is weakly $G$-transitive, but these are sometimes equivalent notions as we
will now show:

\state Proposition \label WeakVsTransitive
  In addition to the assumptions in \ref{StandingHyp}, suppose that $E$ has no \"{sinks}, meaning that $\src \inv
(x)$ is nonempty for every $x$ in $E^0$.  Then, if $E$ is weakly $G$-transitive, it must also be $G$-transitive.

\Proof
  Since $E^0$ is finite, we may choose a minimal element $[x]$ in $\quoapprox $.  Using that $E$ has no sinks, we may
find an infinite sequence of edges
  $$
  \ldots ,\alpha _{-i-1},\alpha _{-i},\ldots ,\alpha _{-2},\alpha _{-1},\alpha _0 \in E^1,
  $$
  such that $\src (\alpha _0)=x$, and $\src (\alpha _{i-1}) = \ran (\alpha _i)$, for every $i\leq 0$.  Since $E^1$ is also finite, there must
be repetitions among the $\alpha _i$, say $\alpha _m = \alpha _n$, where $m<n\leq 0$.  The finite path
  $$
  \gamma =\alpha _m\alpha _{m+1}\ldots \alpha _{n-1},
  $$
  therefore satisfies
  $$
  \src (\gamma ) = \src (\alpha _{n-1}) = \ran (\alpha _n) = \ran (\alpha _m) = \ran (\gamma ),
  $$
  and hence $\gamma $ may be concatenated with itself infinitely many times producing the infinite path
  $$
  \xi =\gamma \gamma \gamma \ldots
  $$

  Given any $y$ in $E^0$, and assuming weak $G$-transitivity, there is some vertex $v$ along $\xi $, such that $v\gg y$.
Since $\xi $ is made of repetitions of $\gamma $, one has that $v=\ran (\alpha _k)$, for some $k$ in the integer interval $[m,n]$. We
then have
  $$
  x = \src (\alpha _0) =
  \src (\alpha _{k}\ldots \alpha _{-2}\alpha _{-1}\alpha _0) \rightharpoonup \ran (\alpha _k) = v \gg y,
  $$
  so $x\gg y$, but since $[x]$ is minimal, we deduce that $[x]=[y]$, which is to say that $x\approx y$.

The conclusion is that  $\quoapprox $ is a singleton, from where $G$-transitivity follows.
  \endProof

Of course the above result has taken advantage of the fact  that $E$ is a finite graph in an essential way, so nothing
like this is to be expected for infinite graphs.

Regardless of the absence of sinks, we have:

\state Theorem \label CharacMinimal
  Given $(G,E,\varphi )$ satisfying \ref{StandingHyp}, one has that the following are equivalent:
  \izitem
  \zitem the standard action of $\SGE $ on $E^\infty $ defined in \ref{ActionOfSGE} is irreducible,
  \zitem $\GpdGE $ is minimal,
  \zitem $E$ is weakly $G$-transitive.

\Proof
  The equivalence between (i) and (ii) is a consequence of \cite[Proposition 5.4]{EPFour}.  We will next show that the above
condition (iii) is equivalent to condition (iii) of \cite[Theorem 5.5]{EPFour}, from where the result will follow.  In
doing so, it is useful to understand how do idempotents in $\EGE $ behave under conjugation by elements in $\SGE $, and
we leave it for the reader to verify that, given $(\alpha ,g,\beta )$ in $\SGE $ and $(\gamma ,1,\gamma )\in \EGE $, one has that
  $$\def \quad {\ }
  (\alpha ,g,\beta )(\gamma ,1,\gamma )(\alpha ,g,\beta )^* = \left \{\matrix {
    (\alpha g\varepsilon ,1,\alpha g\varepsilon ), & \hbox { if } \gamma =\beta \varepsilon , \cr \pilar {15pt}
    (\alpha ,\hfil 1,\hfil \alpha ), & \hbox { if } \gamma \varepsilon =\beta , \cr \pilar {15pt}
    0, & \hbox { otherwise }.
    }\right .
  \equationmark Conjuga
  $$

\noindent (iii)$\Rightarrow $\cite[Theorem 5.5.iii]{EPFour}:
  Given any two nonzero idempotent elements in $\EGE $, necessarily of the form
  $$
  f_\alpha =(\alpha ,1,\alpha ) \and f_\beta =(\beta ,1,\beta ),
  $$
  to employ the notation introduced in \ref{DefineEAlpha},
  we must find an outer cover of $f_\alpha  $ (in the sense of \cite[Definition 2.9]{EPFour}) formed by a finite number of conjugates of $f_\beta $.
  As a first step, notice that $s:=\big(\src (\beta  ),1,\beta  \big)$ lies in $\SGE $ and
  $$
  sf_\beta  s^* = \big(\src (\beta  ),1,\beta  \big)\ (\beta  ,1,\beta  )\ \big(\beta  ,1,\src (\beta  )\big) = \big(\src (\beta  ),1,\src (\beta  )\big) = f_{\src (\beta  )}.
  $$
  Thus,
  anything that  may be  obtained by conjugating $f_{\src (\beta  )}$ by an element $t\in  \SGE $, may also be obtained by
conjugating $f_\beta  $ by $ts$.  It therefore suffices to find an outer cover of $f_\alpha  $ formed by conjugates of $f_{\src (\beta
)}$.

  On the other hand, observe that $f_\alpha  \leq f_{\ran (\alpha  )}$, so any outer cover of $f_{\ran (\alpha  )}$ is necessarily also an
outer cover of $f_\alpha  $.
  This said we see that we may assume, without loss of generality, that $\alpha $ and $\beta $ are vertices.

Our task thus gets simplifyed in the sense that we now need to find an outer cover of $f_x$ made of conjugates of $f_y$,
for any given vertices $x$ and $y$ in $E^0$.

Recall from \ref{DefineCylinder} that the set of all infinite paths with a given prefix $\gamma $ is denoted $\cyl \gamma $.
In case we take $\gamma =x$, we then have that $\cyl x$ is the set of all infinite paths with range $x$.

Thanks to weak $G$-transitivity, for each $\xi $ in $\cyl x$, we may choose a vertex $v_\xi $ along $\xi $ such that
$v_\xi \gg y$.  This is to say that we may write $\xi  = \alpha _\xi \eta _\xi $, where $\alpha _\xi $ is a finite path, $\eta _\xi $ is an infinite path and
$\src (\alpha _\xi )=v_\xi $.

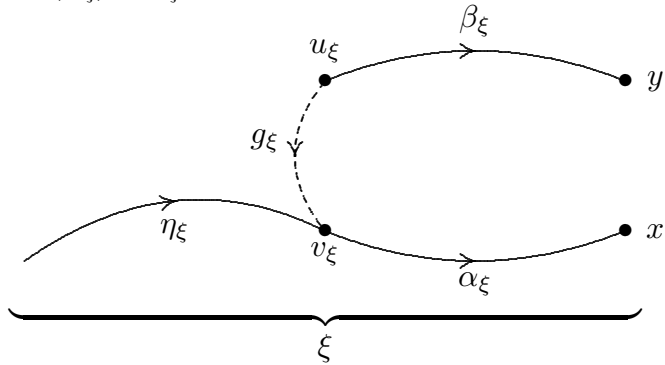
  \beginpicture
  \def \flexa {\arrow <0.17cm> [0.3,1]}
  \setcoordinatesystem units <0.0040truecm, 0.0020truecm> point at 3000 0
  \setplotarea x from -2000 to 2000, y from -1200 to 1000
  \setquadratic
  \put {$\bullet $} at 0 1000
  \put {$\bullet $} at 1000 1000 \put {$y$} at 1100 1000
  \setquadratic
  \put {$u_\xi $} at 0 1200
  \plot 0 1000 500 1200 1000 1000 / \flexa from 500 1200 to 501 1200 \put {$\beta _\xi $} at 500 1400
  \put {$\bullet $} at 0 0 \put {$v_\xi $} at 0 -150
  \put {$\bullet $} at 1000 0 \put {$x$} at 1100 0
  \plot -1000 -200 -500 200 0 0 / \flexa from -500 200 to -497 201 \put {$\eta _\xi $} at -500 0
  \plot 0 0 500 -200 1000 0 / \flexa from 500 -200 to 501 -200 \put {$\alpha _\xi $} at 500 -350
  \flexa from -100 500 to -100 499 \put {$g_\xi $} at -200 600
  \setdashes <2pt>
  \plot 0 1000 -100 500 0 0 /
  \put {$\underbrace {\hbox to 7cm{\hfill }}_{\textstyle \xi }$} at 0 -700
  \endpicture

The fact that $v_\xi \gg y$ may be expressed by saying that $v_\xi \sim u_\xi  \rightharpoonup y$, for some vertex $u_\xi $, so
there exists $g_\xi $ in $G$, and a finite path $\beta _\xi $, such that $g_\xi u_\xi =v_\xi $, $\src (\beta _\xi )=u_\xi $, and $\ran (\beta _\xi )=y$.

Speaking of the cylinders $\cyl {\alpha _\xi }$, it is obvious that $\xi \in \cyl {\alpha _\xi }$, so we see that the collection of cylinders
  $$
  \big \{\cyl {\alpha _\xi }\big \}_{\xi \in \cyl x}
  $$
  is an open cover (in the topological sense of the word) for $\cyl x$.  Since $\cyl x$ is compact, we may extract a
finite subcover, say
  $$
  \cyl x \subseteq  \medcup _{\xi \in F}\cyl {\alpha _\xi },
  \equationmark FiniteCoverForCyl
  $$
  where $F$ is a finite subset of $\cyl x$.  We next claim that $\{f_{\alpha _\xi }\}_{\xi \in F}$ is an outer cover\fn
    {This is in fact a cover but we do not need to worry about this right now.}
  of $f_x$.

  To see this, let $e$ be a nonzero idempotent in $\EGE $, with $e\leq f_x$.  Then $e$ is necessarily given by $e=(\gamma ,1,\gamma )$,
for some finite path $\gamma $ such that $\ran (\gamma )=x$.  Using our standing hypothesis \ref{StandingHyp} according to
which $E$ has no sources, we may prolong $\gamma $ to an infinite path $\eta $, which will then share ranges with $\gamma $, whence
$\eta \in \cyl x$.  By \ref{FiniteCoverForCyl} we then have that $\eta $ lies in $\cyl {\alpha _\xi }$, for some $\xi \in F$.

This implies that both $\alpha _\xi $ and $\gamma $ are prefixes of $\eta $, from where it is easy to see that $\alpha _\xi $ is a prefix of $\gamma $ or
vice-versa.  In particular we conclude that $f_{\alpha _\xi }\Cap f_\gamma $, proving our claim.  Incidentally this could also be
obtained from \cite[Proposition 3.8]{EPFour}.

We next claim that each $f_{\alpha _\xi }$ is a conjugate of $f_y$.  To see this, observe that, since
  $$
  g_\xi \src (\beta _\xi )= g_\xi u_\xi  = v_\xi  = \src (\alpha _\xi ),
  $$
  one has that $s_\xi :=(\alpha _\xi ,g_\xi ,\beta _\xi )$ lies in $\SGE $, and
  $$
  s_\xi  f_y s_\xi ^* =
  (\alpha _\xi ,g_\xi ,\beta _\xi )(y,1,y) (\beta _\xi ,g_\xi ,\alpha _\xi ) =
  (\alpha _\xi ,1,\alpha _\xi ) = f_{\alpha _\xi }.
  $$
  This concludes the proof of condition (iii) of \cite[Theorem 5.5]{EPFour}.

\bigskip \noindent \cite[Theorem 5.5.iii]{EPFour}$\Rightarrow $(iii):
  Given any infinite path $\xi $, and any vertex $y$ in $E^0$, we must show that there is some vertex $v$ along $\xi $ such
that $v\gg y$.  Letting $x=\ran (\xi )$, let us use the hypothesis regarding the nonzero idempotents
  $$
  f_x=(x,1,x)\and f_y = (y,1,y).
  $$
  This is to say that there are $s_1,s_2,\ldots ,s_n$ in $\SGE $, such that $\{s_if_ys_i^*\}_{1\leq i\leq n}$ is an outer cover for $f_x$.
  For each $i$, write $s_i=(\alpha _i,g_i,\beta _i)$, so that
  $$
  s_if_ys_i^* = (\alpha _i,g_i,\beta _i)(y,1,y)(\beta _i,g_i,\alpha _i).
  $$

  Observe that, unless $\ran (\beta _i)=y$, the element displayed above vanishes, so it cannot possibly have any use as a
member of a cover.  We may then safely discard it, being left only with those $\beta _i$ such that that $\ran (\beta _i)=y$.  In
this case, by the second clause in \ref{Conjuga} we have
  $$
  s_if_ys_i^* = (\alpha _i,1,\alpha _i) = f_{\alpha _i}.
  $$

  Unless $\ran (\alpha _i)=x$, notice that $f_{\alpha _i} \perp  f_x$, in which case $f_{\alpha _i}$ may again be discarded as it plays no role
in an outer cover for $f_x$.  We may therefore suppose, without loss of generality that $\ran (\alpha _i)=x$, for all $i$.

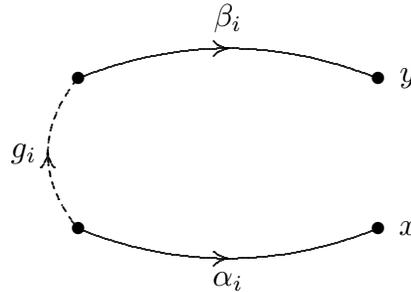
  \beginpicture
  \def \flexa {\arrow <0.17cm> [0.3,1]}
  \setcoordinatesystem units <0.0040truecm, 0.0020truecm> point at 3000 0
  \setplotarea x from -1500 to 3000, y from -500 to 1000
  \setquadratic
  \put {$\bullet $} at 0 1000
  \put {$\bullet $} at 1000 1000 \put {$y$} at 1100 1000
  \setquadratic
  \plot 0 1000 500 1200 1000 1000 / \flexa from 500 1200 to 501 1200 \put {$\beta _i$} at 500 1400
  \put {$\bullet $} at 0 0
  \put {$\bullet $} at 1000 0 \put {$x$} at 1100 0
  \plot 0 0 500 -200 1000 0 / \flexa from 500 -200 to 501 -200 \put {$\alpha _i$} at 500 -350
  \flexa from -100 500 to -100 501 \put {$g_i$} at -180 500
  \setdashes <2pt>
  \plot 0 1000 -100 500 0 0 /
  \endpicture

Given that $(\alpha _i,g,\beta _i)$ lies in $\SGE $, we necessarily have that $\src (\alpha _i)=g_i\src (\beta _i)$.
  Recalling that the infinite path $\xi $, chosen at the beginning of the present argument, has range $x$,
  we claim that $\xi $ is necessarily of the form
  $$
  \xi =\alpha _i\xi ',
  $$
  for some $i$ and some infinite path $\xi '$.  To see this, write
  $$
  \xi =\delta \xi '',
  $$
  where $\xi ''$ is an infinite path and $\delta $ is a finite path whose length exceeds the length of all of the $\alpha _i$.
Observing that $\ran (\delta )=x$, we then have that
  $$
  f_\delta  := (\delta ,1,\delta ) \leq  (x,1,x) = f_x.
  $$

  So, by the covering property we must have $f_\delta \Cap f_{\alpha _i}$, for some i, which implies that either
  $\delta $ is a prefix of $\alpha _i$ or vice versa.  However, due to the fact that $|\delta |>|\alpha _i|$, by construction, the first
alternative cannot hold, meaning that
  $\alpha _i$ is a prefix of $\delta $, and hence also of $\xi $, proving the claim.

It follows that $\src (\alpha _i)$ is a vertex along $\xi $, and it is clear from the above diagram that $\src (\alpha _i)\gg y$.  This
concludes the proof.  \endProof

Combining the above result with \ref{WeakVsTransitive}, we immediately deduce:

\state Corollary
  If, in addition to the assumptions of \ref{CharacMinimal} we have that $E$ has no sinks, then conditions
  \ref{CharacMinimal.i--iii} are also equivalent to:
  \izitem \zitemno 3
  \zitem $E$ is $G$-transitive.

It is interesting to observe that $G$-transitivity, when it holds, is the result of a joint effort by the action of $G$
and the edges, both of which may be seen as pushing vertices around.  However, sometimes only one of the players bear
the responsibility to do the pushing around:
  \initem
  \nitem If $G$ acts transitively on $E^0$, then $E$ is $G$-transitive regardless of the graph. Easy examples of this
situation may be built on a graph formed by a disjoint union of loops, for instance.
  \nitem If $G$ fixes all vertices, then $E$ is (weakly) $G$-transitive if and only if $E$ is (weakly) transitive
\cite [Definition 37.16]{book}. This is the case of Katsura algebras, when seen in the present framework.

\section Essentially principal groupoids

In this section we will discuss conditions under which $\GpdGE$ is an essentially principal groupoid, a condition which
is intimately tied  to the action of $\SGE$ on $E^\infty $ being topologically free.  The reader is referred to
\cite[Section 4]{EPFour} for the definition of the notion of topologically free actions of inverse semigroups, as well as
some of the main tools we shall use here.

\definition
  \initem
  \nitem A \"{circuit}\fn
    {Circuits are also called loops or cycles in the graph C*-algebra literature.  Our preference for \"{circuits} comes
    from the fact that it is the terminology of choice among graph theorists and also because in the established graph
    theory terminology the word \"{loop} refers to a single edge whose source and range coincide.}
   is a finite path $\gamma \in E^*$ of nonzero length such that $\src (\gamma )= \ran (\gamma )$.
  \nitem A \"{$G$-circuit} is a pair $(g,\gamma )$, where $g\in G$, and $\gamma \in E^*$ is a finite path of nonzero length such that $\src
(\gamma )= g\ran (\gamma )$.

\null \hfill \beginpicture \setcoordinatesystem units <0.009truecm, 0.009truecm> point at 1000 1000 \setplotarea x from
-300 to 300, y from -400 to 100
  \put {$\bullet $} at 66 240 \arrow <0.15cm> [0.4,1] from 250 0 to 148 132 \plot 148 132 66 240 / \put {$\gamma _6$} at
175 153
  \put {$\bullet $} at -159 184 \arrow <0.15cm> [0.4,1] from 66 240 to -57 209 \plot -57 209 -159 184 / \put {$\gamma
_5$} at -65 242
  \put {$\bullet $} at -254 -29 \arrow <0.15cm> [0.4,1] from -159 184 to -211 66 \plot -211 66 -254 -29 / \put {$\gamma
_4$} at -242 80
  \put {$\bullet $} at -105 -227 \arrow <0.15cm> [0.4,1] from -254 -29 to -172 -137 \plot -172 -137 -105 -227 / \put
{$\gamma _3$} at -199 -158
  \put {$\bullet $} at 145 -241 \arrow <0.15cm> [0.4,1] from -105 -227 to 32 -234 \plot 32 -234 145 -241 / \put {$\gamma
_2$} at 30 -268
  \put {$\bullet $} at 317 26 \arrow <0.15cm> [0.4,1] from 145 -241 to 239 -94 \plot 239 -94 317 26 / \put {$\gamma _1$}
at 268 -112
  \setdashes <1pt>
  \put {$\bullet $} at 250 0 \arrow <0.15cm> [0.4,1] from 317 26 to 280 11 \plot 280 11 250 0 / \put {$g$} at 267 43
\put {\eightrm A $\scriptstyle G$-circuit.} at 0 -330 \endpicture \hfill \null

Thus, a $G$-circuit needs a little help from the group to close it up.  Notice that a (usual) circuit $\gamma $ may be
concatenated with itself infinitely many times producing an infinite path
  $$
  \xi =\gamma \gamma \gamma \ldots
  $$
  Moreover, if
  $$
  s = \big(\gamma ,1,\src (\gamma )\big),
  $$
  then, regarding the standard action of $\SGE $ on $E^\infty $ defined in \ref{ActionOfSGE}, it is easy to see that
  $
  s\xi =\xi ,
  $
  which is to say that $\xi $ is a fixed point for $s$.  It is also possible to create fixed points from $G$-circuits as
follows:

\state Proposition \label GcircuitGivesFixed
  Given a $G$-circuit $(g,\gamma )$, define
  a sequence $\{\gamma ^n\}_{n\geq 1}$ of finite paths,
  and
  a sequence $\{g_n\}_{n\geq 1}$ of group elements,
  recursively by $\gamma ^1=\gamma $, \ $g_1=g$, and
  $$
  \def \quad {\ \ }
  \left \{\matrix {
  \gamma ^{n+1} &=& g_n\gamma ^n \hfill \cr \pilar {12pt}\stake {8pt}
  g_{n+1} &=& \varphi (g_n,\gamma ^n),
  }\right .
  $$
  for all $n\geq 1$.
  Then
  \izitem
  \zitem $\src (\gamma ^n)=\ran (\gamma ^{n+1})$, for all $n\geq 1$,
  \zitem the concatenation \ $\xi =\gamma ^1\gamma ^2\gamma ^3\ldots $ \ is a well defined infinite path,
  \zitem for every finite path $\beta $ such that $\src (\beta )=\ran (\gamma )$, one has that $s:= (\beta \gamma ,g,\beta )$ lies in $\SGE $, and $\beta \xi $
is a fixed point for $s$.

\Proof In order to prove the case $n=1$ of (i), we have
  $$
  \src (\gamma ^1)=
  \src (\gamma )=
  g\ran (\gamma ) =
  \ran (g\gamma ) =
  \ran (g_1\gamma ^1) =
  \ran (\gamma ^2).
  $$
  For $n\geq 1$, we have
  $$
  \src (\gamma ^{n+1}) =
  \src (g_n\gamma ^n) =
  g_n\src (\gamma ^n) \={Equacoes.vii}
  \varphi (g_n,\gamma ^n)\src (\gamma ^n) =
  g_{n+1}\src (\gamma ^n) \explica={${\scriptstyle(\star)}$} $$$$=
  g_{n+1}\ran (\gamma ^{n+1}) =
  \ran (g_{n+1}\gamma ^{n+1}) =
  \ran (\gamma ^{n+2}),
  $$
  where we have used induction in the step marked with ${\scriptstyle(\star)}$ above.  This proves (i), which in turn immediately implies
(ii).

In order to show that the element $s$ defined in (iii) indeed lies in $\SGE $, it is enough to observe that
  $$
  g\src (\beta ) = g\ran (\gamma ) = \src (\gamma ) = \src (\beta \gamma ).
  $$

Before proving the last part of (iii), we claim that
  $$
  g_n(\gamma ^n\gamma ^{n+1}\ldots \gamma ^{n+k}) = \gamma ^{n+1}\gamma ^{n+2}\ldots \gamma ^{n+k+1} \for n\geq 1 \for k\geq 0.
  $$
  In case $k=0$, this is true by the recursive definition above, and if $k\geq 1$, we have
  $$
  g_n(\gamma ^n\gamma ^{n+1}\ldots \gamma ^{n+k}) =
  (g_n\gamma ^n)\varphi (g_n,\gamma ^n)(\gamma ^{n+1}\ldots \gamma ^{n+k}) \$=
  \gamma ^{n+1}g_{n+1}(\gamma ^{n+1}\ldots \gamma ^{n+k}),
  $$
  and the claim then follows easily by induction.  A useful consequence is that
  $$
  g_1(\gamma ^1\gamma ^2\ldots \gamma ^n) = \gamma ^2\gamma ^3\ldots \gamma ^{n+1} \for n\geq 1,
  $$
  from where we further deduce that
  $$
  g\xi  =
  g_1(\gamma ^1\gamma ^2\gamma ^3\ldots ) = \gamma ^2\gamma ^3\gamma ^4\ldots
  \equationmark ActGXi
  $$

With this we may now tackle the final task:
  $$
  s(\beta \xi ) = (\beta \gamma ,g,\beta )(\beta \xi ) \={ActionOfSGE}
  \beta \gamma g\xi  \={ActGXi}
  \beta \gamma \gamma ^2\gamma ^3\ldots \gamma ^{n+1} =
  \beta \xi .
  \endProof

The above method does not give us all fixed points of every single element $s$ in $\SGE $, but in certain cases it does:

\state Proposition \label DescribeSomeFixedPoints
  Given $s := (\alpha ,g,\beta )$ in $\SGE $, suppose that $|\alpha |>|\beta |$.  Then, regarding the standard action of $\SGE $ on $E^\infty $, one has
that:
  \izitem
  \zitem $s$ admits at most one fixed point,
  \zitem if $s$ admits a fixed point $\zeta $, then there is a $G$-circuit $(g,\gamma )$ such that $\alpha =\beta \gamma $, and $\zeta $ coincides with
the fixed point $\beta \xi $ mentioned in \ref{GcircuitGivesFixed.iii}, constructed from $(g,\gamma )$.

\Proof
  Assuming that $\zeta $ is a fixed point for $s$, it must lie in $\cyl \beta $, so necessarily $\zeta =\beta \xi $, for a suitable infinite
path $\xi $.  We then have
  $$
  \beta \xi  = \zeta  = s\zeta  = (\alpha ,g,\beta )(\beta \xi ) = \alpha g\xi .
  \equationmark FromFixPoint
  $$

  This imples that both $\alpha $ and $\beta $ are prefixes of $\zeta $, so one must be a prefix of the other, but since $|\alpha |>|\beta |$, the
only alternative is that $\beta $ is a prefix of $\alpha $.  We may therefore write
  $$
  \alpha  = \beta \gamma ,
  $$
  for some finite path $\gamma $, which necessarily satisfies
  $$
  g\ran (\gamma ) =g\src (\beta ) = \src (\alpha ) = \src (\gamma ).
  $$
  In other words, $(g,\gamma )$ is a $G$-circuit.  From \ref{FromFixPoint} we also deduce that
  $$
  \beta \xi  = \alpha g\xi  = \beta \gamma g\xi ,
  $$
  so
  $
  \xi =\gamma g\xi .
  $
  Let us now write $\xi =\gamma ^1\gamma ^2\gamma ^3\ldots $, where each $\gamma ^i$ is a finite path with $|\gamma ^i|=|\gamma |$.  Then
  $$
  \gamma ^1\gamma ^2\gamma ^3\ldots = \xi = \gamma g\xi = \gamma g(\gamma ^1\gamma ^2\gamma ^3\ldots ) = \gamma (g_1\gamma ^1)(g_2\gamma ^2)(g_3\gamma ^3)\ldots ,
  $$
  where the $g_i$ are recursively defined by $g_1=g$, and $g_{n+1} = \varphi (g_n,\gamma ^n)$.  It then follows that
  $\gamma ^1=\gamma $, and $\gamma ^{n+1} = g_n\gamma ^n$, for all $n\geq 1$, so we see that the $\gamma ^n$ and the $g_n$ are precisely defined as in
\ref{GcircuitGivesFixed}.  This concludes the proof.
  \endProof

As already announced  we will eventually be interested in determining conditions under which the
standard action of $\SGE $ on $E^\infty $ is topologically free, so the fixed points that will really interest us are the
interior ones.

Under the conditions of the above result, when there is at most one fixed point, the existence of
interior fixed points hinges on whether or not the unique fixed point is isolated in $E^\infty $.  We will  now introduce
certain concepts designed  to study isolated fixed points.

Recall from  \ref{StandingHyp} that our graph $E$ has no sources, meaning that $\ran\inv(x)$ is nonempty for every
vertex $x$.

\definition
  \initem
  \nitem We shall say that a vertex  $x$ in $E^0$ is a \"{simple vertex} if $\ran\inv(x)$ is a singleton.
  \nitem Given a path $\gamma =\gamma _1\gamma _2\ldots \gamma _n$ in $E^*$, where each $\gamma _i$ is  in $E^1$, we will say that $\gamma $ \"{has no entry} if
$\src (\gamma _i)$ is a simple vertex for every $i=1,\ldots ,n$.
  \nitem If the condition above fails, we will say that $\gamma $ \"{has an entry}.

Thus, if a path $\gamma =\gamma _1\gamma _2\ldots \gamma _n$ has no entry, then $\ran\inv\big(\src (\gamma _i)\big)$ is a singleton for every $i$, and we may
obviously guess which is the edge forming this singleton, namely
  $$
  \ran\inv\big(\src (\gamma _i)\big) = \{\gamma _{i+1}\},
  $$
  provided $i<n$.  However the same cannot be said when  $i=n$,  unless $(g,\gamma )$ is a
$G$-circuit, in which case
  $$
  \ran \inv \big(\src (\gamma _n)\big) = \{g\gamma _1\}.
  $$
  The notion of entryless paths will only be  useful when applied to $G$-circuits.

\state Proposition \label IsolatedIffNoEntry
  Under the conditions of \ref{DescribeSomeFixedPoints.ii}, let $(g,\gamma )$ be the $G$-circuit and $\zeta $ be the fixed
point for $s$ mentioned there.  Then the following are equivalent:
  \izitem
  \zitem $\zeta $ is an isolated point in $E^\infty $,
  \zitem $\gamma $ has no entry.

\Proof
  In case $\gamma $ has no entry, writing $\gamma =\gamma _1\gamma _2 \ldots  \gamma _n$, where the $\gamma _i$ are edges, notice that the only infinite path extending $\gamma _1$ is the path
$\xi $ referred to in \ref{GcircuitGivesFixed.ii}.  The fixed point $\zeta =\beta \xi $ mentioned in \ref{DescribeSomeFixedPoints.ii} is therefore the only infinite path extending $\beta \gamma _1$, whence
  $$
  \cyl {\beta \gamma _1} = \{\zeta \},
  $$
  which implies that $\zeta $ is isolated.

Conversely, assuming that $\zeta $ is isolated, there exists a sufficiently long prefix $\varepsilon $ of $\zeta $, such that
  $$
  \cyl \varepsilon  = \{\zeta \}.
  $$
  This means that $\zeta $ is the only infinite path extending $\varepsilon $.  Writing
  $$
  \zeta =\zeta _1\zeta _2\zeta _3\ldots ,
  $$
  where the $\zeta _i$ are edges,   one then has that, for sufficiently large $i$, there is only one edge whose range is $\src (\zeta _i)$.  Letting
  $\{\gamma ^n\}_{n\geq 1}$ and $\{g_n\}_{n\geq 1}$ be the sequences defined in \ref{GcircuitGivesFixed}, we then have that, for
sufficiently large $n$, the $G$-circuit $(g_n,\gamma ^n)$ has no entry.  Since $G$ acts on $E$ by graph automorphisms, we may
easily prove by induction that all $G$-circuits $(g_k,\gamma ^k)$ have no entry, including $(g_1,\gamma ^1)=(g,\gamma )$.  \endProof

Since we are interested in topologically free actions, we would like to avoid isolated fixed points and hence we will be
interested in situations when every $G$-circuit has an entry.  However, given that we are working with finite graphs
only, the action of $G$ on $E$  turns out not to be relevant in this respect.  In precise terms, what we mean is that:

\state Proposition \label TudoIgual
  Under the conditions of \ref{StandingHyp}, the following are equivalent:
  \izitem
  \zitem every $G$-circuit has an entry,
  \zitem every circuit  has an entry.

\Proof
  Since every circuit $\gamma$ gives rise to the $G$-circuit $(1,\gamma)$, it is evident that (i) implies (ii).
  Conversely, assume (ii) and let $\gamma$ be a $G$-circuit.
Leting  $\{\gamma^n\}_{n\geq1}$  and $\{g_n\}_{n\geq1}$ be as in
\ref{GcircuitGivesFixed}, consider the infinite path
$\xi=\gamma^1\gamma^2\gamma^3\ldots$ mentioned in \ref{GcircuitGivesFixed.iii}.  Notice that the $\gamma ^i$ are all in the
orbit of $\gamma $ under the action of $G$, and hence the length of $\gamma ^i$ coincides with that of $\gamma $.   As $E$ is finite,
there is only  a finite number of paths of this length, so there must necessarily be repetitions among the $\gamma ^i$, say
$\gamma ^i=\gamma ^j$, where $i<j$.
  Then
  $$
  \rn{\gamma ^{i+1}} = \sr {\gamma ^i} =   \sr {\gamma ^j},
  $$
  so the path
  $$
  \gamma ^{i+1}\gamma ^{i+2}\ldots \gamma ^j
  $$
  is a circuit, which by hypothesis has an entry.  It is now easy to see that some $\gamma ^k$ must have an entry.  Finally,
since $\gamma ^k$ is in the orbit of $\gamma $ under the action of $G$, then $\gamma $ likewise has an entry, concluding
the proof.
  \endProof

Observe that we have used the finiteness of $E$ in a very strong way above.  Thus, should our theory ever be extended to
infinite graphs, one might have to distinguish between conditions \ref{TudoIgual.i} and \ref{TudoIgual.ii}.

We should point out that a graph in which every circuit has an entry is usually said to satisfy condition (L).

The above results, mainly \ref{DescribeSomeFixedPoints} and \ref{IsolatedIffNoEntry}, may also be used to study the
fixed points for elements $s:=(\alpha ,g,\beta )$ when $|\alpha |<|\beta |$, since such fixed points are precisely the same as the fixed
points of $s^*$, and $s^*$ clearly satisfies the hypothesis of \ref{DescribeSomeFixedPoints}.  However we still have work
to do in order to treat the remaining case $|\alpha |=|\beta |$.

\state Proposition \label FixForG
  Let $s:=(\alpha ,g,\beta )\in \SGE $, whith $|\alpha |=|\beta |$, and suppose that $s$ admits a fixed point.  Then
  \izitem
  \zitem $\alpha =\beta $,
  \zitem the fixed points of $s$ in $E^\infty $ are precisely the elements of the form $\zeta =\beta \xi $, where $\xi $ is an infinite path
such that $\ran (\xi ) = \src (\beta )$, and $g\xi =\xi $.

\Proof Left for the reader.  \endProof

The conclusion of the previous Proposition is that when $|\alpha |=|\beta |$, understanding the fixed points for $s$ requires
understanding the fixed points for the action of $g$ on $E^\infty $.  One may easily describe such fixed points in terms of
the action of $G$ on $E$ and the cocycle $\varphi $, but apparently there is no smart way to control each and every one of
them.  However, since our main interest is in studying topological freeness, we need only focus on large (meaning open)
sets of fixed points:

\state Proposition \label GIntFixPoints
  Suppose that  $s:=(\alpha ,g,\alpha )$ lies in $\SGE $, and  that $\zeta $ is an interior fixed point for $s$.  Then there is a
finite path $\gamma $, such that:
  \izitem
  \zitem $g\gamma =\gamma $,
  \zitem $\src (\alpha ) = \ran (\gamma )$,
  \zitem $\zeta \in \cyl {\alpha \gamma }$,
  \zitem the group element $h: = \varphi (g,\gamma )$ pointwise fixes\fn
    {By this we mean that every point in $\cyl {\src (\gamma )}$ is fixed by $h$.}
  the cylinder $\cyl {\src (\gamma )}$.
  \medskip \noindent Conversely, if $\gamma $ is any finite path satisfying (i), (ii) and (iv), then every $\zeta \in \cyl {\alpha \gamma }$ is a
(necessarily interior) fixed point for $s$.

\Proof
   In particular $\zeta $ a fixed point for $s$ so, by \ref{FixForG} we have
that $\zeta =\alpha \xi $, with $g\xi =\xi $.

Moreover there exists a neighborhood $U$ of $\zeta $ consisting of fixed points for $\zeta $.  Since the cylinders form a basis
for the topology of $E^\infty $, we may assume without loss of generality that $U=\cyl \beta $, for some finite path $\beta $, which we
may assume is as long as we wish, and our wish in this case is simply that $|\beta |>|\alpha |$.

Since $\zeta $ lies in $\cyl \beta $, we have that $\beta $ is a prefix of $\zeta $, so we may write $\zeta =\beta \eta $, for some infinite path $\eta $.  We
then have
  $$
  \beta \eta  = \zeta  = \alpha \xi .
  $$

  Given that $|\beta |>|\alpha |$,
  this implies that $\alpha $ is a prefix of $\beta $, so we write $\beta  = \alpha \gamma $, for a suitable finite path $\gamma $, obviously satisfying
(ii).  Consequently
  $$
  \zeta  = \beta \eta  = \alpha \gamma \eta \in \cyl {\alpha \gamma },
  $$
  proving (iii).  Given any infinite path $\mu \in \cyl {\src (\gamma )}$, we may form the path $\alpha \gamma \mu $, which necessarily lies in
$\cyl {\alpha \gamma }=\cyl \beta $, and hence is fixed under $s$.  Therefore
  $$
  \alpha \gamma \mu  =
  (\alpha ,g,\alpha )(\alpha \gamma \mu ) =
  \alpha g(\gamma \mu ) =
  \alpha \ (g\gamma )\ \big(\varphi (g,\gamma )\mu \big) =
  \alpha \ (g\gamma )\ (h\mu ),
  $$
  whence $\gamma =g\gamma $, proving (i), and $\mu =h\mu $, in turn proving (iv).

In order to prove the last sentence in the statement it is enough to notice that any element in $\cyl {\alpha \gamma }$ is
necessarily of the form $\alpha \gamma \mu $, where $\mu \in \cyl {\src (\gamma )}$, and the last calculation displayed above could be used to
check that $\alpha \gamma \mu $ is fixed under $s$.
  \endProof

Searching for conditions under which the standard action of $\SGE $ on $E^\infty $ is topologically free, one should probably
worry about group elements fixing whole cylinders, as in \ref{GIntFixPoints.iv}.
  The following notion is designed to pinpoint situations under which whole cylinders of the form $\cyl x$ are in fact fixed.

\definition \label DefineSlack
   Given $g\in G$, and $x\in E^0$, we shall say that $g$
   is \"{slack} at $x$, if there is a non-negative integer $n$ such that
   all finite paths $\gamma $ with $\ran (\gamma )=x$, and $|\gamma |\geq n$, are strongly fixed by $g$, as defined in \ref{StrFixed}.

As already discussed at the begining of section \ref{PseudoFreeSect}, if $\gamma $ is strongly fixed by $g$, then $g$ fixes
any finite path extending $\gamma $, and hence also all infinite paths in $\cyl \gamma $.

If $g$ is slack at $x$, and if $n$ is as in \ref{DefineSlack}, notice that
  $$
  \cyl x = \kern -6pt \medcup _{{\buildrel \scriptstyle \ran (\gamma )=x \over {\ |\gamma |=n}}} \kern -6pt \cyl \gamma ,
  $$
  and since each $\gamma $ occuring above is strongly fixed by $g$, we have that $g$ pointwise fixes $\cyl \gamma $, and hence also
the whole cylinder $\cyl x$.

Notice that a path of length zero, namely a vertex $x$, is never strongly fixed by a nontrivial group element $g$, because
  $$
  \varphi (g,x)\={Equacoes.\CocZero} g \neq 1.
  $$

  The concept  of slackness above should therefore be seen as the best replacement for the notion of being  strongly
fixed in case of a vertex.

\medskip
We are now ready for a main result:

\state Theorem \label MainTopFreeSGE
  Under the conditions of \ref{StandingHyp},
  the standard action of\/ $\SGE $ on $E^\infty $ is topologically free if and only if the following two conditions hold:
  \izitem
  \zitem every $G$-circuit has an entry\fn{Recall that this  is the same as saying that every circuit has an entry by \ref{TudoIgual}.},
  \zitem given a vertex $x$, and a group element $g$ fixing every infinite path in $\cyl x$, then necessarily  $g$ is
slack at $x$.

\Proof Suppose (i) and (ii) hold and let $\zeta $ be an interior fixed point for some $s = (\alpha ,g,\beta )$ in $\SGE $.  In order to
prove topological freeness, we need to prove that $\zeta $ is a trivial fixed point for $s$.

\Case 1:
  Let us first assume that $|\alpha | = |\beta |$.  Letting $\gamma $ and $h$ as in \ref{GIntFixPoints}, we then have that $h$
pointwise fixes the cylinder $\cyl {\src (\gamma )}$.  By (ii) we then conclude that $h$ is slack at $\src (\gamma )$, so there is
$n$ such that every finite path of length $n$ and range $\src (\gamma )$ is strongly fixed by $h$.

  By \ref{GIntFixPoints.iii} we have that $\zeta $ lies in $\cyl {\alpha \gamma }$, so we may write $\zeta =\alpha \gamma \xi $, for some infinite path
$\xi $ with $\src (\gamma ) = \ran (\xi )$.
  Denoting by $\varepsilon $ the path formed by the first $n$ edges of $\xi $, we then have that
  $$
  \ran (\varepsilon ) = \ran (\xi ) = \src (\gamma ),
  $$
  so $\varepsilon $ is strongly fixed by $h$, and we may further write
  $$
  \zeta =\alpha \gamma \varepsilon \xi ',
  $$
  for a suitable infinite path $\xi '$.  If follows that $\zeta \in \cyl {\alpha \gamma \varepsilon }$, which is the domain of the idempotent
  $$
  f_{\alpha \gamma \varepsilon }=(\alpha \gamma \varepsilon ,1,\alpha \gamma \varepsilon ).
  $$
  In addition
  $$
  sf_{\alpha \gamma \varepsilon } =
  (\alpha ,g,\alpha )(\alpha \gamma \varepsilon ,1,\alpha \gamma \varepsilon ) =
  \big(\alpha g(\gamma \varepsilon ),\varphi (g,\gamma \varepsilon ),\alpha \gamma \varepsilon \big),
  \equationmark Calcsfage
  $$
  and we claim that the element at the end of the above calculation coincides with $f_{\alpha \gamma \varepsilon }$.  To see this notice that
  $$
  g(\gamma \varepsilon ) =
  (g\gamma )\big(\varphi (g,\gamma )\varepsilon \big) \={GIntFixPoints.i}
  \gamma  h \varepsilon  =
  \gamma  \varepsilon ,
  $$
  while
  $$
  \varphi (g,\gamma \varepsilon ) \={Equacoes.x}
  \varphi \big(\varphi (g,\gamma ),\varepsilon \big) =
  \varphi \big(h,\varepsilon \big) = 1.
  $$

  Plugging the last two identities at the end of \ref{Calcsfage} leads to
  $
  sf_{\alpha \gamma \varepsilon } = f_{\alpha \gamma \varepsilon },
  $
  thus proving that $\zeta $ is a trivial fixed point, as needed.

\Case 2:
  Let us now assume that $|\alpha | > |\beta |$.  By \ref{DescribeSomeFixedPoints} we have that $\zeta $ is the only fixed point
for $s$, necessarily given in terms of a $G$-circuit $(g,\gamma )$, as in \ref{DescribeSomeFixedPoints.iii}.

Being the unique fixed point, as well as an interior member of the set of fixed points, we see that $\zeta $ is isolated in
$E^\infty $.  So $(g,\gamma )$ has no entry by \ref{IsolatedIffNoEntry}, contradicting (i).  This implies that in fact $s$ has
no interior fixed points, so there is nothing to do.

\Case 3:
  The last remaining alternative, namely when $|\alpha | < |\beta |$, may be treated by simply observing that the fixed points for
$s$ are the same as the fixed points for $s^*=(\beta ,g\inv ,\alpha )$, and that $s^*$ fits the previous case studied, so there are
no interior fixed points for $s^*$, either.

\medskip This concludes the proof that (i) and (ii) imply topological freeness.  In order to prove that topological
freeness implies (i), assume the former and suppose by contradiction that a $G$-circuit $(g,\gamma )$ exists with no entry.
Let $x=\ran (\gamma )$, and notice that
  $$
  gx = g\ran (\gamma ) = \src (\gamma ),
  $$
  so the triple
  $
  s := (\gamma ,g,x)
  $
  is seen to lie in $\SGE $.  We may then use \ref{GcircuitGivesFixed} to obtain a fixed point $\zeta $ for $s$, and   by
\ref{IsolatedIffNoEntry} we have that $\zeta $ is an isolated point of $E^\infty $, hence also an interior fixed point.

Working under the assumption of topological freeness, we deduce that $\zeta $ is a trivial fixed point, which is to say
that there is an idempotent $e$ in $\EGE $, whose domain contains $\zeta $, and such that $se=e$.  Observing that $e$ cannot
possibly be zero, we deduce that
  $
  e = (\varepsilon ,1,\varepsilon ),
  $
  for some finite path $\varepsilon $.  We then have that
  $$
  (\varepsilon ,1,\varepsilon )= e = se = (\gamma ,g,x) (\varepsilon ,1,\varepsilon ) = \big(\gamma g\varepsilon ,\varphi (g,\varepsilon ),\varepsilon \big).
  $$
  In particular this implies that $\varepsilon =\gamma g\varepsilon $, so
  $$
  |\varepsilon | = |\gamma g\varepsilon | = |\gamma |+|g\varepsilon | = |\gamma |+|\varepsilon |,
  $$
  whence $|\gamma |=0$, contradicting the fact that $G$-circuits have nonzero length by definition.  This shows that there are
no $G$-circuit without an entry, hence proving (i).

We next show that topological freeness implies (ii).  So we suppose that some $g$ in $G$ pointwise fixes a whole
cylinder $\cyl x$, where $x$ is a vertex.  In particuler we have that $gx=x$, so the element
  $$
  s:= (x,g,x)
  $$
  belongs to $\SGE $, and it clearly also fixes every point in $\cyl x$.  Each $\zeta $ in $\cyl x$ is therefore an interior
fixed point for $s$, hence necessarily a trivial one by hypothesis.  This means that there exists an idempotent $e =
(\gamma ,1,\gamma )\in \E $, such that $\zeta $ lies in the domain of $e$, also known as $\cyl \gamma $, and moreover $se=e$.  Therefore
  $$
  (\gamma ,1,\gamma ) = e = se = (x,g,x) (\gamma ,1,\gamma ) = \big(xg\gamma ,\varphi (g,\gamma ),\gamma \big),
  $$
  from where we deduce that $g\gamma =\gamma $, and $\varphi (g,\gamma )=1$, which is to say that $\gamma $ is strongly fixed by $g$.

Given that $\zeta \in \cyl \gamma $, we have that $\gamma $ is a prefix of $\zeta $, whence $\ran (\gamma ) = \ran (\xi ) = x$, so
  $$
  \zeta \in \cyl \gamma \subseteq \cyl x.
  $$

  We then deduce that $\cyl x$ is the union of the $\cyl \gamma $, where $\gamma $ range in the set of all finite paths strongly
fixed by $g$, with $\ran (\gamma )=x$.  By compactness we may find a finite collection of such finite paths, say
$\gamma _1,\gamma _2,\ldots ,\gamma _k$, such that
  $$
  \cyl x = \medcup _{i=1}^k \cyl {\gamma _i}.
  \equationmark GamaiCoverX
  $$

  We next wish to argue that the above $\gamma _i$'s may be taken so that their length is constant.  To see this
  let $n$ be the length of the longer $\gamma _i$, and observe that, for each $i$, one has that
  $$
  \cyl {\gamma _i} = \kern -10pt \medcup _{{\buildrel \scriptstyle \ran (\varepsilon )=\src (\gamma _i) \over {|\varepsilon |=n-|\gamma _i|}}} \kern -6pt \cyl
{\gamma _i\varepsilon }.
  $$
  Moreover, each $\gamma _i\varepsilon $ occuring above is also strongly fixed by $g$, as seen in the discussion near the beginning of
section \ref{PseudoFreeSect}.  Thus, if we replace each $\gamma _i$ by the set of all $\gamma _i\varepsilon $, where $\varepsilon $ is as above, all of
the properties so far mentioned of the original $\gamma _i$'s will be preserved, and now
  $$
  |\gamma _i\varepsilon | = |\gamma _i| + |\varepsilon | = n.
  $$

Therefore we may and will assume, from now on, that the $\gamma _i$ have a constant length, say $n$.  From \ref{GamaiCoverX}
it is now easy to conclude that the $\gamma _i$ exhaust the set of all finite paths with range $x$ and length $n$.
In fact, if $\alpha $ is such a path, we may extend it to an infinite path of the form $\xi =\alpha \eta $.  Since $\xi \in \cyl x$, then $\xi \in \cyl
{\gamma _i}$, for some $i$, whence $\gamma _i$ is a prefix of $\xi $ and, by considering lengths, we see that $\alpha =\gamma _i$.

The conclusion is that every finite path with length $n$ and range $x$ is strongly fixed by $g$, which is to say that
$g$ is slack at $x$.  \endProof

\state Remark \label RemarkEssPrincTriv
\rm If for any $g\in G\setminus \{1\}$ and for any $x\in E^0$ there exists $\eta \in Z(x)$ such that $g\eta\ne \eta$, then \ref{MainTopFreeSGE.ii} holds trivially. This fact will be used in \ref{CorolKatsuraEssPrinc} and subsequent examples.

\state Remark \rm Regarding \cite[Theorem 4.10.ii]{EPFour}, and letting $\gamma _i$ be as in \ref{GamaiCoverX}, one may show
that $\{f_{\gamma _i}\}_i$ is a cover of $f_x$ consisting of idempotents fixed under $s$ (in the sense of \cite[Definition 4.8.1]{EPFour}).

In case $(G,E, \varphi )$ is pseudo free, and if $g$ is a nontrivial group element, then $g$ admits no strongly fixed paths
by \ref{EssFreePath}, so $g$ will never be slack at any vertex.  Condition
\ref{MainTopFreeSGE.ii} can therefore only be satisfied if no nontrivial group element pointwise fixes a cylinder $\cyl
x$, and hence we have the following immediate consequence of \ref{MainTopFreeSGE}:

\state Corollary \label TopFreeUnderEStarUni
  In addition to the conditions of \ref{StandingHyp}, suppose that $(G,E, \varphi )$ is pseudo free.  Then
  the standard action of\/ $\SGE $ on $E^\infty $ is topologically free if and only if the following two conditions hold:
  \izitem
  \zitem every $G$-circuit has an entry (which is the same as saying that every circuit has an entry by \ref{TudoIgual}),
  \zitem for every $g$ in $G$, with $g\neq 1$, and for every $x$ in $E^0$, there is at least one $\zeta $ in $\cyl x$ such that
$g\zeta \neq \zeta $.

An important case for the theory of self-similar groups is when $G$ acts faithfully\fn
  {Meaning that if $g\xi =\xi $, for all $\xi $ in $E^\infty $, then $g=1$.}
  on $E^\infty $, and $E$ is a graph with a single vertex.

\state Corollary \label TopFreeUnderSSGroup
  Under the conditions of \ref{StandingHyp}, suppose moreover that:
  \iaitem
  \aitem $E$ has a single vertex, and at least two edges,
  \aitem $G$ acts faithfully on $E^\infty $.
  \medskip \noindent Then
  the standard action of\/ $\SGE $ on $E^\infty $ is topologically free.

\Proof
  In the present situation the conditions of  \ref{MainTopFreeSGE} become trivially true because:
  (i) all path are circuits and all circuits have entries, and
  (ii) there is only one $\cyl x$ to consider, namely the whole space $E^\infty $, and by faithfulness no nontrivial group
element acts trivially on $E^\infty $.
  \endProof

As the title of the present section suggests, our main interest is in determining conditions for $\GpdGE$ to be an
essentially principal groupoid.  Having understood topological freeness, an immediate consequence of
\cite[Theorem 4.7]{EPFour} is:

\state Corollary \label TopFreeSameEssPrin
  Under the assumptions of \ref{StandingHyp}, one has that $\GpdGE$ is essentially principal if and only if
\ref{MainTopFreeSGE.i\&ii} hold.

Two other similar results  could be stated giving conditions for $\GpdGE$ to be essentially principal, by combining
\cite[Theorem 4.7]{EPFour} with either
\ref{TopFreeUnderEStarUni} or
\ref{TopFreeUnderSSGroup}, but we will refrain  from doing it here since the reader can easily guess them.

\section Local contractivity for $\SGE$

In \cite[Section 6]{EPFour} local contractivity for groupoids and for actions of inverse
semigroups is studied.  We will now use these results to characterize local contractivity for the tight groupoid
associated to an inverse semigroup $\S$.

\state Theorem \label LocContrVsTopFree
  Under the conditions of \ref{StandingHyp}, one has that the following are equivalent:
  \izitem
  \zitem $\SGE$ is a locally contracting inverse semigroup,
  \zitem the standard action $\theta :\SGE\curvearrowright E^\infty $ is locally contracting,
  \zitem $\GpdGE$ is a locally contracting groupoid,
  \zitem every circuit in $E$ has an entry.

\Proof As already mentioned in section \ref{TightGpdSectn}, every tight filter in $\EGE$ is an ultra-filter, so the
equivalence between (i) and (ii) follows from \cite[Theorem 6.5]{EPFour}.

\iskip (ii)$\Rightarrow$(iii): Follows immediately from \cite[Proposition 6.3]{EPFour}.

\iskip (iii)$\Rightarrow$(iv):
We will prove this by contraposition,  that is,
assuming  the existence of a circuit $\gamma $ without an entry, we will show that $\GpdGE$ is not locally contracting.

Our task is actually very easy. Given an entryless circuit $\gamma $, the path $\xi =\gamma \gamma \gamma \ldots $ is an isolated point, whence
$U:=\{\xi \}$ is an open subset of $E^\infty $.  Viewing the latter as the unit space of $\GpdGE$, as usual, and plugging $U$ into
\cite[Definition 6.1]{EPFour},  clearly there can be no open set $V$, and bissection $S$, as mentioned there, simply
because a chain of nonempty subsets
  $$
  S\clos VS\inv\subsetneqq V\subseteq  U
  $$
  cannot possibly exist withing a singleton such as $U$. This shows that $\GpdGE$ is not locally contracting, as desired.

\iskip
(iv)$\Rightarrow$(i):
  Assuming that every circuit has an entry, we will show local contractivity of $\SGE$ via  \cite[Proposition 6.7]{EPFour}.
  Given a nonzero idempotent  $e$ in $\EGE$, write $e=(\mu ,1,\mu )$ for some finite path $\mu $.  Using that $E$ has no sources,
we may find an  infinite path $\xi =\xi _1\xi _2\xi _3\ldots $, such that
  $\sr \mu  = \rn \xi $.
  Since $E$ is a finite graph, there must be  repetitions amongst the  $\xi _i$, say $\xi _i=\xi _j$, for
some  $i<j$.  Letting
  $$
  \alpha =\xi _1\xi _2\ldots \xi _i
  \and
  \gamma =\xi _{i+1}\xi _{i+2}\ldots \xi _j,
  $$
  notice that
  $$
  \sr\gamma  = \sr {\xi _j} = \sr {\xi _i} = \rn {\xi _{i+1}} =\rn \gamma ,
  $$
  so $\gamma $ is a circuit.  It is also clear that $\mu \alpha $ and $\mu \alpha \gamma $ are well defined paths.  Noticing that
  $$
  \sr \alpha  = \rn \gamma  = \sr \gamma ,
  $$
  we have that $s:=(\mu \alpha \gamma ,1,\mu \alpha )$ lies in $\SGE$.  Moreover,  setting
  $$
  \beta _1=\mu \alpha ,
  $$
  and using the notation introduced in \ref{DefineEAlpha}, we have
  $$
  sf_{\beta _1}s^* = (\mu \alpha \gamma ,1,\mu \alpha )(\mu \alpha ,1,\mu \alpha )  (\mu \alpha \gamma ,1,\mu \alpha )^* \={Conjuga} (\mu \alpha \gamma ,1,\mu \alpha \gamma ) \leq f_{\beta _1},
  $$
  thus  verifying \cite[Proposition 6.7.ii]{EPFour}.
  By hypothesis  $\gamma $  has an entry, so we may find a path $\gamma '$, with $\rn {\gamma '} = \rn \gamma $, which is not a prefix of $\gamma $, or
vice versa.  Setting
  $$
  \beta _0 = \mu \alpha \gamma ',
  $$
  we then have
  $$\def\quad{ }
  \matrix{ 0&\neq & f_{\mu \alpha \gamma '} &\leq & f_{\mu \alpha } &&&\leq & f_{\mu } & \ \ \Rightarrow \cr
  \pilar{18pt}
  0&\neq & f_{\beta _0} &\leq &  f_{\beta _1} &=&  s^*s &\leq &  e,}
  $$
  verifying \cite[Proposition 6.7.i]{EPFour}.
  Focusing now on \cite[Proposition 6.7.iii]{EPFour} notice that
  $$
  f_{\beta _0}s = (\mu \alpha \gamma ',1,\mu \alpha \gamma ')(\mu \alpha \gamma ,1,\mu \alpha ) = 0,
  $$
  precisely because $\gamma $ and $\gamma '$ are not each other's prefix.
  So evidently   $f_{\beta _0}sf_{\beta _1}=0$, proving the last condition in \cite[Proposition 6.7]{EPFour}, and hence that $\SGE$ is
locally contracting, thus proving (i).
\endProof

It is worth noticing that many  results of \cite{EPFour} used in the above proof,
  such as
  \cite[Proposition 6.3]{EPFour},
  \cite[Theorem 6.5]{EPFour} and
  \cite[Proposition 6.7]{EPFour},
  comparing local contractivity for groupoids, inverse semigroups, and actions,
  are either one way implications only, or the converse depends on special conditions.  Nevertheless, the situation in
which we are working has fortunately allowed for a downright equivalence of the various manifestations of contractivity.

However, this result  should be taken with a certain skepticism. First of all it is well known that the above condition on circuits
is not sufficient for local contractivity for the groupoid associated to \"{infinite} graphs \cite{KPR}.  Considering that
\"{finite} graphs are special cases of our theory (just take the acting group to be the trivial group), it is not
unreasonable to believe that our results admit natural  generalizations to infinite graphs,
  but then a characterization of local contractivity for the corresponding groupoid will certainly not follow from the
fact that every circuit has an entry, since this is false for infinite graphs, as mentioned above.

  Secondly, observe that the condition on the existence of entries for circuits completely ignores the group $G$, but, again, a
generalization to infinite graphs will probably depend on the action.  A hypothesis such as ``every vertex connects
to a $G$-circuit with an entry'', to paraphrase the main hypothesis of \cite[Lemma 3.8]{KPR}, is probably more realistic in the
conjectured infinite graph scenario.

\section {Simplicity and pure infiniteness for $\OGE $}

In this section we use the results in the previous sections to characterize when $\OGE $ is simple and purely
infinite. The central results are the following:

\state Theorem \label ThmSimple
  Assume that $(G,E, \varphi )$ satisfies \ref{StandingHyp}, that $G$ is amenable, and that for every $g\in G$ there are
at most finitely many minimal strongly fixed paths for $g$.   Then  $\OGE $ is simple if and only if the following
conditions are satisfied:
  \iaitem
  \aitem $E$ is weakly-$G$-transitive.
  \aitem Every $G$-circuit has an entry.
  \aitem Given a vertex $x$, and a group element $g$ fixing $Z(x)$ pointwise, then necessarily $g$ is a slack at $x$.

  \Proof
  By \ref{MainHausdorff}, the groupoid $\GpdGE $ is Hausdorff. Clearly $\GpdGE $ is \'etale with second countable unit
space. By \ref{Newamenabuilidade}, $\GpdGE $ is amenable. Then, by \ref{UniversalTightAlgebra}, $\OGE \cong C^*(\GpdGE)
= C_{r}^*(\GpdGE)$. By \cite [Theorem 5.1]{SimpleGroupoid}, $\OGE $ is simple if and only if $\GpdGE $ is minimal and
essentially principal. Since minimality of $\GpdGE $ is equivalent to (a) by \ref{CharacMinimal}, and essential
principality of $\GpdGE $ is equivalent to (b\&c) by \ref{TopFreeSameEssPrin}, the result holds.
  \endProof

With respect to pure infiniteness, we have:

\state Theorem \label Thmpi
Let $(G,E, \varphi )$ be under \ref{StandingHyp}, and let $G$ be an amenable group. If $\GpdGE $ is essentially principal,  then every hereditary
subalgebra of $\OGE $ contains an infinite projection.

  \Proof
  By the same argument as in \ref{ThmSimple}, $\OGE \cong C^*(\GpdGE) = C_{r}^*(\GpdGE)$. By \ref{MainTopFreeSGE.i}, every circuit of $E$ has an entry. Thus, $\GpdGE $ is locally contracting by \ref{LocContrVsTopFree}. Hence, by \cite [Proposition 2.4]{AdelaR}, every nonzero hereditary sub-C*-algebra of $\OGE$ contains an infinite projection, as desired.
  \endProof

As an immediate consequence we have

\state Corollary \label CorolPinfSimple
  If $(G,E, \varphi )$ satisfies \ref{StandingHyp}, the group $G$ is amenable, and
$\GpdGE $ is Hausdorff then, whenever  $\OGE $ is simple,  it is necessarily also  purely infinite (simple).

  \Proof
  By \ref{MainTopFreeSGE} and \ref{ThmSimple}, $\GpdGE $ is essentially principal. Thus, by \ref{Thmpi}, every nonzero hereditary sub-C*-algebra of $\OGE$ contains an infinite projection. Hence, $\OGE $ is purely infinite simple, as desired.
  \endProof

\section {Revisiting Nekrashevych algebras}

In this section we will analyze Nekrashevych algebras from our point of view.

The Nekrashevych C*-algebra $\mathcal {O}_{(G,X)}$, associated to a self-similar action of a group $G$ on a finite
alphabet $X$ \cite {NC}, is a direct example of our definition (see \ref{Nekrashevych}). Here, the graph $E$ is the rose of $n$ petals
for $n=\vert X\vert \geq 2$, so that the action on vertices is trivial, and the action is faithful.  Since $\vert
E^0\vert =1$, we have the following facts:
  \initem
  \nitem $E$ is $G$-transitive, whence $\GpdGE $ is minimal by \ref{CharacMinimal}.
  \nitem $\GpdGE $ is essentially principal by \ref{TopFreeUnderSSGroup} and \cite[Theorem 14.7]{EPFour}. In particular, $\GpdGE $ is locally contracting by \ref{MainTopFreeSGE} and \ref{LocContrVsTopFree}.

\medskip

\noindent Thus, if $\GpdGE $ is Hausdorff, we conclude that $C_r^*(\GpdGE)$ is a purely infinite simple C*-algebra by \cite[Theorem 6.8]{EPFour}. Hausdorffness of $\GpdGE $ is equivalent, according to \ref{MainHausdorff}, to the existence of at most finitely many minimal strongly fixed paths for every $g\in G$.

In this sense, it is interesting to remark that Nekrashevych also gave a presentation of its algebra as a groupoid
C*-algebra associated to a groupoid of germs of an inverse semigroup $\mathcal S$ \cite[Section 5]{NC}. While $\mathcal
S$ turns out to be $\SGE $, the notion of germ that he used is the one adopted by Arzumanian and Renault \cite{ArzRen}, which differs from the one we used, due to Patterson \cite[Page 140]{pat}. Luckily, both definitions coincide when the action of $\SGE$ on $E^{\infty}$ is topologically free, which is the case of Nekrashevych triples, as we noticed above. So, Nekrashevych's groupoid and $\GpdGE $ coincide, and the characterization of Hausdorffness we obtained in \ref{MainHausdorff} coincide with that given by Nekrashevych \cite[Lemma 5.4]{NC}.

In order to obtain a characterization of (pure infinite) simplicity for $\mathcal {O}_{(G,X)}$, we need to keep control of whether $\mathcal {O}_{(G,X)}$ is nuclear. So, it only remains to determine when $\GpdGE $ is amenable, which implies that $C^*(\GpdGE )=C_r^*(\GpdGE )$. By \ref{Newamenabuilidade}, if $G$ is an amenable group, then $\GpdGE $ is an amenable groupoid. Thus, we obtain

\state Proposition \label NekraOne
If $(G,X,\varphi)$ is a Nekrashevych triple, with $G$ an  amenable group and $\GpdGE$ a  Hausdorff groupoid, then $\mathcal {O}_{(G,X)}$ is a nuclear, separable, purely infinite simple C*-algebra.

Here, Nekrashevych's approach differs from ours. In \cite{NC} he stated a sufficient condition for the amenability of
$\GpdGE$, which apparently does not require the group $G$ to be  amenable. The condition relies on two concepts associated
to self-similar groups: self-replication and contractiveness (see \cite{Nmsn} or \cite{NC} for definitions of these
concepts). Nekrashevych \cite[Theorem 5.6]{NC} proved that if $\GpdGE$ is Hausdorff and $(G,X)$ is self-replicating and
contractive, then $\GpdGE$ is of polynomial growth \cite{Nmsn}, and thus it its amenable by \cite[Proposition
3.2.32]{AnanRen}.

\section {Revisiting Katsura algebras}

\label SectKatAlg In this section we will analyze Katsura algebras from our point of view.

We will quickly recall the definition and basic properties of Katsura algebras that will be needed in the sequel. This
is borrowed from \cite {KatsuraOne}.

\definition \label {DefKatAlgData}
  Let $N\in \Ninf $, let $A\in M_N(\Zplus )$ and $B\in M_N(\Z )$ be row-finite matrices. Define a set $\OmA $ by $$\OmA
:=\{ (i,j)\in \{ 1, 2, \dots ,N\}\times \{ 1, 2, \dots ,N\} : A_{i,j}\geq 1 \}.$$ For each $i\in \{ 1, 2, \dots ,N\}$,
define a set $\OmA (i)\subset \{ 1, 2, \dots ,N\}$ by $$\OmA (i):=\{ j \in \{ 1, 2, \dots ,N\}: (i,j)\in \OmA \}.$$
Notice that, by definition, $\OmA (i)$ is finite for all $i$.  Finally, fix the following relation:
  \medskip \quad (0) $\OmA (i)\ne \emptyset $ for all $i$, and $B_{i,j}=0$ for $(i,j)\not \in \OmA $.

With these data we can define Katsura algebras

\definition \label {DefKatAlgAlgebra}
Define $\OAB $ to be the universal C*-algebra generated by
mutually orthogonal projections $\{ q_i\}_{i=1}^N$, partial unitaries $\{ u_i\}_{i=1}^N$ with $u_iu_i^*=u_i^*u_i=q_i$,
and partial isometries $\{ s_{i,j,n}\}_{(i,j)\in \OmA , n\in \Z }$ satisfying the relations:
  \izitem
  \zitem $s_{i,j,n}u_j=s_{i,j, n+A_{i,j}}$ and $u_is_{i,j,n}=s_{i,j, n+B_{i,j}}$ for all $(i,j)\in \OmA $ and $n\in \Z
$.
  \zitem $s_{i,j,n}^*s_{i,j,n}=q_j$ for all $(i,j)\in \OmA $ and $n\in \Z $.
  \zitem $q_i=\sum \limits _{j\in \OmA (i)}\sum \limits _{n=1}^{A_{i,j}}s_{i,j,n}s_{i,j,n}^*$ for all $i$.

\state Remark  \label KatsuraConditions
  \rm Now, the following facts holds: \uplevel \itm The C*-algebra $\OAB $ is separable,
nuclear, in the UCT class \cite [Proposition 2.9]{KatsuraOne}.  \itm If the matrices $A,B$ satisfy the following
additional properties:
  \uplevel
  \itm $A$ is irreducible, and
  \itm $\vert A_{ii}\vert\geq 2$ and $B_{i,i}=1$ for every $1\leq i\leq N$,
  \dnlevel
  \medskip
  then the C*-algebra $\OAB $ is purely infinite simple, and hence a Kirchberg algebra \cite [Proposition
2.10]{KatsuraOne}.
  \itm The $K$-groups of $\OAB $ are \cite [Proposition 2.6]{KatsuraOne}:
    \uplevel
    \itm $K_0(\OAB )\cong \mbox {coker}(I-A)\oplus \mbox {ker}(I-B)$, and
    \itm $K_1(\OAB )\cong \mbox {coker}(I-B)\oplus \mbox {ker}(I-A)$.
    \dnlevel
  \itm Every Kirchberg algebra can be represented, up to isomorphism, by an algebra $\OAB $ for matrices $A,B$
satisfying the conditions in \ref{KatsuraConditions.2} \cite [Proposition 4.5]{KatsuraTwo}.
  \dnlevel

As we have seen in \ref{KatsuraExample}, unital Katsura algebras are natural examples of our construction. So, it is easy to use our results in order to characterize some properties, like  simplicity or pure infinite simplicity, in terms of matrices $A$ and $B$. This work has been previously done in \cite{EP}, but the approach we chose there was fairly more direct and computational, so that the conditions appearing there were less elegant and clear than the ones we will present here.

Across this section, we will say that a triple $(\Z, E, \varphi)$ is a Katsura triple if there exist finite matrices $A,B$ satisfying \ref{DefKatAlgData} such that the triple associated to the algebra $\OAB$ is $(\Z, E, \varphi)$; in particular, $E$ is the graph whose adjacency matrix is $A$. Also, we will fix the following agreement: let $\xi$ be either in $E^*$ or in $E^{\infty}$, i.e.
$$\xi= e_{i_1,i_2,n_1}e_{i_2,i_3,n_2}\cdots e_{i_k,i_{k+1},n_k} \mbox{ or }\xi=e_{i_1,i_2,n_1}e_{i_2,i_3,n_2}\cdots e_{i_k,i_{k+1},n_k}\cdots ,$$
then, for any $r\in \N$ we define
$$B_{{\xi}_{\vert r}}:= \prod \limits _{t=1}^{r}{B_{i_t,i_{t+1}}}\text{ and }A_{{\xi}_{\vert r}}:= \prod \limits _{t=1}^{r}{A_{i_t,i_{t+1}}}.$$

The first step to work out the corresponding results to the ones we obtained for the general setting is to  determine when a finite path $\alpha\in E^*$ is fixed under the action of an element $l\in \Z$.

\state Lemma \label LemFixedUSubI
Let $(\Z ,E, \varphi )$ be a Katsura triple.  Given an element $\alpha$ of $E^*$ of length $r$ and an
integer $l\in \Z $, the following are equivalent:
  \initem
  \nitem $\alpha $ is fixed under the action of $l$.
  \nitem For every $1\leq j\leq r$ the element $K_j:=l \displaystyle \frac{B_{{\alpha}_{\vert j}}}{A_{{\alpha}_{\vert j}}}$ belongs to $\Z $.

  \Proof
Set $\alpha =e_{i_1,i_2,n_1}e_{i_2,i_3,n_2}\cdots e_{i_r,i_{r+1},n_r}$.  By definition of $(\Z,E,\varphi) $, $\alpha =l\alpha $ if and only if there exists a sequence $(K_j)_{j\geq 0}\subseteq \Z
$ such that: \initem
  \izitem
  \zitem $K_0=l$.
  \zitem For every $1\leq j\leq r$, $n_{j-1}+K_{j-1}B_{i_j,i_{j+1}}=n_{j-1}+K_{j}A_{i_j,i_{j+1}}$.  \medskip \noindent Notice
that (ii) is equivalent to ask $K_{j-1}B_{i_j,i_{j+1}}=K_{j}A_{i_j,i_{j+1}}$ for every $j\geq 1$.

Now, for $j=1$ we have $K_0 B_{i_1,i_{2}}=l B_{i_1,i_{2}}=K_{1}A_{i_1,i_{2}}$, so that $K_1=l\displaystyle \frac
{B_{i_1,i_{2}}}{A_{i_1,i_{2}}}$. Now, suppose that for $1\leq t\leq j-1$ we have proved that $K_t:=l \displaystyle \frac{B_{{\alpha}_{\vert t}}}{A_{{\alpha}_{\vert t}}}$. Hence
$$K_{j}A_{i_j,i_{j+1}}=K_{j-1}B_{i_j,i_{j+1}}=l \displaystyle \frac{B_{{\alpha}_{\vert j-1}}}{A_{{\alpha}_{\vert j-1}}}\cdot B_{i_j,i_{j+1}},$$ so that $K_j=l \displaystyle \frac{B_{{\alpha}_{\vert j}}}{A_{{\alpha}_{\vert j}}}$. This completes the proof.
  \endProof

Now, we are ready to characterize pseudo freeness for a Katsura triple $(\Z ,E, \varphi )$.

\state Lemma \label KatsuPF
Let $(\Z , E, \varphi )$ be a Katsura triple. Then,
the following are equivalent:
\initem
\nitem $(\Z ,E, \varphi )$  is pseudo free.
\nitem $B_{i,j}=0$ if and only if $(i,j)\not \in \OmA $.

\Proof
Let $\alpha =e_{i_1,i_2,n_1}e_{i_2,i_3,n_2}\cdots e_{i_r,i_{r+1},n_r}$ of $E^*$, and let $l\in Z$. By \ref{LemFixedUSubI}, $l\alpha=\alpha$ exactly when the elements $K_j:=l \displaystyle \frac{B_{{\alpha}_{\vert j}}}{A_{{\alpha}_{\vert j}}}$ belongs to $\Z $ for every $1\leq j\leq r$. Since $\varphi (\l,\alpha_{\vert j})=K_j$, $\varphi (\l,\alpha)=0$ exactly when $K_j=0$ for some $j\leq r$. Thus, the situation reduces to $K_{j-1}e_{i_j,i_{j+1}, n_l}=e_{i_j,i_{j+1}, n_l}$ and $K_j=0$ for some $1\leq j\leq r$, which corresponds to the equation $n_j+K_{j-1}B_{i_j,i_{j+1}}=n_j$. And this occurs exactly when $B_{i_j,i_{j+1}}=0$, so we are done.
\endProof

Which these results in mind, we are ready to characterize when $\GZE $ is Hausdorff

\state Theorem \label KatsuHausdorff
Let $(\Z , E, \varphi )$ be a  Katsura triple. Then,
the following are equivalent: \initem \nitem $\GZE $ is Hausdorff.  \nitem Whenever $(i,j)\in \Omega _A$ with $B_{i,j}=0$,
then for any $l\in \Z $ there exist finitely many finite paths $\alpha \in E^*$ with $d(\alpha)=i$ such that $l \displaystyle \frac{B_{{\alpha}_{\vert t}}}{A_{{\alpha}_{\vert t}}}\in \Z $ for every $1\leq t\leq r-1$.

  \Proof
  The result holds by \ref{MainHausdorff} and \ref{KatsuPF}.
  \endProof

The next step is to determine the minimality of  $\GZE $.

\state Theorem \label KatsuraMinimal
Let $(\Z ,E, \varphi )$ be a Katsura triple. Then, the following are equivalent:
\initem
\nitem $\GZE $ is minimal.
\nitem The adjacency matrix $A$ of $E$ is irreducible.

   \Proof
First notice that $E$ has no sinks by \ref{DefKatAlgData.(0)}. Moreover, the action of $\Z$ on $E$ fixes all the vertices. Then, by \ref{CharacMinimal}, $\GZE $ is minimal if and only if $E$ is transitive, which is equivalent to the matrix $A$ being irreducible, so we are done.
   \endProof

Now, we will give a characterization of when $\GZE $ is essentially principal.

\state Theorem \label KatsuraEssPrinc
Let $(\Z ,E, \varphi )$ be a Katsura triple. Then, the following are equivalent:
\uplevel \itm $\GZE $ is essentially principal.  \itm \notext \uplevel
  \itm Every circuit in $E$ has an entry.
  \itm If $1\leq i\leq N$, $l\in \Z$, and for any $\xi\in Z(i)$ the elements $l\displaystyle \frac{B_{{\xi}_{\vert n}}}{A_{{\xi}_{\vert n}}}\in \Z$ for all $n\in \N$, then there exists $m\in \N$ such that $B_{{\xi}_{\vert m}}=0$ for all $\xi\in Z(i)$.
  \dnlevel \dnlevel

   \Proof
Since the action of $\Z$ fixes all the vertices of $E$, (2a) is \ref{MainTopFreeSGE.i}. On the other side, (2b) is exactly \ref{MainTopFreeSGE.ii} because of \ref{LemFixedUSubI} and \ref{KatsuPF}. Thus, the result is consequence of \ref{TopFreeSameEssPrin}.
   \endProof

We can obtain an easy sufficient condition for $\GZE $ being essentially principal.

\state Corollary \label CorolKatsuraEssPrinc
Let $(\Z ,E, \varphi )$ be a Katsura triple. If
\initem
\nitem Every circuit of $E$ has an entry, and
\nitem For every $1\leq i\leq N$ and every $l\in \Z$ there exists $\eta \in Z(i)$ such that $\lim\limits_{n\rightarrow \infty} l\displaystyle\frac{B_{\eta_{\vert n}}}{A_{\eta_{\vert n}}} =0$,
\medskip
\noindent then $\GZE $ is essentially principal.

   \Proof
   By \ref{LemFixedUSubI}, condition (2) implies that $l\eta\ne \eta$ for any $l\in \Z$, whence the triple $(Z,E,\varphi)$ trivially satisfies \ref{MainTopFreeSGE.ii}, as remarked in \ref{RemarkEssPrincTriv}.
   \endProof

  Corollary \ref{CorolKatsuraEssPrinc} applies when we have a pair of finite matrices $A,B$ under \ref{DefKatAlgData}, such that for every $1\leq i\leq N$ we have $\vert A_{ii}\vert \geq 2$ and $B_{ii}< \vert A_{ii}\vert$. In particular, $\GZE$ is essentialy principal for Katsura systems $(\Z,E,\varphi)$ satisfying \ref{KatsuraConditions.2}.

Also, it is immediate to characterize when $\GZE$ is locally contracting.

\state Theorem \label KatsuraLocContract
Let $(\Z ,E, \varphi )$ be a Katsura triple. Then, the following are equivalent:
\initem
\nitem $\GZE $ is locally contracting.
\nitem Every circuit of $E$ has an entry.

   \Proof
This is \ref{LocContrVsTopFree}.
   \endProof

Finally, we have the following fact

\state Proposition \label KatsuraAmenable
If $(\Z, E, \varphi)$ is a Katsura triple, then $\GZE$ is an amenable groupoid.

   \Proof
Since $\Z$ is an amenable group, \ref{Newamenabuilidade} applies.
   \endProof

Now, we are ready to characterize simplicity of the algebra $\OAB$, as follows

\state Theorem \label ThmSimpleOld
Let $(\Z ,E,\varphi )$ be a Katsura triple such that $\GZE$ is Hausdorff (see \ref{KatsuHausdorff}). Then, the
following are equivalent: \uplevel \itm \notext
  \uplevel
  \itm The matrix $A$ is irreducible.
  \itm Every circuit of $E$ has an entry.
  \itm If $1\leq i\leq N$, $l\in \Z$, and for any $\xi\in Z(i)$ the elements $l\displaystyle \frac{B_{{\xi}_{\vert n}}}{A_{{\xi}_{\vert n}}}\in \Z$ for all $n\in \N$, then there exists $m\in \N$ such that $B_{{\xi}_{\vert m}}=0$.
   \dnlevel \itm $\OAB $ is simple.
  \dnlevel

   \Proof
   This is exactly \ref{ThmSimple} for the Katsura triple $(\Z ,E,\varphi )$, because of \ref{KatsuraMinimal}, \ref{KatsuraEssPrinc} and \ref{KatsuraAmenable}.
   \endProof

In particular, when $\GZE$ is Hausdorff and $\OAB$ is simple, the $\GZE$ is locally contracting by \ref{KatsuraLocContract} and \ref{ThmSimpleOld.1b}. Hence, we have

\state Corollary \label KatsuraPurInfSimple
If $(\Z ,E,\varphi )$ is a Katsura triple such that $\GZE$ is Hausdorff and $\OAB$ is simple, then $\OAB$ is purely infinite simple.

   \Proof
   This is by \ref{CorolPinfSimple}.
   \endProof

\state Remark
\rm Notice that, because of \ref{CorolKatsuraEssPrinc}, Katsura's condition \ref{KatsuraConditions.2} for $\OAB$ being a purely infinite simple C*-algebra derive directly from \ref{ThmSimpleOld} and \ref{KatsuraPurInfSimple} when $\GZE$ is Hausdorff. Moreover, when $\GZE$ is Hausdorff, \ref{ThmSimpleOld} provides a characterization of simplicity for $\OAB$, improving Katsura's results on that direction, where only sufficient conditions are given \cite{KatsuraOne}.

We close this section by presenting a couple of examples. The first one illustrates the difference between $(G,E,
\varphi )$ being pseudo free and $\GpdGE $ being Hausdorff, and also the difference between the action of $\SGE $ on
$E^{\infty }$ being topologically free and the action of $G$ on $E^{\infty }$ being topologically free.

\state Example \rm \label KatsuNoPSbutHaus
Set $N=2$, and consider the matrices $A=\left (\matrix { 2 & 1 \cr 1 & 2
}\right )$ and $B=\left (\matrix { 1 & 0 \cr 0 & 1 }\right )$. Let $(\Z , E, \varphi )$ be the associated Katsura triple. Then, we have the following: \initem \nitem Since $\Z $ is amenable, then so is $\GZE $.  \nitem
Since $A$ is irreducible, $\GZE $ is minimal.  \nitem Every circuit in $E$ has an entry.  \nitem Since $A_{1,2}\ne 0$ and $B_{1,2}=0$, $(\Z ,E, \varphi )$ is not pseudo free by
  \ref{KatsuPF}.  \nitem Notice that the only possible quotient values $\displaystyle \frac
{B_{i,j}}{A_{i,j}}$ are $\displaystyle \frac {B_{1,1}}{A_{1,1}}=\displaystyle \frac {B_{2,2}}{A_{2,2}}=\displaystyle
\frac {1}{2}$ and $\displaystyle \frac {B_{1,2}}{A_{1,2}}=\displaystyle \frac {B_{2,1}}{A_{2,1}}=\displaystyle \frac
{0}{1}=0$. Then, for any $l\in \Z $, it is clear that there exists only finitely many minimal strongly fixed paths for $l$. Thus, $\GZE $ is Hausdorff by
  \ref{KatsuHausdorff}.  \nitem Moreover, by the argument in point
$(5)$, the only infinite paths fixed by the action of $\Z $ are the ones associated to minimal strongly fixed paths, and thus trivial. Hence, the
action of $\SZE $ on $E^{\infty }$ is topologically free. But the action of $\Z $ is not topologically free, since every
element of $\Z $ fix the cylinders $Z(e_{1,2,1})$ and $Z(e_{2,1,1})$.

\noindent Notice that $\OAB $ is purely infinite simple by \ref{ThmSimpleOld} and \ref{KatsuraPurInfSimple}.\vspace{.2truecm}

The second example shows that $(G,E,\varphi )$ being pseudo free do not imply that the action of
$\SGE $ is topologically free.

\state Example \rm \label {KatsuPSbutNoTopFree} Set $N=1$, and set $A=B=(n)$ for any $n\geq 2$. Let $(\Z , E, \varphi )$ be
the associated Katsura triple. Then, we have the following: \initem \nitem Since $\Z $ is
amenable, then so is $\GZE $.  \nitem Since $A$ is irreducible, $\GZE $ is minimal.
\nitem Every circuit in $E$ has an entry.  \nitem Since $A=B$, $(\Z , E, \varphi )$ is pseudo free, whence in
particular $\GZE $ is Hausdorff.  \nitem Since $A=B$, the action of $\Z $ on $E$ is trivial. Thus, the action of $\SZE $ on
$E^{\infty }$ (and of $\Z$) cannot be topologically free, because there are no slacks.

\references

\Article AdelaR
  C. Anantharaman-Delaroche;
  Purely infinite $C^*$-algebras arising form dynamical systems;
  Bull. Soc. Math. France, 125 (1997), no. 2, 199-225

\Bibitem AnanRen
  C. Anantharaman-Delaroche and J. Renault;
  Amenable groupoids;
  Monogr. Enseign. Math. 36, Universit\'e de Gen\`eve, 2000

\Bibitem ArzRen
  V. Arzumanian, J. Renault;
  Examples of pseudogroups and their $C^*$-algebras;
  Operator algebras and quantum field theory (Rome, 1996), 93-104, Int. Press, Cambridge, MA, 1997

\Article SimpleGroupoid
  J. Brown, L. O. Clark, C. Farthing and A. Sims;
  Simplicity of algebras associated to \'etale groupoids;
  Semigroup Forum, 88 (2014), 433-452

\Bibitem BO
  N. P. Brown and N. Ozawa;
  C*-algebras and finite-dimensional approximations;
  Graduate Studies in Mathematics, 88, American Mathematical Society, 2008

\Article actions
  R. Exel;
  Inverse semigroups and combinatorial C*-algebras;
  Bull. Braz. Math. Soc., 39 (2008), no. 2, 191-313

\Article nhausdorff
    R. Exel;
    Non-Hausdorff \'etale groupoids;
    Proc. Amer. Math. Soc., 139 (2011), no. 3, 897-907

\Bibitem book
  R. Exel;
  Partial Dynamical Systems, Fell Bundles and Applications;
  {Licensed under a Creative Commons Attribution-ShareAlike 4.0 International License, 351pp, 2014. \hfill  Available online from\break
mtm.ufsc.br/$\sim$exel/publications. \rm PDF file md5sum: bc4cbce3debdb584ca226176b9b76924}

\Article ExelLaca
  R. Exel and M. Laca;
  Cuntz-Krieger algebras for infinite matrices;
  J. reine angew. Math., 512 (1999), 119-172

\Bibitem EP
  R. Exel and E. Pardo;
  Representing Kirchberg algebras as inverse semigroup crossed products;
  arXiv:1303.6268 [math.OA], 2013

\Bibitem EPTwo
  R. Exel and E. Pardo;
  Graphs, groups and self-similarity;
  arXiv:1307.1120 [math.OA], 2013

\Bibitem EPFour
  R. Exel and E. Pardo;
  The tight groupoid of an inverse semigroup;
  arXiv:1408.5278 [math.OA], 2014

\Bibitem ExelStar
  R. Exel and C. Starling;
  Self-similar graph C*-algebras and partial crossed products;
  arXiv:1406. 1086 [math.OA], 2014

\Article ExelVesshik
  R. Exel and A. Vershik;
  C*-algebras of irreversible dynamical systems;
  Canadian Mathematical Journal, 58 (2006), 39-63

\Article Grig
  R. I. Grigorchuk;
  On Burnside's problem on periodic groups;
  Funct. Anal. Appl., 14 (1980), 41-43

\Article GS
  N. D. Gupta and S. N. Sidki;
  On the Burnside problem for periodic groups;
  Math. Z., 182 (1983), 385-388

\Article KatsuraFundRes
  T. Katsura;
  A class of C*-algebras generalizing both graph algebras and homeomorphism C*-algebras. I. Fundamental results;
  Trans. Amer. Math. Soc., 356 (2004), no. 11, 4287-4322

\Article KatsuraOne
  T. Katsura;
  A construction of actions on Kirchberg algebras which induce given actions on their $K$-groups;
  J. reine angew. Math., 617 (2008), 27-65

\Article KatsuraTwo
  T. Katsura;
  A class of $C^*$-algebras generalizing both graph algebras and homeomorphism $C^*$-algebras IV, pure infiniteness;
  J. Funct. Anal., 254 (2008), 1161-1187

\Article KPRR
  A. Kumjian, D. Pask, I. Raeburn and J.  Renault;
  Graphs, groupoids, and Cuntz-Krieger algebras;
  J. Funct. Anal., 144 (1997), 505-541

\Article KPR
  A. Kumjian, D. Pask and I. Raeburn;
  Cuntz-Krieger algebras of directed graphs;
  Pacific J. Math., 184 (1998), no. 1, 161-174

\Bibitem Lawson
  M. V. Lawson;
  Inverse semigroups, the theory of partial symmetries;
  World Scientific, 1998

\Article LawsonCompactable
  M. V. Lawson;
  Compactable semilattices;
  Semigroup Forum, 81 (2010), no. 1, 187-199

\Article NekraJO
  V. Nekrashevych;
  Cuntz-Pimsner algebras of group actions;
  J. Operator Theory, 52 (2004), 223-249

\Bibitem Nmsn
  V. Nekrashevych;
  Self-similar groups;
  Mathematical Surveys and Monographs, \bf 117, \rm  Amer. Math. Soc., Providence, RI, 2005

\Article NC
  V. Nekrashevych;
  C*-algebras and self-similar groups;
  J. reine angew. Math., 630 (2009), 59-123

\Bibitem pat
  A. L. T. Paterson;
  Groupoids, inverse semigroups, and their operator algebras;
  Birkh\umlaut auser, 1999

\Bibitem Pedersen
  G. K. Pedersen;
  C*-algebras and Their Automorphism Groups;
  Academic Press, 1979

\Article Pimsner
  M. V. Pimsner;
  A class of C*-algebras generalizing both Cuntz-Krieger algebras and crossed products by ${\bf Z}$;
  Fields Inst. Commun., 12 (1997), 189-212

\Bibitem Raeburn
  I. Raeburn;
  Graph algebras;
  CBMS Regional Conference Series in Mathematics, \bf 103 \rm (2005), pp. vi+113

\Article RenaultCartan
  J. Renault;
  Cartan subalgebras in $C^*$-algebras;
  Irish Math. Soc. Bull., 61 (2008), 29-63

\Bibitem Starling
  C. Starling;
  Boundary quotients of C*-algebras of right LCM semigroups;
  arXiv:1409.1549\break\hfill [math.OA]

\Article Steinberg
  B. Steinberg;
  A groupoid approach to discrete inverse semigroup algebras;
  Adv. Math., 223 (2010), 689-727

\endgroup

\bigskip \noindent \tensc
  Departamento de Matem\'atica;
  Universidade Federal de Santa Catarina;
  88010-970 Florian\'opolis SC;
  Brazil
  \hfill \break \tt (ruyexel@gmail.com)

\bigskip \noindent \tensc
  Departamento de Matem\'aticas, Facultad de Ciencias;
  Universidad de C\'adiz, Campus de Puerto Real;
  11510 Puerto Real (C\'adiz);
  Spain
  \hfill \break \tt (enrique.pardo@uca.es)

  \bye

\bye